\documentclass[reqno, a4paper]{amsart}

\usepackage[margin=1in]{geometry}
\usepackage[algo2e, linesnumbered, noline, noend, ruled]{algorithm2e}
\usepackage{amsfonts}
\usepackage{amsmath}
\usepackage{amsthm}
\usepackage[foot]{amsaddr}
\usepackage{amssymb}
\usepackage[english]{babel}
\usepackage{bbm}
\usepackage{booktabs}
\usepackage{braket}
\usepackage{float}
\usepackage[T1]{fontenc}
\usepackage{graphicx}
\usepackage[utf8]{inputenc}
\usepackage{lmodern}
\usepackage{nicefrac}
\usepackage[retain-explicit-plus, retain-zero-exponent, per-mode=symbol, output-exponent-marker=\text{e}]{siunitx}
\usepackage{subcaption}
\usepackage[table]{xcolor}
\usepackage{tikz}
\usepackage{cite}
\usepackage{hyperref}

\usetikzlibrary{positioning, arrows, backgrounds}
\definecolor{tabblue}{RGB}{0, 107, 164}
\definecolor{taborange}{RGB}{255, 128, 14}
\definecolor{tablightgray}{RGB}{171,171,171}
\definecolor{tabdarkgray}{RGB}{89,89,89}
\definecolor{tablightblue}{RGB}{95,158,209}

\newtheorem{Theorem}{Theorem}[section]

\newtheorem{Definition}[Theorem]{Definition}

\newtheorem*{Remark*}{Remark}

\newcommand{\qe}[1]{``#1''}
\newcommand{\tr}[1]{\mathrm{tr} #1 }

\newcommand{\R}{\mathbb{R}}

\newcommand{\norm}[2]{\left\lvert \left\lvert #1 \right\rvert \right\rvert_{#2}}
\newcommand{\abs}[1]{\left\lvert #1 \right\rvert}
\newcommand{\dx}[1]{\ \mathrm{d}#1}

\makeatletter
\DeclareOldFontCommand{\rm}{\normalfont\rmfamily}{\mathrm}
\DeclareOldFontCommand{\sf}{\normalfont\sffamily}{\mathsf}
\DeclareOldFontCommand{\tt}{\normalfont\ttfamily}{\mathtt}
\DeclareOldFontCommand{\bf}{\normalfont\bfseries}{\mathbf}
\DeclareOldFontCommand{\it}{\normalfont\itshape}{\mathit}
\DeclareOldFontCommand{\sl}{\normalfont\slshape}{\@nomath\sl}
\DeclareOldFontCommand{\sc}{\normalfont\scshape}{\@nomath\sc}
\makeatother

\numberwithin{equation}{section}

\newcommand{\tablegray}{gray!25}

\newcommand{\subin}{\mathrm{in}}
\newcommand{\subwall}{\mathrm{wall}}
\newcommand{\subout}{\mathrm{out}}
\newcommand{\subfluid}{\mathrm{f}}
\newcommand{\subdarcy}{\mathrm{d}}
\newcommand{\subchannels}{\mathrm{mc}}
\newcommand{\subdes}{\mathrm{des}}
\newcommand{\subfd}{\mathrm{fd}}

\newcommand{\dimred}[1]{\tilde{#1}}
\newcommand{\porous}[1]{#1^{\mathrm{por}}}

\newcommand{\viscosity}{\mu}
\newcommand{\density}{\rho}
\newcommand{\hcapacity}{C_p}
\newcommand{\conductivity}{\kappa}
\newcommand{\htc}{\alpha}
\newcommand{\height}{h}
\newcommand{\massin}{m_\mathrm{in}}
\newcommand{\hfs}{h_\mathrm{fs}}
\newcommand{\volumefrac}{\varphi}
\newcommand{\normal}{n}
\newcommand{\permeability}{K}

\newcommand{\velocity}{u}
\newcommand{\pressure}{p}
\newcommand{\temperature}{T}
\newcommand{\solution}{U}
\newcommand{\testvelo}{\hat{v}}
\newcommand{\testpres}{\hat{q}}
\newcommand{\testtemp}{\hat{S}}

\newcommand{\testadvelo}{\hat{u}}
\newcommand{\testadpres}{\hat{p}}
\newcommand{\testadtemp}{\hat{T}}

\newcommand{\testadsolution}{\hat{U}}
\newcommand{\testsolution}{\hat{V}}
\newcommand{\advelo}{v}
\newcommand{\adpres}{q}
\newcommand{\adtemp}{S}
\newcommand{\adsolution}{P}

\newcommand{\fspacetrial}{V}
\newcommand{\fspacetest}{V_0}
\newcommand{\pspace}{P}
\newcommand{\tspacetrial}{W}
\newcommand{\tspacetest}{W_0}

\newcommand{\statespace}{\mathcal{U}}
\newcommand{\adjointspace}{\mathcal{P}}

\newcommand{\laplace}{\Delta}
\newcommand{\grad}{\nabla}
\newcommand{\divergence}[1]{\mathrm{div}\left(#1\right)}
\newcommand{\integral}[1]{\int_{#1}}

\newcommand{\tdiv}[1]{\mathrm{div}_\Gamma \left( #1 \right)}

\newcommand{\shapelagrangian}{\mathcal{G}}

\newcommand{\costfunction}{J}
\newcommand{\reducedcostfunction}{j}

\newcommand{\vectorfield}{\mathcal{V}}
\newcommand{\shapegradient}{\mathcal{G}}
\newcommand{\weighttemp}{\lambda_1}
\newcommand{\weightvelo}{\lambda_2}
\newcommand{\weightreg}{\lambda_3}
\newcommand{\flow}{\Phi_t}
\newcommand{\searchdirection}{\mathcal{V}^\mathrm{s}}

\newcommand{\fluxdes}{Q_\subdes}
\newcommand{\flux}{Q}

\begin{document}

{\footnotesize
	\begin{center}
		This is a preprint. The final version of this article can be found at \url{https://doi.org/10.1002/zamm.202000166}.
	\end{center}
}

\title{Model Hierarchy for the Shape Optimization of a Microchannel Cooling System}
\author{Sebastian Blauth$^{*,1,2}$}
\address{$^*$ Corresponding author}
\address{$^1$ Fraunhofer ITWM, Kaiserslautern, Germany}
\email{\href{mailto:sebastian.blauth@itwm.fraunhofer.de}{sebastian.blauth@itwm.fraunhofer.de}}
\author{Christian Leith\"auser$^1$}
\email{\href{mailto:christian.leithaeuser@itwm.fraunhofer.de}{christian.leithaeuser@itwm.fraunhofer.de}}
\author{Ren\'e Pinnau$^2$}
\email{\href{mailto:pinnau@mathematik.uni-kl.de}{pinnau@mathematik.uni-kl.de}}
\address{$^2$ TU Kaiserslautern, Kaiserslautern, Germany}

\begin{abstract}
	We model a microchannel cooling system and consider the optimization of its shape by means of shape calculus. A three-dimensional model covering all relevant physical effects and three reduced models are introduced. The latter are derived via a homogenization of the geometry in 3D and a transformation of the three-dimensional models to two dimensions. A shape optimization problem based on the tracking of heat absorption by the cooler and the uniform distribution of the flow through the microchannels is formulated and adapted to all models. We present the corresponding shape derivatives and adjoint systems, which we derived with a material derivative free adjoint approach. To demonstrate the feasibility of the reduced models, the optimization problems are solved numerically with a gradient descent method. A comparison of the results shows that the reduced models perform similarly to the original one while using significantly less computational resources.
	
	
\end{abstract}

 \maketitle



\vspace{-1cm}
\section{Introduction}
\label{sec:introduction}

For small devices, such as chemical microreactors and electronic equipment, it is critical to have a cooling system that is able to absorb a lot of heat over a small surface area since the performance of these devices is directly related to their operating temperature. For this purpose, cooling systems based on microchannels have been used (see, e.g., \cite{review_channels, khan2011review, naqiuddin2018overview} and the references therein). Their heat transfer coefficient is large due to the high specific surface area of the microchannels and, as they are rather small, they do not increase the overall size of the heat-emitting device too much. Furthermore, the effectiveness of the cooling system critically depends on the uniform distribution of the coolant among all channels. Otherwise, localized zones of high temperature, so-called hot spots, can occur and potentially damage the device. 

The purpose of this paper is to investigate the shape optimization of such a microchannel cooling system. In the literature, this task already received a lot of attention, e.g., in \cite{rao2016dimensional, foli2006optimization, salimpour2013constructal, pan2008optimal, kubo_topology, chentopo, andreasen, soghrati2013computational}. In most of these publications, geometrical properties, such as the shape of the boundary or the topology of the geometry, are parametrized and the optimization of the design is carried out over these parameters, leading to finite dimensional optimization problems. This approach suffers from the obvious drawback that only shapes representable by the parametrization can be reached in the optimization. A more general approach for design optimization consists of using shape or topological sensitivity analysis. These techniques are based on the so-called shape derivative, which measures the sensitivity of a shape due to infinitesimal deformations, and the topological derivative, which measures the sensitivity of a geometry with respect to the insertion of an infinitesimally small hole, see, e.g., \cite{sokolowski_zolesio, delfour_zolesio} for shape calculus and \cite{novotny_sokolowski} for topological sensitivity analysis. In recent years these techniques have been applied to many industrial problems, e.g., the shape design of polymer spin packs \cite{hohmann, leith, leith2, leith3}, electric motors \cite{gangl_shape, gangl_topo}, acoustic horns \cite{schmidt_horn, berggren}, automobiles \cite{othmer, itakura, ferraris}, aircrafts \cite{schmidt2013three, gauger, lyu} or pipe systems \cite{hohmann_erosion, pipe1, pipe2}. To the best of our knowledge, the optimization of a microchannel cooling system by means of shape calculus has only been investigated in our earlier work \cite{blauth}, where we rigorously analyzed the theoretical aspects of this problem. 

To model the cooling system mathematically, we introduce the following models: A three-dimensional model representing the most important physics as well as three reduced models. One of them is a three-dimensional porous medium model based on a homogenization of the domain, similar to those used in, e.g., \cite{chen2007forced, kim1999forced, kim2000local}. For the other two reduced models we use a dimension reduction technique similar to the one of \cite{dimred} that transforms the previous, three-dimensional models to two-dimensional ones. A numerical comparison of the reduced models with the original one shows that they capture the most important physical effects properly, while reducing the computational resources needed considerably.

For the shape optimization, we introduce a cost functional based on the absorption of heat and the uniform distribution of flow through the microchannels. This is adapted to all reduced models using analogous techniques. Subsequently, we present the shape derivatives and adjoint systems for all optimization problems, which we obtained using the material derivative free Lagrange approach of \cite{sturm}. We use these to solve the optimization problems numerically with the help of a gradient descent method, which is again adjusted to all models. A comparison of the optimized geometries suggests that the models give similar results for the shape optimization problem while being substantially more efficient.

This paper is structured as follows: We introduce our three-dimensional model of the cooling system as well as the shape optimization problem in Section~\ref{sec:problem_formulation}. In the three subsequent sections we then present the reduced models respectively. Additionally, the corresponding adaptations for the shape optimization problem are discussed. The dimension reduction technique is introduced and used to derive a two-dimensional model in Section~\ref{sec:dimension_reduction}. Section~\ref{sec:porous_medium} provides the details for the three-dimensional porous medium model of the cooler, to which the dimension reduction approach is then applied in Section~\ref{sec:darcy_2D}, so that we also get a two-dimensional porous medium model. The accuracy of the reduced models is compared numerically in Section~\ref{sec:numerical_state}. Finally, we discuss the implementation as well as the results of the numerical shape optimization in Section~\ref{sec:numerical_optimization}, where we again focus on comparing the different models to each other.

\section{Problem Formulation}
\label{sec:problem_formulation}

First, we describe the geometry of the cooling system and the three-dimensional mathematical model as well as the shape optimization problem. Subsequently, we give a short introduction to shape calculus and then present the shape derivative and adjoint system for the optimization problem.

\subsection{Description of the Geometry}
\label{sec:description_of_the_geometry}

\begin{figure}[b]
	\centering
	\begin{subfigure}[b]{0.49\textwidth}
		\begin{tikzpicture}
			\node at (0,0) {\includegraphics[width=\textwidth, trim=11cm 0cm 0cm 0cm, clip]{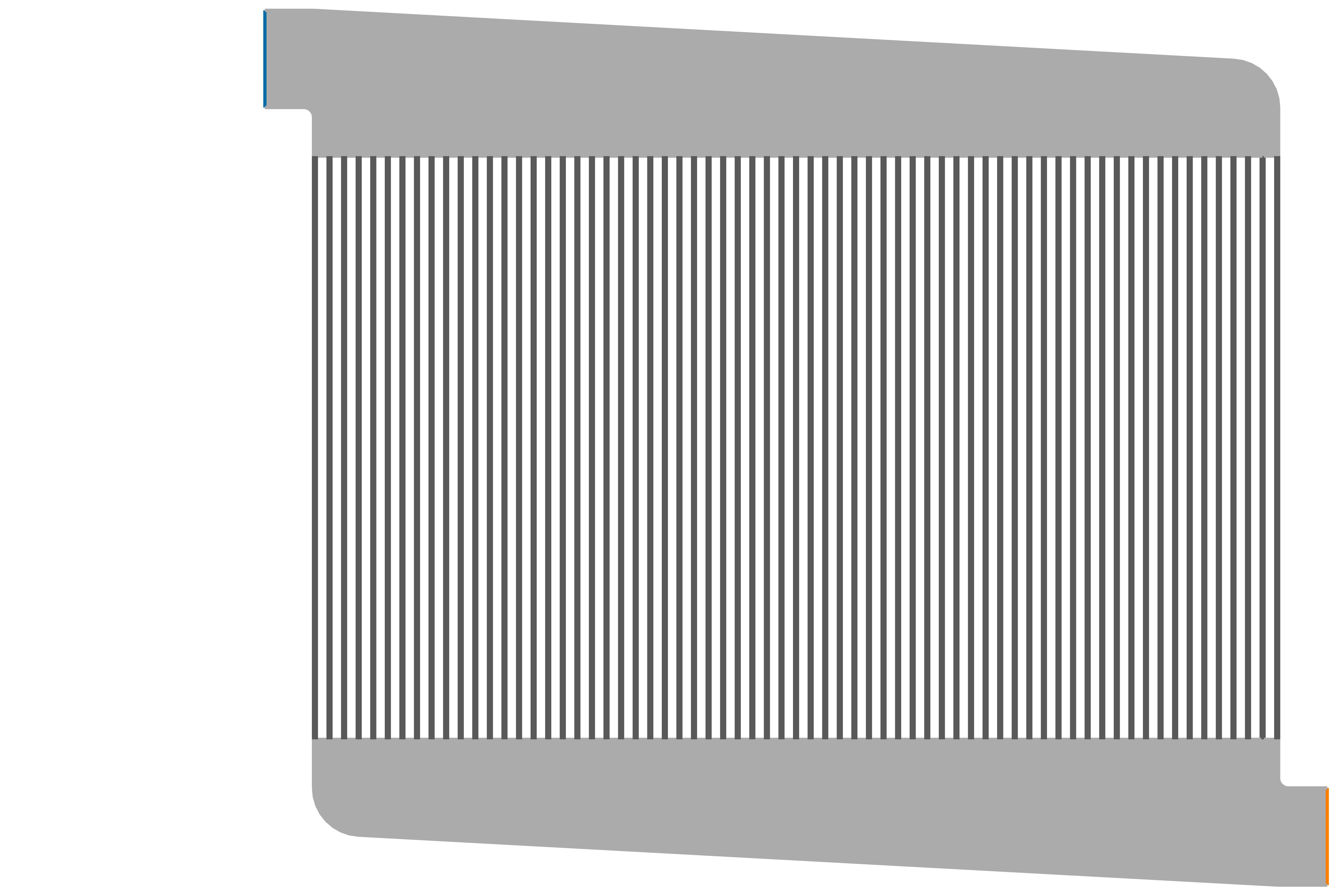}};
			\node at (0, 2.5) {\LARGE $\dimred{\Omega}$};
			\node at (-3, 2.75) {\LARGE \color{tabblue} $\dimred{\Gamma}_\subin$};
			\draw[tabblue, ->, line width=1] (-3.75, 2.75) -- (-3.375, 2.75);
			\node at (4.75, -2.75) {\LARGE \color{taborange} $\dimred{\Gamma}_\subout$};
			\draw[taborange, ->, line width=1] (3.9, -2.75) -- (4.3, -2.75);
			\node at (3, 3.5) {\LARGE \color{tablightgray} $\dimred{\Gamma}_\subwall$};
			\draw[tablightgray, ->, line width=1] (2.75, 3.25) -- (2.7, 2.9);
		\end{tikzpicture}
		\caption{The complete geometry.}
	\end{subfigure}
	\hfill
	\begin{subfigure}[b]{0.4\textwidth}
		\begin{tikzpicture}
			\node at (0,0) {\includegraphics[width=\textwidth, trim= 1cm 0cm 0.5cm 0cm, clip]{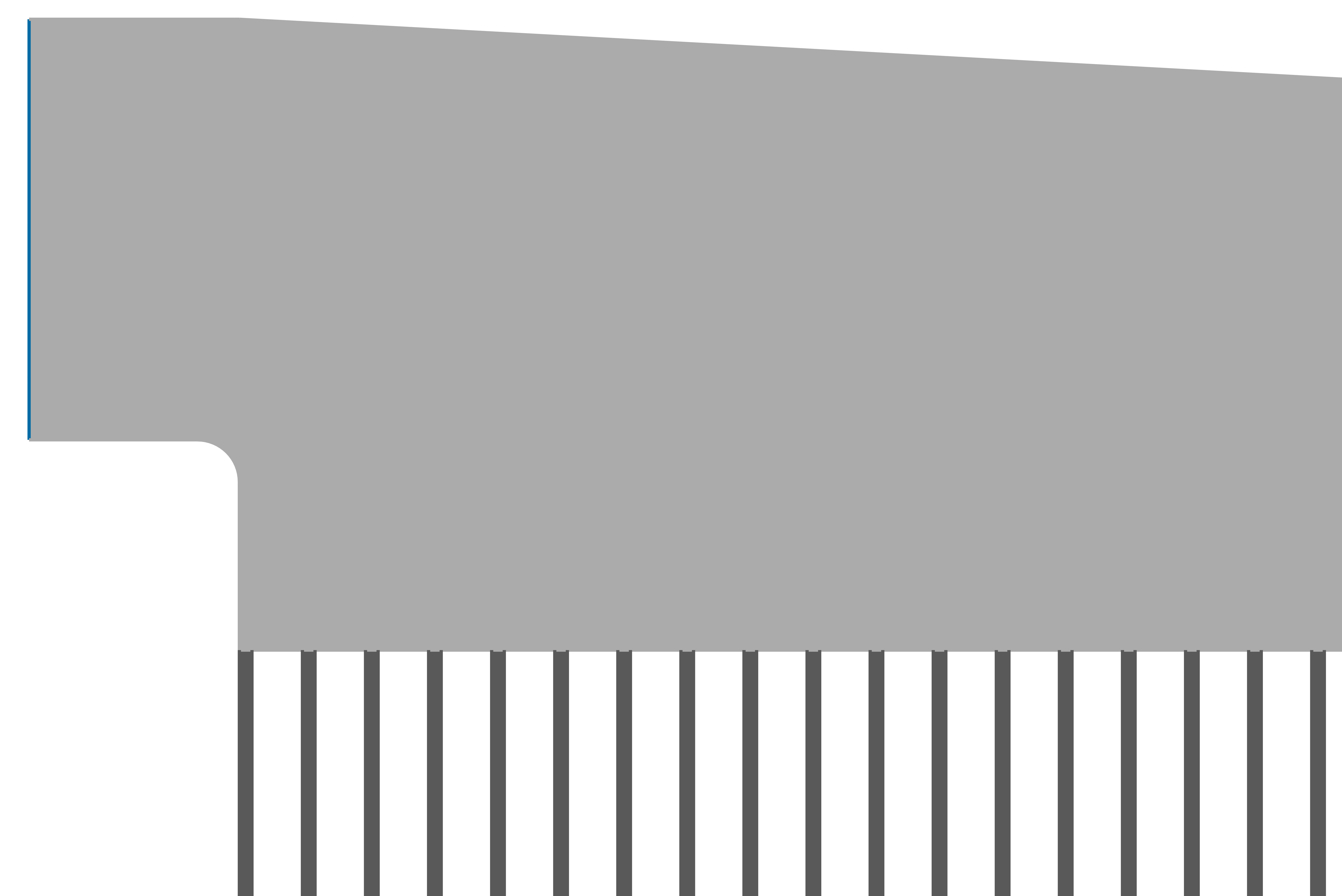}};
			\node at (1,0) {\LARGE \color{tabdarkgray} $\dimred{\Omega}_\subchannels$};
			\draw[tabdarkgray, ->, line width=1] (0.95, -0.45) -- (0.95, -0.85);
			\draw[tabdarkgray, ->, line width=1] (0.7, -0.45) -- (0.35, -0.85);
			\draw[tabdarkgray, ->, line width=1] (1.2, -0.45) -- (1.55, -0.85);
			\node at (-2, 1) {\LARGE \color{tabblue} $\dimred{\Gamma}_\subin$};
			\draw[tabblue, ->, line width=1] (-3, 1) -- (-2.375, 1);
			\draw[tabblue, -, line width=2] (-3.15, 2.095) -- (-3.15, 0.03);
			\node at (2, 3) {\LARGE \color{tablightgray} $\dimred{\Gamma}_\subwall$};
			\draw[tablightgray, ->, line width=1] (2, 2.6) -- (1.95, 2.1);
			\node at (-2.8, -1.75) {\LARGE \color{tabdarkgray} $\dimred{\Gamma}_\subchannels$};
			\draw[tabdarkgray, ->, line width=1] (-2.6, -1.65) -- (-2.2, -1.65);
			\node at (0, -3) {};
		\end{tikzpicture}		
		\caption{Zoom, showing inlet and some microchannels.}
	\end{subfigure}
	\caption{Two-dimensional domain of the cooling system $\dimred{\Omega}$.}
	\label{figure:geometry_full}
\end{figure}
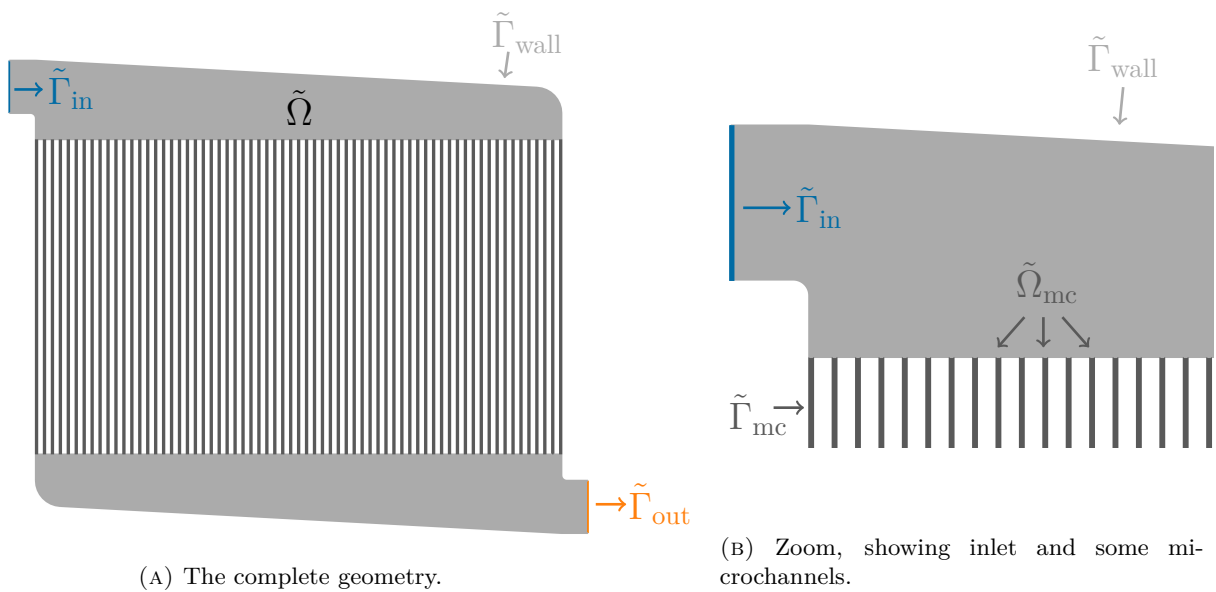


To model the cooling system we consider only the domain of the coolant $\Omega \subset \R^3$, which is assumed to be encased by metal. This metal conducts the heat from the heat source, to which the entire cooler is attached, to the coolant. For simplicity, we call $\Omega$ the geometry or domain of the cooling system throughout the rest of this paper. The domain $\Omega$ has the structure $\Omega = \dimred{\Omega} \times (0, \height)$, where $\dimred{\Omega} \subset \R^2$ is a two-dimensional domain and $\height > 0$ is the constant height of the geometry, i.e., $\Omega$ is an extrusion of $\dimred{\Omega}$ along the $z$-axis. The boundary of $\Omega$ is denoted by $\Gamma$ and is divided into three parts: The inlet $\Gamma_\subin$, where the coolant enters the system, the wall boundary $\Gamma_\subwall$, where heat transfer from the heat source to the cooling system takes place, and the outlet $\Gamma_\subout$, where the cooling liquid leaves the domain. The boundaries inherit the structure of $\Omega$, i.e., $\Gamma_\subin = \dimred{\Gamma}_\subin \times (0, \height)$, $\Gamma_\subout = \dimred{\Gamma}_\subout \times (0, \height)$, and  $\Gamma_\subwall \setminus \left( \Set{z=0} \cup \set{z=\height} \right) = \dimred{\Gamma}_\subwall \times (0, \height)$. The structure of $\Gamma_\subwall$ comes from the fact that the planes $\Set{z=0}$ and $\Set{z=\height}$, i.e., the top and bottom of the cooler and part of $\Gamma_\subwall$, cannot be represented as boundaries of $\dimred{\Omega}$. Of course, it holds that $\dimred{\Gamma} = \dimred{\Gamma}_\subin \cup \dimred{\Gamma}_\subwall \cup \dimred{\Gamma}_\subout$. Additionally, we denote by $\Omega_\subchannels$ the subdomain corresponding to the microchannels, whose wall boundaries are denoted by $\Gamma_\subchannels \subset \Gamma_\subwall$. This situation is depicted in Figure~\ref{figure:geometry_full}, where the domain $\dimred{\Omega}$ is shown. Note, that as the height of the geometry $\height$ is very small compared to the diameter of $\dimred{\Omega}$, we only show slices through $\Omega$ at $z = \nicefrac{\height}{2}$ for all three-dimensional problems. Finally, we remark that such cooling systems can be manufactured, e.g., using wet chemical etching \cite{etching} or 3D-printing \cite{printing}.

\subsection{Mathematical Model}

As we only investigate the steady state of the system, we model all physical processes using stationary equations. The viscous dominated flow of the cooling fluid, which is caused by the slow flow velocities and the small size of the geometry, is modeled by the Stokes equations
\begin{equation}
	\label{eq:stokes}
	\left\lbrace \quad
	\begin{alignedat}{2}
		-\viscosity \laplace \velocity + \grad \pressure &= 0 \quad &&\text{ in } \Omega,\\
		\divergence{\velocity} &= 0 \quad &&\text{ in } \Omega,\\
		\velocity &= \velocity_\subin \quad &&\text{ on } \Gamma_\subin,\\
		\velocity &= 0 \quad &&\text{ on } \Gamma_\subwall,\\
		\viscosity \partial_\normal \velocity - \pressure \normal &= 0 \quad &&\text{ on } \Gamma_\subout,
	\end{alignedat}
	\right.
\end{equation}
where $\velocity$ and $\pressure$ denote the fluid velocity and pressure, respectively. Furthermore, $\viscosity$ denotes the dynamic viscosity of the coolant and $\partial_\normal v$ is the normal derivative of a function $v$ given by $\partial_\normal v = Dv\ \normal$, where $Dv$ is the Jacobian of $v$ and $\normal$ is the outward unit normal on $\Gamma$. We model the flow of the coolant into $\Omega$ with the fixed inflow condition $\velocity = \velocity_\subin$ on $\Gamma_\subin$. Here, $\velocity_\subin$ is the inflow velocity which is chosen to be in the direction of the inward facing normal and for the sake of simplicity we assume that $\velocity_\subin\in H^1(\Omega)^3$ with $\velocity_\subin = 0$ on $\Gamma_\subwall$, as in our theoretical analysis \cite{blauth}. For the wall boundary $\Gamma_\subwall$ we use the no-slip condition and on the outlet $\Gamma_\subout$ we use the do-nothing condition that models the unimpaired outflow of the fluid (see, e.g., \cite{john}).

The temperature of the coolant changes due to both conduction and convection. This is modeled by the convection-diffusion equation
\begin{equation}
	\label{eq:convection-diffusion}
	\left\lbrace \quad
	\begin{alignedat}{2}
		- \grad \cdot \left( \conductivity \grad \temperature \right) + \density \hcapacity\ \velocity \cdot \grad \temperature &= 0 \quad &&\text{ in } \Omega,\\
		\temperature &= \temperature_\subin \quad &&\text{ on } \Gamma_\subin,\\
		\conductivity \partial_\normal \temperature + \htc \left( \temperature - \temperature_\subwall \right) &= 0 \quad &&\text{ on } \Gamma_\subwall,\\
		\conductivity \partial_\normal \temperature &= 0 \quad &&\text{ on } \Gamma_\subout,
	\end{alignedat}
	\right.
\end{equation}
where $\temperature$ denotes the fluid's temperature and $\velocity$ is its velocity, i.e., the solution of \eqref{eq:stokes}. Furthermore, $\conductivity, \density, \text{ and } \hcapacity$ denote the coolant's thermal conductivity, density, and specific heat capacity, respectively. The temperature of the inflowing fluid is fixed to $\temperature_\subin$ on $\Gamma_\subin$ by a Dirichlet condition. Analogously to before, we assume that $\temperature_\subin \in H^1(\Omega)$ (cf. \cite{blauth}).
The heat transfer over the wall boundary is modeled by the Robin boundary condition on $\Gamma_\subwall$, where $\htc$ denotes the heat transfer coefficient between coolant and wall. For simplicity, we assume that the heat source generates a temperature distribution that yields a constant temperature $\temperature_\subwall$ on $\Gamma_\subwall$, but non-constant temperature distributions could also be considered. The behavior of the temperature at $\Gamma_\subout$ is modeled by a homogeneous Neumann condition. 
The values of the physical parameters for our particular problem setting are given in Table~\ref{table:physical_constants}. They are chosen to resemble a viscous coolant that could be used for the application in a chemical microreactor. The heat transfer coefficient is chosen so that the amount of energy absorbed by the system resembles the amount of heat generated by the chemical reaction. In particular, we consider a setting similar to the one described in \cite{methanation}, where the Sabatier process in microreactors is investigated. 


\begin{table}[b]
	\centering
	{\footnotesize
	\rowcolors{2}{\tablegray}{white}
	\begin{tabular}{l S c l S}
		\toprule
		parameter [unit] & {value} & \hspace{1.5em} & parameter [unit] & {value} \\
		\midrule
		dynamic viscosity $\viscosity$ [\si{\kilogram \per \meter \per \second}] & 3e-4 & & inflow temperature $\temperature_\subin$ [\si{\celsius}] & 300 \\
		density $\density$ [\si{\kilogram \per \cubic \meter}] & 700 & & wall temperature $\temperature_\subwall$ [\si{\celsius}] & 400 \\
		mass inflow $\massin$ [\si{\kilogram \per \second}] & 6e-5 & & heat transfer coefficient $\htc$ [\si{\watt \per \square \meter \per \kelvin}] & 10 \\
		thermal conductivity $\conductivity$ [\si{\watt \per \meter \per \kelvin}] & 0.1 & & height of the geometry $\height$ [\si{\meter}] & 3e-4 \\
		specific heat capacity $\hcapacity$ [\si{\joule \per \kilogram \per \kelvin}] & 2000 & & & \\
		\bottomrule
	\end{tabular}
	\caption{Physical parameters for the cooling fluid.}
	\label{table:physical_constants}
	}
\end{table}

As the notion of strong solutions of PDEs is usually too strict for the means of PDE constrained optimization problems (see, e.g., \cite{troeltzsch, hinze_pinnau_ulbrich2}), we consider equations \eqref{eq:stokes} and \eqref{eq:convection-diffusion} in their weak form. To this end, we define the (affine) Sobolev spaces
\begin{equation}
	\label{eq:spaces_3D}
	\left\lbrace\quad
		\begin{aligned}
			\fspacetest &:= \set{\testvelo \in H^1(\Omega)^3 | \testvelo=0 \text{ on } \Gamma_\subin \cup \Gamma_\subwall}, \qquad \fspacetrial := \velocity_\subin + \fspacetest, \qquad \pspace := L^2(\Omega), \\
			\tspacetest &:= \set{\testtemp \in H^1(\Omega) | \testtemp=0 \text{ on } \Gamma_\subin}, \qquad \tspacetrial := \temperature_\subin + \tspacetest, \\
			\statespace &:= \fspacetrial \times \pspace \times \tspacetrial, \qquad \adjointspace := \fspacetest \times \pspace \times \tspacetest.
		\end{aligned}
	\right.
\end{equation}
The weak form of \eqref{eq:stokes} and \eqref{eq:convection-diffusion} is then given by
\begin{equation}
	\label{eq:weak_state}
	\left\lbrace \quad
	\begin{aligned}
		&\text{Find } \solution = (\velocity, \pressure, \temperature) \in \statespace \text{ such that } \\
		&\qquad \integral{\Omega} \viscosity D \velocity : D \testvelo - \pressure\ \divergence{\testvelo} - \testpres\ \divergence{\velocity} + \conductivity \grad \temperature \cdot \grad \testtemp + \density \hcapacity\ \velocity\cdot \grad \temperature\ \testtemp \dx{x} \\
		&\qquad + \integral{\Gamma_\subwall} \htc (\temperature - \temperature_\subwall) \testtemp \dx{s} = 0 \\
		&\text{for all } \testsolution = (\testvelo, \testpres, \testtemp) \in \adjointspace.
	\end{aligned}
	\right.
\end{equation}
If $\Omega$ has a Lipschitz boundary and the data is sufficiently smooth, it is not too difficult to show that problem \eqref{eq:weak_state} has a unique weak solution that depends continuously on the data (see, e.g., \cite{blauth, ern_guermond, roos}).

\subsection{The Optimization Problem}
\label{sec:optimization_problem}

Our goal is to improve the cooling system by means of shape optimization, i.e., we want to optimize it by changing its shape without altering its topology. To do so, we assume that we can only change the shape of $\Gamma_\subwall$ and that $\Gamma_\subin$ and $\Gamma_\subout$ remain fixed. This is reasonable since the cooling system is connected to other parts via in- and outlet. Furthermore, we also do not change the underlying structure of the geometry since we assume that the height $\height$ of the cooler is fixed, and we also only allow the microchannels to change in length.

The quality of the cooling system is measured by the following criteria. First, it should absorb a specific amount of heat so that a finely detailed cooling of the heat source is possible. Second, we want to find a geometry that distributes the coolant uniformly to the microchannels as this helps preventing hot spots. To model this, we define the cost functional as
\begin{equation}
	\label{eq:def_cost_function}
	\costfunction(\Omega, \solution) := \weighttemp \costfunction_1(\Omega, \solution) + \weightvelo \costfunction_2(\Omega,\solution) + \weightreg \costfunction_3(\Omega,\solution),
\end{equation}
with
\begin{equation*}
\costfunction_1(\Omega,\solution) := \Big( \flux(\Omega, \solution) - \fluxdes \Big)^2, \quad \costfunction_2(\Omega, \solution) := \integral{\Omega_\subchannels} \abs{\velocity - \velocity_\subdes}^2 \dx{x}, \quad \costfunction_3(\Omega, \solution) := \integral{\Gamma} 1 \dx{s}.
\end{equation*}
Here, $\solution = (\velocity, \pressure, \temperature)$ is the vector of state variables and the weights $\lambda_i$ are nonnegative. Further, we define 
\begin{equation*}
	\flux(\Omega, \solution) := \integral{\Gamma_\subwall} \htc (\temperature_\subwall - \temperature) \dx{s}.
\end{equation*}
Finally, we consider the following shape optimization problem
\begin{equation}
	\label{eq:def_opt_problem}
	\min_{\Omega, \solution}\ \costfunction(\Omega, \solution) \quad \text{ subject to } \eqref{eq:weak_state}.
\end{equation}

For the term $\costfunction_1$ we have the following considerations. The heat flux on $\Gamma_\subwall$ is given by $-\conductivity \grad \temperature$ as the fluid velocity vanishes there. Thus, the overall energy entering the domain is given by $\integral{\Gamma_\subwall} \conductivity \grad \temperature \cdot \normal \dx{s}$. Using the Robin boundary condition on $\Gamma_\subwall$ (cf. \eqref{eq:convection-diffusion}), we can express the amount of energy entering the domain via $\flux(\Omega, \solution)$. Therefore, by minimizing $\costfunction_1$ we try to get a cooling system that absorbs a particular amount of heat, given by $\fluxdes$. This is particularly useful for cooling devices that rely on a specific operating temperature, such as chemical reactors.


The term $\costfunction_2$ has the purpose of minimizing the distance of the fluid velocity $\velocity$ to some desired velocity $\velocity_\subdes$ in $L^2(\Omega_\subchannels)$, where $\Omega_\subchannels$ is the domain corresponding to the microchannels (cf. Section~\ref{sec:description_of_the_geometry}). In particular, we choose $\velocity_\subdes$ as the velocity profile corresponding to the case of uniformly distributed flow through each channel, which can be computed as Poiseuille flow for a channel with square base (see, e.g., \cite[Chapter 3.4]{bruus}). This term of the cost functional leads to a uniform flow distribution which helps to minimize the occurrence of hot spots.

Finally, the term $\costfunction_3$ is a so-called perimeter regularization (see, e.g., \cite{welker, sturm}), that penalizes the increase of surface area. This is used since it can grant the existence of a minimizer, as discussed in \cite{perimeter_regularization}.

As we assume that the state system \eqref{eq:weak_state} has a unique weak solution for every domain $\Omega$ with a Lipschitz boundary, we now introduce the so-called reduced cost functional $\reducedcostfunction$ by
\begin{equation}
	\label{eq:def_reduced_cost_function}
	\reducedcostfunction(\Omega) = \costfunction(\Omega, \solution(\Omega)),
\end{equation}
where $\solution(\Omega)$ denotes the solution of \eqref{eq:weak_state} in $\Omega$. With this, it is evident that \eqref{eq:def_opt_problem} is equivalent to the reduced optimization problem
\begin{equation*}
	\min_{\Omega}\ \reducedcostfunction(\Omega),
\end{equation*}
where the PDE constraint is formally eliminated.

\subsection{Shape Calculus}
\label{sec:shape_calculus}

To solve the shape optimization problem numerically, we calculate the shape derivative of the reduced cost functional \eqref{eq:def_reduced_cost_function}, which is then used in a gradient descent method (see also Section~\ref{sec:numerical_optimization}). For a detailed introduction to shape calculus we refer to the monographs \cite{delfour_zolesio, sokolowski_zolesio}. To compute the shape derivative, we utilize the so-called speed method, which is also known as velocity method. It uses the flow of a vector field $\vectorfield$ to deform a domain $\Omega \subset D$, where $D \subset{\R^d}$ is a so-called hold-all domain, i.e., a domain that contains all domains under consideration, and $d$ denotes the dimension of the geometry. For a given vector field $\vectorfield \in C^k_0(D; \R^d)$ with $k\geq 1$, i.e., the space of all $k$-times continuously differentiable functions from $D$ to $\R^d$ having compact support in $D$, the evolution of a point $x\in \Omega$ under the flow of $\vectorfield$ is given by the solution of the ODE
\begin{equation*}
	\dot{x}(t) = \vectorfield(x(t)), \quad x(0) = x.
\end{equation*}
We know that this system has a unique solution $x(t)$ for $t \in [0, \tau]$ if $\tau > 0$ is sufficiently small. Hence, we define the flow of $\vectorfield$ as the diffeomorphism
\begin{equation*}
	\flow \colon \R^d \to \R^d; \quad x \mapsto \flow (x) = x(t).
\end{equation*}
Now we define the shape derivative as in \cite{sturm}, where we write $2^D := \Set{\Omega | \Omega \subset D}$ for the power set of $D$.
\begin{Definition}
	Let $\mathcal{S} \subset 2^D$, $\costfunction\colon \mathcal{S} \to \R$ and $\Omega \in \mathcal{S}$. Further, let $\vectorfield \in C^k_0(D;\R^d)$ with $k \geq 1$ and associated flow $\flow$ and suppose that $\flow(\Omega) \in \mathcal{S}$ for all $t\in [0,\tau]$ with $\tau > 0$ sufficiently small. The Eulerian semi-derivative of $\costfunction$ at $\Omega$ in direction $\vectorfield$ is given by the limit
	 \begin{equation*}
	 	d\costfunction(\Omega)[\vectorfield] := \lim\limits_{t\searrow 0} \frac{\costfunction(\flow(\Omega)) - \costfunction(\Omega)}{t},
	 \end{equation*}
	 if it exists. The function $\costfunction$ is called shape differentiable at $\Omega \in \mathcal{S}$ if the above limit exists for all $\vectorfield \in C^\infty_0(D;\R^d)$ and if the mapping 
	 \begin{equation*}
	 	C^\infty_0(D;\R^d) \to \R; \quad \vectorfield \mapsto d\costfunction(\Omega)[\vectorfield],
	 \end{equation*}
	 is linear and continuous. Then, $d\costfunction(\Omega)[\vectorfield]$ is called the shape derivative of $\costfunction$ at $\Omega$ in direction $\vectorfield$.
\end{Definition}

\subsection{Formal Shape Optimization}
\label{sec:formal_shape_optimization}

Our goal is to calculate the shape derivative of the reduced cost functional \eqref{eq:def_reduced_cost_function}, which is done with an adjoint approach (see, e.g., \cite{sturm, delfour_zolesio} for more details). The rigorous verification of the assumptions for this approach is beyond the scope of this paper and can be found in our earlier work \cite{blauth}. We define the Lagrangian associated to \eqref{eq:def_reduced_cost_function} by
\begin{equation*}
	\left\lbrace \quad
	\begin{alignedat}{2}
		\mathcal{L}(\Omega, \solution, \adsolution) =\ &\costfunction(\Omega, \solution) + \integral{\Omega} \viscosity D \velocity : D \advelo - \pressure\ \divergence{\advelo} - \adpres\ \divergence{\velocity} \dx{x} \\
		+\ &\integral{\Omega} \conductivity \grad \temperature \cdot \grad \adtemp + \density \hcapacity\ \velocity\cdot \grad \temperature\ \adtemp \dx{x} + \integral{\Gamma_\subwall} \htc (\temperature - \temperature_\subwall) \adtemp \dx{s},
	\end{alignedat}
	\right.
\end{equation*}
where $\solution = (\velocity, \pressure, \temperature) \in \statespace$ and $\adsolution = (\advelo, \adpres, \adtemp) \in \adjointspace$. In particular, we have $\mathcal{L}(\Omega, \solution(\Omega), \psi) = \reducedcostfunction(\Omega)$ for all $\psi \in \adjointspace$ as $\solution(\Omega)$ is the solution of the state system \eqref{eq:weak_state}. Further, for a sufficiently smooth vector field $\vectorfield $ with associated flow $\flow$, we define $\Omega_t := \flow(\Omega)$. With this, we introduce the associated shape Lagrangian by
\begin{equation}
	\label{eq:shape_lagrangian}
	\shapelagrangian \colon [0,\tau] \times \statespace \times \adjointspace \to \R; \quad \shapelagrangian(t, \varphi, \psi) = \mathcal{L}(\Omega_t, \varphi \circ \flow^{-1}, \psi \circ \flow^{-1}),
\end{equation}
which corresponds to the original Lagrangian $\mathcal{L}$ on the perturbed domain $\Omega_t$. We denote by $\solution_t = \solution(\Omega_t)$ the solution of the state system on the domain $\Omega_t$ and define $\solution^t = \solution_t \circ \Phi_t$. With this, a similar argumentation as in \cite{sturm_new_trends} shows that we have $\shapelagrangian(t, \solution^t, \psi) = \reducedcostfunction(\Omega_t)$ for all $\psi\in \adjointspace$, because $\flow$ is a diffeomorphism. Hence, we can calculate the shape derivative of the reduced cost functional using
\begin{equation*}
	d\reducedcostfunction(\Omega)[\vectorfield] = \lim_{t\searrow 0} \frac{\reducedcostfunction(\flow(\Omega)) - \reducedcostfunction(\Omega)}{t} = \left. \frac{d}{dt} \reducedcostfunction(\Omega_t) \right|_{t=0} = \left. \frac{d}{dt} \shapelagrangian(t, \solution^t, \psi) \right|_{t=0}.
\end{equation*}
Additionally, under the assumptions of \cite[Theorem~3.1]{sturm} it holds that
\begin{equation}
	\label{eq:sturm_derivative}
	\left. \frac{d}{dt} \shapelagrangian(t, \solution^t, \psi) \right|_{t=0} = \left. \partial_t \shapelagrangian(t, \solution, \adsolution) \right|_{t=0},
\end{equation}
where $\solution$ is the solution of the state system, i.e., \eqref{eq:weak_state}, and $\adsolution$ is the solution of the adjoint system
\begin{equation*}
	\text{Find } \adsolution \in \adjointspace \text{ such that} \qquad \partial_\varphi \shapelagrangian(0, \solution, \adsolution)[\hat{\varphi}] = 0 \qquad \text{for all } \hat{\varphi} \in \adjointspace.
\end{equation*}
It is straightforward to see that this adjoint system is equivalent to
\begin{equation}
	\label{eq:weak_adjoint}
	\left\lbrace \quad
	\begin{aligned}
		&\text{Find } \adsolution = (\advelo, \adpres, \adtemp) \in \adjointspace \text{ such that } \\
		&\qquad \integral{\Omega} \conductivity \grad \adtemp \cdot \grad \testadtemp + \density \hcapacity\ \velocity \cdot \grad \testadtemp \adtemp \dx{x} + \integral{\Gamma_\subwall} \htc \adtemp \testadtemp \dx{s} - 2 \weighttemp  \Big( \flux(\Omega, \solution(\Omega)) - \fluxdes \Big) \integral{\Gamma_\subwall} \htc \testadtemp \dx{s} \\
		&\qquad + \integral{\Omega} \viscosity\ D \advelo : D \testadvelo - \adpres\ \divergence{\testadvelo} - \testadpres\ \divergence{\advelo} + \density \hcapacity \testadvelo \cdot \grad \temperature \adtemp \dx{x} =  -2 \weightvelo \integral{\Omega_\subchannels} (\velocity - \velocity_\subdes) \cdot \testadvelo \dx{x} \\
		&\text{for all } \testadsolution = (\testadvelo, \testadpres, \testadtemp) \in \adjointspace,
	\end{aligned}
	\right.
\end{equation}
where $(\velocity, \pressure, \temperature) = \solution(\Omega)$ is the weak solution of the state system \eqref{eq:weak_state}. To calculate the shape derivative using \eqref{eq:sturm_derivative}, we use the transformation formula to rewrite $\shapelagrangian(t,\solution, \adsolution)$, where we pull-back the integrals to the domain $\Omega$. To this end, we define 
\begin{equation*}
	\xi(t) := \det\left( D\Phi_t \right), \quad \omega(t) := \xi(t) \abs{D \Phi_t^{-T} \normal}, \quad A(t) := \xi(t) D\Phi_t^{-1} D\Phi_t^{-T},\quad B(t) := \xi(t) D\Phi_t^{-T},
\end{equation*}
where $D\Phi_t^{-1}$ denotes the inverse of $D\Phi_t$, and not the Jacobian of $\Phi_t^{-1}$.
Additionally, we note that
\begin{equation*}
	D (f \circ \Phi_t^{-1}) = (Df\ D\Phi_t^{-1})\circ \Phi_t^{-1}, \quad \integral{\Phi_t(\Omega)} \hspace{-1em} f \dx{x} = \integral{\Omega} f\circ \Phi_t\ \xi(t) \dx{x} \quad \integral{\Phi_t(\Gamma)} \hspace{-1em} f \dx{s} = \integral{\Gamma} f\circ \Phi_t\ \omega(t) \dx{s},
\end{equation*}
for all sufficiently smooth functions $f$ and vector fields $\vectorfield$ for $t \in [0,\tau]$ with $\tau>0$ being sufficiently small (see, e.g., \cite{sturm, delfour_zolesio}). With this, we rewrite $\shapelagrangian$ as follows
\begin{equation}
	\label{eq:pulled_back}
	\left\lbrace \quad
	\begin{aligned}
		&\shapelagrangian(t, \solution, \adsolution) \\
		=\ &\weighttemp \left( \integral{\Gamma_\subwall} \htc(\temperature_\subwall - \temperature)\ \omega(t) \dx{s} - \fluxdes \right)^2 + \weightvelo \integral{\Omega_\subchannels} \abs{\velocity - \velocity_\subdes \circ \Phi_t}^2 \xi(t) \dx{x} \\
		&+ \weightreg \integral{\Gamma} \omega(t) \dx{s} + \integral{\Omega} \viscosity \left(D\velocity\ A(t)\right) : D\advelo \dx{x} - \integral{\Omega} \pressure\ \tr\left(D\advelo\ B(t)^T\right) \dx{x} - \integral{\Omega} \adpres\ \tr\left(D\velocity\ B(t)^T\right) \dx{x} \\
		&+ \integral{\Omega} \conductivity \left( A(t) \grad \temperature \right) \cdot \grad \adtemp \dx{x} + \integral{\Omega} \density\hcapacity \velocity \cdot (B(t)\ \grad T)\ \adtemp \dx{x} + \integral{\Gamma_\subwall} \htc (\temperature - \temperature_\subwall) \adtemp\ \omega(t) \dx{s}.
	\end{aligned}
	\right.
\end{equation}
Further, due to \cite{delfour_zolesio, sturm} we have the following derivatives
\begin{equation*}
	\xi'(0^+) = \divergence{\vectorfield}, \quad \omega'(0^+) = \tdiv{\vectorfield}, \quad A'(0^+) = \divergence{\vectorfield} I - 2\varepsilon(\vectorfield), \quad B'(0^+) = \divergence{\vectorfield} I - D\vectorfield^T,
\end{equation*}
where we write $f'(0^+) = \lim_{t \searrow 0} \nicefrac{1}{t} (f(t) - f(0))$ for some differentiable function $f$. Further, $I$ denotes the identity matrix in $\R^d$, $\varepsilon(\vectorfield) = \nicefrac{1}{2} ( D\vectorfield + D\vectorfield^T )$ is the symmetric part of the Jacobian of $\vectorfield$ and $\mathrm{div}_\Gamma$ denotes the tangential divergence defined as $\tdiv{\vectorfield} = \divergence{\vectorfield} - D\vectorfield\ \normal \cdot \normal$. Finally, differentiating \eqref{eq:pulled_back} w.r.t. $t$ and using \eqref{eq:sturm_derivative} proves the following theorem (cf. \cite{blauth}).

\begin{Theorem}
	\label{prop:full_3D}
	Under the assumptions of \cite[Theorem~3.1]{sturm} the shape derivative of \eqref{eq:def_reduced_cost_function} at $\Omega \subset D$ in direction $\vectorfield$ is given by
	\begin{equation}
		\label{eq:shape_derivative}
		\left\lbrace \quad
		\begin{aligned}
			d\reducedcostfunction(\Omega)[\vectorfield] =\ &2\weighttemp \Big(\flux(\Omega, \solution(\Omega)) - \fluxdes \Big) \integral{\Gamma_\subwall} \htc (\temperature_\subwall - \temperature) \tdiv{\vectorfield} \dx{s} \\
			&+ \weightvelo \integral{\Omega_\subchannels} \abs{\velocity - \velocity_\subdes}^2 \divergence{\vectorfield} - 2(\velocity - \velocity_\subdes) \cdot D\velocity_\subdes\ \vectorfield \dx{x} + \weightreg \integral{\Gamma} \tdiv{\vectorfield} \dx{s} \\
			&+ \integral{\Omega} \viscosity \left( D\velocity (\divergence{\vectorfield} I - 2\varepsilon(\vectorfield))\right) : D\advelo \dx{x} - \integral{\Omega} \pressure\ \tr \left( D\advelo \left( \divergence{\vectorfield}I - D\vectorfield \right) \right) \dx{x} \\
			&- \integral{\Omega} \adpres\ \tr \left( D\velocity \left( \divergence{\vectorfield}I - D\vectorfield \right) \right) \dx{x} + \integral{\Omega} \conductivity \left( \left( \divergence{\vectorfield}I - 2\varepsilon(\vectorfield) \right) \grad \temperature \right) \cdot \grad \adtemp \dx{x} \\
			&+ \integral{\Omega} \density \hcapacity \velocity \cdot \left( \left( \divergence{\vectorfield}I - D\vectorfield^T \right) \grad \temperature \right) \adtemp \dx{x} + \integral{\Gamma_\subwall} \htc \left( \temperature - \temperature_\subwall \right) \adtemp \tdiv{\vectorfield} \dx{s},
		\end{aligned}
		\right.
	\end{equation}
	where $(\velocity, \pressure, \temperature) = \solution(\Omega)$ is the solution of the state system \eqref{eq:weak_state} and $(\advelo, \adpres, \adtemp) = \adsolution(\Omega)$ is the solution of the adjoint system \eqref{eq:weak_adjoint}.
\end{Theorem}

\section{The Dimension Reduction Technique}
\label{sec:dimension_reduction}

For our first reduced model, we introduce a dimension reduction technique, similar to the idea discussed in \cite{dimred} for the topology optimization of Stokes flow which was also used, e.g., in \cite{yan2019topology, sigmund}. The technique transforms the three-dimensional state system to a two-dimensional one, making its numerical solution considerably easier. Afterwards, we discuss this technique in the context of the shape optimization problem. 

\subsection{Description of the Model}

The main idea of this model is to exploit the structure of our domain, i.e., we use that $\Omega = \dimred{\Omega} \times (0, \height)$. As this model is based on the weak formulation of the PDEs, we now introduce the following two-dimensional (affine) function spaces in analogy to the ones given in \eqref{eq:spaces_3D}
\begin{equation*}
	\left\lbrace \quad
	\begin{aligned}
		\dimred{\fspacetest} &:= \Set{\dimred{\testvelo} \in H^1(\dimred{\Omega})^2 | \dimred{\testvelo}=0 \text{ on } \dimred{\Gamma}_\subin \cup \dimred{\Gamma}_\subwall}, \qquad \dimred{\fspacetrial} := \dimred{\velocity}_\subin + \dimred{\fspacetest}, \qquad \dimred{\pspace} := L^2(\dimred{\Omega}), \\
		\dimred{\tspacetest} &:= \Set{\dimred{\testtemp} \in H^1(\dimred{\Omega}) | \dimred{\testtemp}=0 \text{ on } \dimred{\Gamma}_\subin}, \qquad \dimred{\tspacetrial} := \dimred{\temperature}_\subin + \dimred{\tspacetest}, \\
		\dimred{\statespace} &:= \dimred{\fspacetrial} \times \dimred{\pspace} \times \dimred{\tspacetrial}, \qquad \dimred{\adjointspace} := \dimred{\fspacetest} \times \dimred{\pspace} \times \dimred{\tspacetest},
	\end{aligned}
	\right.
\end{equation*}
where we assume that $\dimred{\velocity}_\subin \in H^1(\dimred{\Omega})^2$ with $\dimred{\velocity}_\subin = 0$ on $\dimred{\Gamma}_\subwall$ and $\dimred{\temperature}_\subin \in H^1(\dimred{\Omega})$.

The key observation for the dimension reduction of the Stokes equation is the following. For viscous flow between parallel, infinite planes we get a parabolic velocity profile in analogy to two-dimensional Poiseuille flow (cf. \cite{bruus}). This situation is not very different from our setting. We still consider the flow of a viscous fluid between two parallel plates, albeit the geometry is more complex. Due to this observation we make the following assumption: there is no fluid velocity in $z$-direction and the $x$ and $y$ components of the velocity are parabolic in $z$, i.e., we describe the fluid velocity $\velocity = [\velocity_1, \velocity_2, \velocity_3]^T$ as 
\begin{equation}
	\label{eq:parabolic_profile_ansatz}
	\velocity(x,y,z) = \frac{4}{\height^2} z(\height-z)
	\begin{bmatrix}
		\dimred{\velocity}_1(x,y) & \dimred{\velocity}_2(x,y) & 0
	\end{bmatrix}^T = \frac{4}{\height^2} z(\height-z)\begin{bmatrix}
		\dimred{\velocity}(x,y) & 0
	\end{bmatrix}^T,
\end{equation}
where $\dimred{\velocity} = [\dimred{\velocity}_1, \dimred{\velocity}_2]^T \in \dimred{\fspacetrial}$ is the 2D velocity. The factor $\nicefrac{4}{\height^2}$ is chosen so that $\dimred{\velocity}$ represents the maximum fluid velocity. Furthermore, we see that the no-slip boundary conditions at $\Set{z=0}$ and $\Set{z=\height}$ are satisfied by \eqref{eq:parabolic_profile_ansatz}. For the fluid temperature and pressure we assume that they are constant along the $z$-direction, i.e., 
\begin{equation}
	\label{eq:constant_profile_ansatz}
	\pressure(x,y,z) = \dimred{\pressure}(x,y) \qquad \text{ and } \qquad \temperature(x,y,z) = \dimred{\temperature}(x,y),
\end{equation}
where $\dimred{\pressure} \in \dimred{\pspace}$ and $\dimred{\temperature} \in \dimred{\tspacetrial}$.
For the pressure this is reasonable as we do not have a flow in $z$-direction (cf. \eqref{eq:parabolic_profile_ansatz}) and, hence, we easily get $\partial_z \pressure = 0$ from the Stokes equation \eqref{eq:stokes}. For the temperature the situation is more complicated and, in general, it is not constant in $z$. However, as our problem is convection-dominated and there is no convection in $z$-direction due to \eqref{eq:parabolic_profile_ansatz}, we observe that the variation in temperature is significantly higher in the $x$\nobreakdash-$y$-plane than it is in $z$-direction, and the assumptions \eqref{eq:constant_profile_ansatz} are justified. Naturally, we assume that $\dimred{\velocity}_\subin$ and $\dimred{\temperature}_\subin$ are obtained from $\velocity_\subin$ and $\temperature_\subin$ analogously to \eqref{eq:parabolic_profile_ansatz} and \eqref{eq:constant_profile_ansatz}.

As we want to derive a two-dimensional weak formulation, we introduce the functions
\begin{equation*}
	\testvelo(x,y,z) = \frac{4}{\height^2} z(\height-z)
	\begin{bmatrix}
		\dimred{\testvelo}_1(x,y) & \dimred{\testvelo}_2(x,y) & 0
	\end{bmatrix}^T,
\end{equation*}
where $\dimred{\testvelo} = [\dimred{\testvelo}_1, \dimred{\testvelo}_2]^T \in \dimred{\fspacetest}$ as well as
\begin{equation*}
    \testpres(x,y,z) = \dimred{\testpres}(x,y) \quad \text{ and } \quad \testtemp(x,y,z) = \dimred{\testtemp}(x,y),
\end{equation*}
in analogy to \eqref{eq:parabolic_profile_ansatz} and \eqref{eq:constant_profile_ansatz}.

Using these functions in the weak form \eqref{eq:weak_state} and performing the integration w.r.t.\ $z$ yields the two-dimensional system
\begin{equation}
	\label{eq:weak_state_dimension_reduction}
	\left\lbrace \quad
	\begin{aligned}
		&\text{Find } \dimred{\solution} = (\dimred{\velocity}, \dimred{\pressure}, \dimred{\temperature}) \in \dimred{\statespace} \text{ such that } \\
		&\qquad \integral{\dimred{\Omega}} \frac{8\height}{15}\ \viscosity D \dimred{\velocity} : D \dimred{\testvelo} + \frac{16}{3\height}\ \viscosity\ \dimred{\velocity}\cdot \dimred{\testvelo} - \frac{2\height}{3}\ \dimred{\pressure}\ \divergence{\dimred{\testvelo}} - \frac{2\height}{3}\ \dimred{\testpres}\ \divergence{\dimred{\velocity}} \dx{x} \\
		&\qquad + \integral{\dimred{\Omega}} \height\ \conductivity \grad \dimred{\temperature} \cdot \grad \dimred{\testtemp} + \frac{2\height}{3}\ \density \hcapacity\ \dimred{\velocity}\cdot \grad \dimred{\temperature}\ \dimred{\testtemp} + 2 \htc (\dimred{\temperature} - \temperature_\subwall) \dimred{\testtemp} \dx{x} + \integral{\dimred{\Gamma}_\subwall} \height\ \htc (\dimred{\temperature} - \temperature_\subwall) \dimred{\testtemp} \dx{s} = 0 \\
		&\text{for all } \dimred{\testsolution} = (\dimred{\testvelo}, \dimred{\testpres}, \dimred{\testtemp}) \in \dimred{\adjointspace}.
	\end{aligned}
	\right.
\end{equation}
Note, that we use $\dimred{\temperature}_\subin = \temperature_\subin$ for the inflow temperature as this is assumed to be constant. For the inflow velocity we choose a parabolic profile like in \eqref{eq:parabolic_profile_ansatz} so that we are compatible with both the dimension reduction and the no-slip boundary condition on the top and bottom of the geometry.

\subsection{Application to the Shape Optimization Problem}

We now describe the application of the dimension reduction technique to the shape optimization problem \eqref{eq:def_opt_problem}. For this, we again use the parabolic profile \eqref{eq:parabolic_profile_ansatz} for the velocity as well as the constant profile \eqref{eq:constant_profile_ansatz} for pressure and temperature. Note, that we also use the parabolic profile for $\velocity_\subdes$ so that this state is reachable by the dimension reduction model. Substituting these into the cost functional \eqref{eq:def_cost_function} and carrying out the integration over $z$ yields the dimension-reduced cost functional
\begin{equation*}
	\dimred{\costfunction}(\dimred{\Omega}, \dimred{\solution}) = \weighttemp \Big( \dimred{\flux}(\dimred{\Omega}, \dimred{\solution}) - \fluxdes \Big)^2 + \weightvelo \integral{\dimred{\Omega}_\subchannels} \frac{8\height}{15} \abs{\dimred{\velocity} - \dimred{\velocity}_\subdes}^2 \dx{x} + \weightreg \left( \integral{\dimred{\Gamma}} \height \dx{s} + \integral{\dimred{\Omega}} 2 \dx{x} \right)
\end{equation*}
where $\dimred{\solution} = (\dimred{\velocity}, \dimred{\pressure}, \dimred{\temperature})$ and $\dimred{\flux}$ is defined as
\begin{equation*}
\dimred{\flux}(\dimred{\Omega}, \dimred{\solution}) = \integral{\dimred{\Gamma}_\subwall} \height\ \htc (\temperature_\subwall - \dimred{\temperature}) \dx{s} + \integral{\dimred{\Omega}} 2\htc (\temperature_\subwall - \dimred{\temperature} ) \dx{x}.
\end{equation*}
Note, that the terms $\integral{\dimred{\Omega}} 2 \dx{x}$ and $\integral{\dimred{\Omega}} 2\htc (\temperature_\subwall - \dimred{\temperature} ) \dx{x}$ arise due to the integration over lid and bottom of the three-dimensional domain. 
The corresponding reduced optimization problem reads
\begin{equation*}
	\min_{\dimred{\Omega}}\ \dimred{\reducedcostfunction}(\dimred{\Omega}) = \dimred{\costfunction}(\dimred{\Omega}, \dimred{\solution}(\dimred{\Omega})),
\end{equation*}
where $\dimred{\solution}(\dimred{\Omega})$ is the solution of \eqref{eq:weak_state_dimension_reduction} on $\dimred{\Omega}$.

To derive the shape derivative for this model, we have to make an additional assumption on the vector fields we use to deform the domain. As stated in Section~\ref{sec:optimization_problem}, we want to keep the height of the cooling system fixed. Therefore, we restrict ourselves to deformations generated by vector fields whose $z$-component vanishes, i.e.,
\begin{equation}
	\label{eq:ansatz_deformation}
	\vectorfield(x,y,z) = \begin{bmatrix}
	\dimred{\vectorfield}_1(x,y) & \dimred{\vectorfield}_2(x,y) & 0
	\end{bmatrix}^T
	= \begin{bmatrix}
	\dimred{\vectorfield}(x,y) & 0
	\end{bmatrix}^T,
\end{equation}
where $\dimred{\vectorfield} \in C^k_0(D;\R^{2})$. Any other vector field with non-vanishing $z$-component only introduces a reparametrization along the $z$-axis that cannot change the geometry as $\height$ is fixed. Using this, we apply the dimension reduction technique to the shape Lagrangian defined in \eqref{eq:shape_lagrangian} and, after lengthy calculations, end up with the shape derivative
\begin{equation}
	\label{eq:shape_derivate_dr}
	\left\lbrace \quad
	\begin{aligned}
		&d\dimred{\reducedcostfunction}(\dimred{\Omega})[\dimred{\vectorfield}] \\
		=\ &2 \weighttemp \Big( \dimred{\flux}(\dimred{\Omega}, \dimred{\solution}(\dimred{\Omega})) - \fluxdes \Big) \left( \integral{\dimred{\Gamma}_\subwall} \height\ \htc (\temperature_\subwall - \dimred{\temperature})\ \tdiv{\dimred{\vectorfield}} \!\!\! \dx{s} + \integral{\dimred{\Omega}} 2\htc (\temperature_\subwall - \dimred{\temperature})\ \divergence{\dimred{\vectorfield}} \!\!\! \dx{x} \right) \\
		%
		%
		&+ \weightvelo \integral{\dimred{\Omega}_\subchannels} \frac{8 \height}{15}\ \abs{\dimred{\velocity} - \dimred{\velocity}_\subdes}^2 \divergence{\dimred{\vectorfield}} - \frac{16 \height}{15}\ (\dimred{\velocity} - \dimred{\velocity}_\subdes) \cdot D\dimred{\velocity}_\subdes\ \dimred{\vectorfield} \dx{x} + \weightreg \integral{\dimred{\Gamma}} \height\ \tdiv{\dimred{\vectorfield}} \dx{s} \\
		&+ \weightreg \integral{\dimred{\Omega}} 2\ \divergence{\dimred{\vectorfield}} \dx{x} + \integral{\dimred{\Omega}} \frac{8\height}{15}\ \viscosity \left( D\dimred{\velocity} (\divergence{\dimred{\vectorfield}} I - 2\varepsilon(\dimred{\vectorfield})) \right) : D\dimred{\advelo} + \frac{16}{3 \height}\ \viscosity \dimred{\velocity} \cdot \dimred{\advelo}\ \divergence{\dimred{\vectorfield}} \dx{x} \\
		&- \integral{\dimred{\Omega}} \frac{2\height}{3}\ \dimred{\pressure}\ \tr \left( D\dimred{\advelo} \left( \divergence{\dimred{\vectorfield}}I - D\dimred{\vectorfield} \right) \right) \dx{x} - \integral{\dimred{\Omega}} \frac{2 \height}{3}\ \dimred{\adpres}\ \tr \left( D\dimred{\velocity} \left( \divergence{\dimred{\vectorfield}}I - D\dimred{\vectorfield} \right) \right) \dx{x} \\
		&+ \integral{\dimred{\Omega}} \height\ \conductivity \left( \left( \divergence{\dimred{\vectorfield}}I - 2\varepsilon(\dimred{\vectorfield}) \right) \grad \dimred{\temperature} \right) \cdot \grad \dimred{\adtemp} \dx{x} + \integral{\dimred{\Omega}} \frac{2\height}{3}\ \density \hcapacity \dimred{\velocity} \cdot \left( \left( \divergence{\dimred{\vectorfield}}I - D\dimred{\vectorfield}^T \right) \grad \dimred{\temperature} \right) \dimred{\adtemp} \dx{x} \\
		&+ \integral{\dimred{\Gamma}_\subwall} \height\ \htc \left( \dimred{\temperature} - \temperature_\subwall \right) \dimred{\adtemp}\ \tdiv{\dimred{\vectorfield}} \dx{s} + \integral{\dimred{\Omega}} 2\htc \left( \dimred{\temperature} - \temperature_\subwall \right) \dimred{\adtemp}\ \divergence{\dimred{\vectorfield}} \dx{x},
	\end{aligned}
	\right.
\end{equation}
where $(\dimred{\velocity}, \dimred{\pressure}, \dimred{\temperature}) = \dimred{\solution}(\dimred{\Omega})$ is the weak solution of \eqref{eq:weak_state_dimension_reduction} and $(\dimred{\advelo}, \dimred{\adpres}, \dimred{\adtemp}) = \dimred{P}(\dimred{\Omega})$ solves the adjoint system
\begin{equation*}
	\left\lbrace \quad
	\begin{aligned}
		&\text{Find } \dimred{\adsolution} = (\dimred{\advelo}, \dimred{\adpres}, \dimred{\adtemp}) \in \dimred{\adjointspace} \text{ such that}\\
		&\qquad\ \ \integral{\dimred{\Omega}} \height\ \conductivity\ \grad \dimred{\adtemp} \cdot \grad \dimred{\testadtemp} + \frac{2\height}{3} \density\hcapacity\ \dimred{\velocity} \cdot \grad \dimred{\testadtemp} \dimred{\adtemp} + 2\htc \dimred{\adtemp} \dimred{\testadtemp} \dx{x} + \integral{\dimred{\Gamma}_\subwall} \height\ \htc \dimred{\adtemp} \dimred{\testadtemp} \dx{s}\\
		&\qquad\ \ + \integral{\dimred{\Omega}} \frac{8\height}{15}\viscosity D \dimred{\advelo} : D \dimred{\testadvelo} + \frac{16}{3\height} \viscosity \dimred{\advelo} \cdot \dimred{\testadvelo} - \frac{2\height}{3} \dimred{\adpres}\ \divergence{\dimred{\testadvelo}} - \frac{2\height}{3}\dimred{\testadpres}\ \divergence{\dimred{\advelo}} + \frac{2\height}{3}\density\hcapacity\ \dimred{\testadvelo} \cdot \grad \dimred{\temperature} \dimred{\adtemp} \dx{x} \\
		&\qquad = 2\weighttemp \Big( \dimred{\flux}(\dimred{\Omega}, \dimred{\solution}(\dimred{\Omega})) - \fluxdes \Big) \left( \integral{\dimred{\Gamma}_\subwall} \height\ \htc \dimred{\testadtemp} \dx{s} + \integral{\dimred{\Omega}} 2\htc \dimred{\testadtemp} \dx{x} \right) -\frac{16 \height}{15}\weightvelo \integral{\dimred{\Omega}_\subchannels} (\dimred{\velocity} - \dimred{\velocity}_\subdes) \cdot \dimred{\testadvelo} \dx{x}\\
		&\text{for all } \dimred{\testsolution} = (\dimred{\testadvelo}, \dimred{\testadpres}, \dimred{\testadtemp}) \in \dimred{\adjointspace}.
	\end{aligned}
	\right.
\end{equation*}

\begin{Remark*}
	Note, that the only simplification made for the dimension-reduction model is that we suppose that the state variables are of the form \eqref{eq:parabolic_profile_ansatz} and \eqref{eq:constant_profile_ansatz}. Hence, we could derive the shape derivative given in \eqref{eq:shape_derivate_dr} also by applying the dimension-reduction technique directly to the adjoint system \eqref{eq:weak_adjoint} as well as to the shape derivative \eqref{eq:shape_derivative}. Thus, the optimization commutes with the dimension reduction technique.
\end{Remark*}

\section{Porous Medium Model}
\label{sec:porous_medium}

Both previously described models completely resolve all microchannels which complicates a physically meaningful discretization of the geometry. To circumvent this, we now introduce our second reduced model, for which we model the microchannels as a porous medium. This allows us to use a much simpler geometry which substantially eases the numerical solution of the equations. For a detailed introduction to porous medium models and their application to microchannels we refer, e.g., to \cite{nield_bejan, kaviany} and \cite{chen2007forced, kim1999forced, kim2000local}, respectively.

\subsection{Description of the Geometry for the Darcy Model}

\begin{figure}[b]
	\centering
	\begin{tikzpicture}
		\node at (0,0) {\includegraphics[width=0.55\textwidth, trim=11cm 0cm 0cm 0cm, clip]{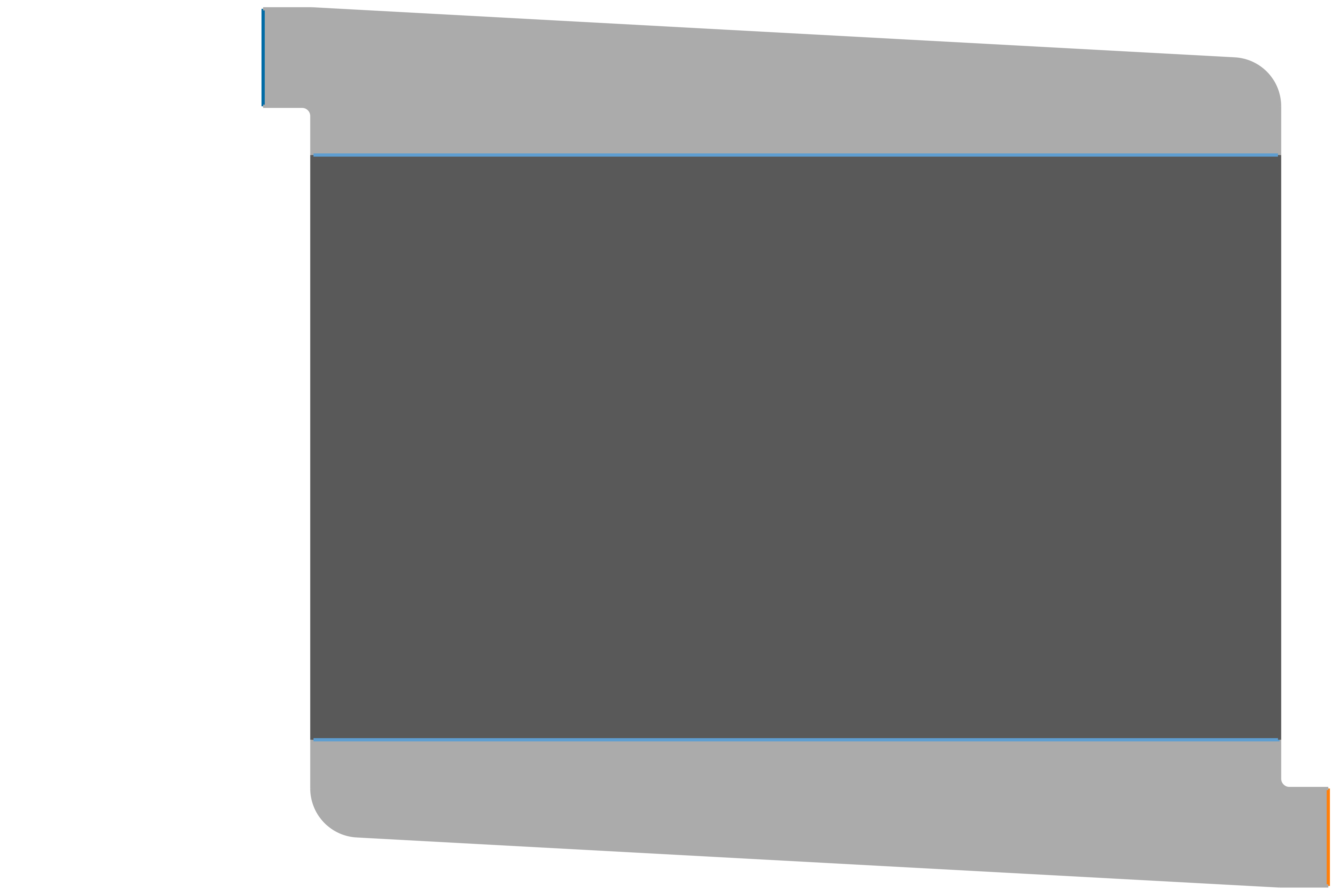}};
		\node at (0,0) {\LARGE \color{white} $\porous{\dimred{\Omega}}_\subdarcy$};
		\node at (0,2.75) {\LARGE $\porous{\dimred{\Omega}}_\subfluid$};
		\node at (0, -2.75) {\LARGE $\porous{\dimred{\Omega}}_\subfluid$};
		\node at (-5.5,3.25) {\LARGE \color{tabblue} $\porous{\dimred{\Gamma}}_\subin$};
		\draw[tabblue, ->, line width=1]  (-5,3.25) -- (-4.5,3.25);
		\node at (5.4, -3.25) {\LARGE \color{taborange} $\porous{\dimred{\Gamma}}_\subout$};
		\draw[taborange, ->, line width=1] (4.4, -3.25) -- (4.9, -3.25);
		\node at (5.125, 3.25) {\LARGE \color{tablightgray} $\porous{\dimred{\Gamma}}_\subwall$};
		\draw[tablightgray, ->, line width=1] (4.5, 3.25) -- (4.0, 3.15);
		\node at (-5, 0) {\LARGE \color{tabdarkgray} $\porous{\dimred{\Gamma}}_\subdarcy$};
		\draw[tabdarkgray, ->, line width=1] (-4.5, 0) -- (-4, 0);
		\node at (2.5, 0) {\LARGE \color{tablightblue} $\porous{\dimred{\Gamma}}_\subfd$};
		\draw[tablightblue, ->, line width=1] (2.5, 0.5) -- (2.5, 2.25); 
		\draw[tablightblue, ->, line width=1] (2.5, -0.5) -- (2.5, -2.25); 
	\end{tikzpicture}
	\caption{Two-dimensional geometry for the Darcy model $\porous{\dimred{\Omega}}$, partitioned into the fluid part $\porous{\dimred{\Omega}}_\subfluid$ (light gray) and porous medium part $\porous{\dimred{\Omega}}_\subdarcy$ (dark gray).}
	\label{figure:geometry_averaged}
\end{figure}


First, we introduce the structure of the domain $\porous{\Omega}$ that we consider for the Darcy model. It consists of the two subdomains $\porous{\Omega}_\subfluid = \Omega \setminus \Omega_\subchannels$, the domain without the microchannels that remains unaltered, and $\porous{\Omega}_\subdarcy$, which is the homogenized geometry corresponding to $\Omega_\subchannels$. The boundary of $\porous{\Omega}$ is denoted by $\porous{\Gamma}$ and consists of the following parts: For the in- and outlet we have $\porous{\Gamma}_\subin = \Gamma_\subin$ and $\porous{\Gamma}_\subout = \Gamma_\subout$, and the wall boundary is given as $\porous{\Gamma}_\subwall = \Gamma_\subwall \setminus \Gamma_\subchannels$. Additionally, we obtain a new boundary $\porous{\Gamma}_\subdarcy$ that is the outer boundary of $\porous{\Omega}_\subdarcy$, and the interfaces between $\porous{\Omega}_\subfluid$ and $\porous{\Omega}_\subdarcy$ are denoted by $\porous{\Gamma}_\subfd$. Naturally, the geometry of the Darcy model has a similar structure to the geometry of the full model, i.e., we have $\porous{\Omega} = \porous{\dimred{\Omega}} \times (0, \height)$. An analogous decomposition also holds for the two subdomains as well as the boundaries (cf. Section~\ref{sec:description_of_the_geometry}). The corresponding two-dimensional domain $\porous{\dimred{\Omega}}$ depicting the above situation can be seen in Figure~\ref{figure:geometry_averaged}.

\subsection{Description of the Model}

To couple the fluid equations on $\porous{\Omega}_\subfluid$, where we have the usual Stokes system as in Section~\ref{sec:problem_formulation}, and $\porous{\Omega}_\subdarcy$, where we consider the porous medium, we choose the Brinkman equation (cf. \cite{nield_bejan}). This allows an implicit coupling through transmission conditions that vanish in the weak formulation of the problem (cf. \eqref{eq:weak_darcy_state}). In particular, the Brinkman equation, without the transmission conditions, reads
\begin{equation}
	\label{eq:brinkman}
	\left\lbrace \quad
	\begin{alignedat}{2}
		- \viscosity' \laplace \velocity + \grad \pressure + \viscosity \permeability^{-1} \velocity &= 0 \quad &&\text{ in } \porous{\Omega}_\subdarcy,\\
		\divergence{\velocity} &= 0 \quad &&\text{ in } \porous{\Omega}_\subdarcy,\\
		\velocity \cdot \normal &= 0 \quad &&\text{ on } \porous{\Gamma}_\subdarcy,\\
		\viscosity' \partial_\normal \velocity \times \normal &= 0 \quad &&\text{ on } \porous{\Gamma}_\subdarcy,
	\end{alignedat}
	\right.
\end{equation}
where $\permeability$ denotes the permeability tensor and $\viscosity'$ denotes the effective viscosity. Here, $\velocity$ describes the averaged (or Darcy) velocity and not a physical fluid velocity. The permeability tensor $\permeability$ already models the friction arising from the no-slip boundary condition in the microchannels completely so that we must not prescribe a no-slip condition on $\porous{\Gamma}_\subdarcy$ for the Brinkman equation. Instead we choose the slip condition $\velocity \cdot \normal = 0$ and $\viscosity' \partial_\normal \velocity \times \normal = 0$ on $\porous{\Gamma}_\subdarcy$ so that the fluid cannot leave the domain and does not experience additional friction. For the permeability tensor we observe that the flow through the channels is basically one-dimensional and in $y$-direction. To model this anisotropy, $\permeability$ has the structure $\hat{K} \cdot \mathrm{diag}(\varepsilon, 1, \varepsilon)$, where $\hat{K}$ is the permeability in $y$-direction and $1 \gg \varepsilon > 0$ is a relaxation parameter needed for the invertibility of $\permeability$ (cf. \eqref{eq:brinkman}). We compute $\hat{K}$ numerically using Darcy's law (see, e.g., \cite{nield_bejan, kaviany}) and choose $\varepsilon = \num{1e-5}$. Further, as $\nicefrac{1}{\hat{\permeability}}$ is very large, we neglect the influence of the effective viscosity and set $\viscosity' = \viscosity$ as originally suggested in \cite{brinkman}.


To describe the transmission conditions for the coupling of \eqref{eq:brinkman} with the Stokes system, we denote by $v^\subfluid$ and $v^\subdarcy$ the restriction of a function $v$ on $\porous{\Omega}_\subfluid$ and $\porous{\Omega}_\subdarcy$, respectively. Then, the coupling conditions are given by the continuity of both the velocity and the normal component of the viscous stress tensor, i.e.,
\begin{equation}
	\label{eq:transmission_fluid}
	\velocity^\subfluid - \velocity^\subdarcy = 0 \text{ on } \porous{\Gamma}_\subfd \qquad \text{ and } \qquad \left(\viscosity \partial_{\hat{\normal}} \velocity^\subfluid - \pressure^\subfluid \hat{\normal}\right) - \left(\viscosity \partial_{\hat{\normal}} \velocity^\subdarcy - \pressure^\subdarcy \hat{\normal}\right) = 0 \text{ on } \porous{\Gamma}_\subfd,
\end{equation}
where $\hat{\normal}$ denotes the outer unit normal to $\porous{\Omega}_\subfluid$.
\begin{table}[b]
	\centering
	{\footnotesize
	\rowcolors{2}{\tablegray}{white}
	\begin{tabular}{l l}
		\toprule
		parameter [unit] & value\\
		\midrule
		porosity $\volumefrac$ [\si{1}] & \num{2.02e-1} \\
		permeability $\hat{\permeability}$ [\si{\cubic \meter}] & \num{3.16e-9} \\
		interfacial heat transfer coefficient $\hfs$ [\si{\watt \per \kelvin \per \cubic \meter}] & \num{2.63e+4} \\
		\bottomrule
	\end{tabular}
	\caption{Parameters for the Darcy model.}
	\label{table:porous_parameters}
	}
\end{table}
To summarize, the Darcy model for the fluid reads
\begin{equation}
\label{eq:darcy_stokes}
\left\lbrace \quad
\begin{alignedat}{2}
-\viscosity \laplace \velocity + \grad \pressure &= 0 \quad &&\text{ in } \porous{\Omega}_\subfluid,\\
-\viscosity \laplace \velocity + \grad \pressure + \viscosity \permeability^{-1} \velocity &= 0 \quad &&\text{ in } \porous{\Omega}_\subdarcy,\\
\divergence{\velocity} &= 0 \quad &&\text{ in } \porous{\Omega},\\
\velocity &= \velocity_\subin \quad &&\text{ on } \porous{\Gamma}_\subin,\\
\velocity &= 0 \quad &&\text{ on } \porous{\Gamma}_\subwall,\\
\velocity \cdot \normal &= 0 \quad &&\text{ on } \porous{\Gamma}_\subdarcy,\\
\viscosity \partial_\normal \velocity \times \normal &= 0 \quad &&\text{ on } \porous{\Gamma}_\subdarcy, \\
\viscosity \partial_\normal \velocity - \pressure\normal &= 0 \quad &&\text{ on } \porous{\Gamma}_\subout, \\
\velocity^\subfluid - \velocity^\subdarcy &= 0 \quad &&\text{ on } \porous{\Gamma}_\subfd, \\
\left(\viscosity \partial_{\hat{\normal}} \velocity^\subfluid - \pressure^\subfluid \hat{\normal}\right) - \left(\viscosity \partial_{\hat{\normal}} \velocity^\subdarcy - \pressure^\subdarcy \hat{\normal}\right) &= 0 \quad &&\text{ on } \porous{\Gamma}_\subfd.
\end{alignedat}
\right.
\end{equation}

For the modeling of heat transfer in the porous medium we modify the so-called local thermal non-equilibrium approach that can be found, e.g., in \cite[Chapter 2.2]{nield_bejan}. This model uses two equations, one for the liquid and one for the solid phase. However, we assume that the solid's temperature is constant, in analogy to the constant wall temperature $\temperature_\subwall$ (cf. Section~\ref{sec:problem_formulation}), and are left with the equation for the fluid phase. In particular, the heat transfer equation in $\porous{\Omega}_\subdarcy$ reads
\begin{equation}
	\label{eq:ltne}
	\left\lbrace \quad
	\begin{alignedat}{2}
		- \grad \cdot (\volumefrac \conductivity \grad \temperature) + \density \hcapacity\ \velocity \cdot \grad \temperature + \hfs (\temperature - \temperature_\subwall) &= 0 \quad &&\text{ in } \porous{\Omega}_\subdarcy,\\
		\volumefrac \conductivity\ \partial_\normal \temperature &= 0 \quad &&\text{ on } \porous{\Gamma}_\subdarcy,
	\end{alignedat}
	\right.
\end{equation}
where $\volumefrac$ is the porosity of the porous medium, i.e., the fraction of total volume occupied by the fluid, and $\hfs$ denotes the interfacial heat transfer coefficient. The latter can, again, be computed numerically using the formulas given in, e.g., \cite{ihtc1, ihtc2}.
The homogeneous Neumann boundary condition in \eqref{eq:ltne} is used for a similar reason as the slip boundary condition in \eqref{eq:brinkman}. The heat transfer between solid and fluid phase is completely contained in the term $\hfs (\temperature - \temperature_\subwall)$ and, therefore, we must not have an additional heat source on the boundary $\porous{\Gamma}_\subdarcy$. The parameters for the porous medium model which we computed numerically using a single microchannel as representative elementary volume are given in Table~\ref{table:porous_parameters}.

Finally, for the transmission conditions coupling the temperature in $\porous{\Omega}_\subdarcy$ with the one in $\porous{\Omega}_\subfluid$ we proceed analogously and require the continuity of the temperature and the normal component of the heat flux over the interface $\porous{\Gamma}_\subfd$, i.e.,
\begin{equation*}
	\temperature^\subfluid - \temperature^\subdarcy = 0 \text{ on } \porous{\Gamma}_\subfd \qquad \text{ and } \qquad \conductivity\ \partial_{\hat{\normal}} \temperature^\subfluid - \volumefrac \conductivity\ \partial_{\hat{\normal}} \temperature^\subdarcy = 0 \quad \text{ on } \porous{\Gamma}_\subfd.
\end{equation*}
Note, that demanding the continuity of the heat flux' conductive part over the interface is already sufficient for the continuity of the entire flux, as we get the continuity of its convective part directly from \eqref{eq:transmission_fluid} and the continuity of the temperature. To summarize, our model for heat transfer in the porous medium reads
\begin{equation*}
	\left\lbrace \quad
	\begin{alignedat}{2}
	- \grad \cdot \left( \conductivity \grad \temperature \right) + \density \hcapacity\ \velocity \cdot \grad \temperature &= 0 \quad &&\text{ in } \porous{\Omega}_\subfluid,\\
	- \grad \cdot \left( \volumefrac\conductivity \grad \temperature \right) + \density \hcapacity\ \velocity \cdot \grad \temperature + \hfs (\temperature - \temperature_\subwall) &= 0 \quad &&\text{ in } \porous{\Omega}_\subdarcy, \\
	\temperature &= \temperature_\subin \quad &&\text{ on } \porous{\Gamma}_\subin,\\
	\conductivity \partial_\normal \temperature + \htc \left( \temperature - \temperature_\subwall \right) &= 0 \quad &&\text{ on } \porous{\Gamma}_\subwall,\\
	\conductivity \volumefrac\ \partial_\normal \temperature &= 0 \quad &&\text{ on } \porous{\Gamma}_\subdarcy,\\
	\conductivity \partial_\normal \temperature &= 0 \quad &&\text{ on } \porous{\Gamma}_\subout,\\
	\temperature^\subfluid - \temperature^\subdarcy &= 0 \quad &&\text{ on } \porous{\Gamma}_\subfd,\\
	\conductivity\ \partial_{\hat{\normal}} \temperature^\subfluid - \volumefrac\conductivity\ \partial_{\hat{\normal}} \temperature^\subdarcy &= 0 \quad &&\text{ on } \porous{\Gamma}_\subfd,
	\end{alignedat}
	\right.
\end{equation*}
where $\velocity$ solves \eqref{eq:darcy_stokes}. 
For the weak form of the Darcy model we introduce the (affine) spaces
\begin{equation*}
	\left\lbrace \quad
	\begin{aligned}
		\porous{\fspacetest} &:= \Set{\testvelo \in H^1(\porous{\Omega})^3 | \testvelo=0 \text{ on } \porous{\Gamma}_\subin \cup \porous{\Gamma}_\subwall \text{ and } \testvelo \cdot \normal = 0 \text{ on } \porous{\Gamma}_\subdarcy}, \quad \porous{\fspacetrial} := \velocity_\subin + \porous{\fspacetest}, \\
		\porous{\pspace} &:= L^2(\porous{\Omega}), \qquad \porous{\tspacetest} := \Set{\testtemp \in H^1(\porous{\Omega}) | \testtemp=0 \text{ on } \porous{\Gamma}_\subin}, \qquad \porous{\tspacetrial} := \temperature_\subin + \porous{\tspacetest}, \\
		\porous{\statespace} &:= \porous{\fspacetrial} \times \porous{\pspace} \times \porous{\tspacetrial}, \qquad \porous{\adjointspace} = \porous{\fspacetest} \times \porous{\pspace} \times \porous{\tspacetest},
	\end{aligned}
	\right. 
\end{equation*}
where we assume that $\velocity_\subin \in H^1(\porous{\Omega})^3$ with $\velocity_\subin = 0$ on $\porous{\Gamma}_\subwall$ and $\velocity_\subin \cdot \normal = 0$ on $\porous{\Gamma}_\subdarcy$ as well as $\temperature_\subin \in H^1(\porous{\Omega})$ in analogy to Sections~\ref{sec:problem_formulation} and \ref{sec:dimension_reduction}. Finally, the weak form of the Darcy model is given by
\begin{equation}
	\label{eq:weak_darcy_state}
	\left\lbrace \quad 
	\begin{aligned}
		&\text{Find } \solution = (\velocity, \pressure, \temperature) \in \porous{\statespace} \text{ such that } \\
		&\qquad \integral{\porous{\Omega}} \viscosity D \velocity : D \testvelo - \pressure\ \divergence{\testvelo} - \testpres\ \divergence{\velocity} \dx{x} + \integral{\porous{\Omega}_\subdarcy} \viscosity \permeability^{-1} \velocity \cdot \testvelo \dx{x} \\
		&\qquad + \integral{\porous{\Omega}_\subfluid} \conductivity \grad \temperature \cdot \grad \testtemp \dx{x} + \integral{\porous{\Omega}_\subdarcy} \volumefrac\conductivity \grad \temperature \cdot \grad \testtemp \dx{x} + \integral{\porous{\Omega}} \density \hcapacity\ \velocity\cdot \grad \temperature\ \testtemp \dx{x} \\
		&\qquad + \integral{\porous{\Omega}_\subdarcy} \hfs(\temperature - \temperature_\subwall) \testtemp \dx{x} + \integral{\porous{\Gamma}_\subwall} \htc (\temperature - \temperature_\subwall) \testtemp \dx{s} = 0 \\
		&\text{for all } \testsolution = (\testvelo, \testpres, \testtemp) \in \porous{\adjointspace}.
	\end{aligned}
	\right.
\end{equation}

\subsection{Application to the Shape Optimization Problem}

While the problem for $\porous{\Omega}_\subfluid$ stays like in Section~\ref{sec:problem_formulation}, we now have to modify the cost functional in the porous region $\porous{\Omega}_\subdarcy$ to reflect the changes made in modeling the channels as porous medium. Therefore, the cost functional \eqref{eq:def_cost_function} takes the following form for the Darcy model. The term $\costfunction_1$ can be expressed by
\begin{equation*}
\porous{\costfunction}_1(\porous{\Omega}, \solution) = \Big(\porous{\flux}(\porous{\Omega}, \solution) - \fluxdes \Big)^2,
\end{equation*}
where
\begin{equation*}
	\porous{\flux}(\porous{\Omega}, \solution) = \integral{\porous{\Gamma}_\subwall} \htc (\temperature_\subwall - \temperature) \dx{s} + \integral{\porous{\Omega}_\subdarcy} \hfs (\temperature_\subwall - \temperature) \dx{x}.
\end{equation*}
This is due to the fact that the heat transfer in $\porous{\Omega}_\subdarcy$ is given by the volume source term $\hfs (\temperature - \temperature_\subwall)$ for the Darcy model. Note, that we again write $\solution = (\velocity, \pressure, \temperature)$. Further, $\costfunction_2$ can be represented by the tracking type functional
\begin{equation*}
	\porous{\costfunction}_2(\porous{\Omega}, \solution) = \integral{\porous{\Omega}_\subdarcy} \abs{\velocity - \porous{\velocity}_\subdes}^2 \dx{x}.
\end{equation*}
The difference now lies in the interpretation of the velocity. Whereas we have a physical velocity in $\Omega_\subchannels$ and, thus, also prescribe a physical one for $\velocity_\subdes$, the fluid velocity in $\porous{\Omega}_\subdarcy$ is an averaged one. Hence, the desired velocity $\porous{\velocity}_\subdes$ has to be the mean of $\velocity_\subdes$ over $\porous{\Omega}_\subdarcy$. The perimeter regularization is modeled analogously to before, i.e., we use
\begin{equation*}
	\porous{\costfunction}_3(\porous{\Omega}, \solution) = \integral{\porous{\Gamma}} 1 \dx{s}.
\end{equation*}
Finally, the cost functional $\porous{\costfunction}$ for the Darcy model reads
\begin{equation*}
	\porous{\costfunction}(\porous{\Omega}, \solution) = \weighttemp \porous{\costfunction}_1(\porous{\Omega}, \solution) + \weightvelo \porous{\costfunction}_2(\porous{\Omega}, \solution) + \weightreg \porous{\costfunction}_3(\porous{\Omega}, \solution).
\end{equation*}

As before, we define the reduced cost functional by $\porous{\reducedcostfunction}(\porous{\Omega}) = \porous{\costfunction}(\porous{\Omega}, \solution(\porous{\Omega}))$, where $\solution(\porous{\Omega})$ denotes the solution of \eqref{eq:weak_darcy_state} on $\porous{\Omega}$ and the reduced optimization problem is given by
\begin{equation*}
	\min_{\porous{\Omega}}\ \porous{\reducedcostfunction}(\porous{\Omega}).
\end{equation*}
Applying the same techniques as in Section~\ref{sec:shape_calculus} we derive the following shape derivative for $\porous{\reducedcostfunction}$
\begin{equation*}
	\left\lbrace \quad
	\begin{aligned}
		&d\porous{\reducedcostfunction}(\porous{\Omega})[\vectorfield] \\
		=\ &2 \weighttemp \Big( \porous{\flux}(\porous{\Omega}, \solution(\porous{\Omega})) - \fluxdes \Big) \integral{\porous{\Gamma}_\subwall} \htc (\temperature_\subwall - \temperature) \tdiv{\vectorfield} \dx{s} \\ 
		&+ 2 \weighttemp \Big( \porous{\flux}(\porous{\Omega}, \solution(\porous{\Omega})) - \fluxdes \Big) \integral{\porous{\Omega}_\subdarcy} \hfs (\temperature_\subwall - \temperature) \divergence{\vectorfield} \dx{x} \\
		%
		%
		&+ \weightvelo \integral{\porous{\Omega}_\subdarcy} \abs{\velocity - \porous{\velocity}_\subdes}^2 \divergence{\vectorfield} \dx{x} + \weightreg \integral{\porous{\Gamma}} \tdiv{\vectorfield} \dx{s} \\
		&+ \integral{\porous{\Omega}} \viscosity \left( D\velocity \left( \divergence{\vectorfield}I - 2\varepsilon(\vectorfield) \right) \right) : D\advelo \dx{x} - \integral{\porous{\Omega}} \pressure\ \tr \left( D\advelo \left( \divergence{\vectorfield}I - D\vectorfield \right) \right) \dx{x} \\
		&- \integral{\porous{\Omega}} \adpres\ \tr \left( D\velocity \left( \divergence{\vectorfield}I - D\vectorfield \right) \right) \dx{x} + \integral{\porous{\Omega}_\subdarcy} \viscosity \permeability^{-1} \velocity \cdot \advelo\ \divergence{\vectorfield} \dx{x} \\
		&+ \integral{\porous{\Omega}_\subfluid} \conductivity \left( \left( \divergence{\vectorfield}I - 2\varepsilon(\vectorfield) \right) \grad \temperature \right) \cdot \grad \adtemp  \dx{x} + \integral{\porous{\Omega}_\subdarcy} \volumefrac \conductivity \left( \left( \divergence{\vectorfield}I - 2\varepsilon(\vectorfield)  \right) \grad \temperature \right) \cdot \grad \adtemp \dx{x} \\
		&+ \integral{\porous{\Omega}} \density \hcapacity\ \velocity \cdot \left( \left( \divergence{\vectorfield}I - D\vectorfield^T \right) \grad \temperature \right) \adtemp \dx{x} + \integral{\porous{\Gamma}_\subwall} \htc \left( \temperature - \temperature_\subwall \right) \adtemp\ \tdiv{\vectorfield} \dx{s} \\
		&+ \integral{\porous{\Omega}_\subdarcy} \hfs \left( \temperature - \temperature_\subwall \right) \adtemp\ \divergence{\vectorfield} \dx{x},
	\end{aligned}
	\right.
\end{equation*}
where $(\velocity, \pressure, \temperature) = \solution(\porous{\Omega})$ is the solution of \eqref{eq:weak_darcy_state} and $(\advelo, \adpres, \adtemp) = \adsolution(\porous{\Omega})$ is the solution of the following adjoint system
\begin{equation*}
	\left\lbrace \quad
	\begin{aligned}
		&\text{Find } \adsolution = (\advelo, \adpres, \adtemp) \in \porous{\adjointspace} \text{ such that } \\
		&\qquad \integral{\porous{\Omega}_\subfluid} \conductivity \grad \adtemp \cdot \grad \testadtemp \dx{x} + \integral{\porous{\Omega}_\subdarcy} \volumefrac \conductivity\ \grad \adtemp \cdot \grad \testadtemp + \hfs \adtemp \testadtemp \dx{x} + \integral{\porous{\Omega}} \density \hcapacity\ \velocity \cdot \grad \testadtemp \adtemp \dx{x} + \integral{\porous{\Gamma}_\subwall} \hspace{-0.5em} \htc \adtemp \testadtemp \dx{s} \\
		&\qquad\ +\integral{\porous{\Omega}} \viscosity\ D \advelo : D \testadvelo - \adpres\ \divergence{\testadvelo} - \testadpres\ \divergence{\advelo} + \density \hcapacity \testadvelo \cdot \grad \temperature \adtemp \dx{x} + \integral{\porous{\Omega}_\subdarcy} \viscosity \permeability^{-1} \advelo \cdot \testadvelo \dx{x} \\
		&\qquad = 2\weighttemp \Big( \porous{\flux}(\porous{\Omega}, \solution(\porous{\Omega})) - \fluxdes \Big) \left( \integral{\porous{\Gamma_\subwall}} \htc \testadtemp \dx{s} + \integral{\porous{\Omega}_\subdarcy} \hfs \testadtemp \dx{x} \right)\\
		&\qquad\ \ -2 \weightvelo \integral{\porous{\Omega}_\subdarcy} (\velocity - \porous{\velocity}_\subdes) \cdot \testadvelo \dx{x} \\
		&\text{for all } \testadsolution = (\testadvelo, \testadpres, \testadtemp) \in \porous{\adjointspace}.
	\end{aligned}
	\right.
\end{equation*}

\section{Dimension Reduction Applied to the Darcy Model}
\label{sec:darcy_2D}

For our final reduced model, we apply the dimension reduction technique of Section~\ref{sec:dimension_reduction} to the porous medium model derived in the previous section. 

\subsection{Description of the Model}

For the application of the dimension reduction approach, we now use two different profiles for the fluid velocity. In $\porous{\Omega}_\subfluid$ we have a physical velocity and, hence, we use a parabolic profile as in \eqref{eq:parabolic_profile_ansatz}. In contrast, in $\porous{\Omega}_\subdarcy$ we have an averaged velocity and, therefore, we assume to have a constant profile for the velocity there. I.e., we use
\begin{equation*}
	\velocity = \frac{6}{\height^2} z(\height-z)
	\begin{bmatrix}
		\dimred{\velocity}_1(x,y) & \dimred{\velocity}_2(x,y) & 0
	\end{bmatrix}^T \text{ in } \porous{\Omega}_\subfluid \quad \text{ and } \quad
	\velocity = 
	\begin{bmatrix}
		\dimred{\velocity}_1(x,y) & \dimred{\velocity}_2(x,y) & 0
	\end{bmatrix}^T \text{ in } \porous{\Omega}_\subdarcy.
\end{equation*}
Note, that we use the scaling factor of $\nicefrac{6}{\height^2}$ so that $\dimred{\velocity} = [\dimred{\velocity}_1\ \dimred{\velocity}_2]^T$ represents the mean fluid velocity in both parts of the domain, making a coupling with transmission conditions similarly to \eqref{eq:transmission_fluid} possible. For the pressure and temperature we use the same approach as in \eqref{eq:constant_profile_ansatz} and assume that they are constant in $z$. We introduce the (affine) Sobolev spaces
\begin{equation*}
	\left\lbrace \quad
	\begin{aligned}
		\porous{\dimred{\fspacetest}} &:= \Set{\dimred{\testvelo} \in H^1(\porous{\dimred{\Omega}})^2 | \dimred{\testvelo}=0 \text{ on } \porous{\dimred{\Gamma}}_\subin \cup \porous{\dimred{\Gamma}}_\subwall \text{ and } \dimred{\testvelo} \cdot \normal = 0 \text{ on } \porous{\dimred{\Gamma}}_\subdarcy}, \hspace{0.75em} \porous{\dimred{\fspacetrial}} := \dimred{\velocity}_\subin + \porous{\dimred{\fspacetest}}, \\
		\porous{\dimred{\pspace}} &:= L^2(\porous{\dimred{\Omega}}), \qquad \porous{\dimred{\tspacetest}} := \Set{\dimred{\testtemp} \in H^1(\porous{\dimred{\Omega}}) | \dimred{\testtemp}=0 \text{ on } \porous{\dimred{\Gamma}}_\subin}, \qquad \porous{\dimred{\tspacetrial}} := \dimred{\temperature}_\subin + \porous{\dimred{\pspace}},\\
		\porous{\dimred{\statespace}} &:= \porous{\dimred{\fspacetrial}} \times \porous{\dimred{\pspace}} \times \porous{\dimred{\tspacetrial}}, \qquad \porous{\dimred{\adjointspace}} := \porous{\dimred{\fspacetest}} \times \porous{\dimred{\pspace}} \times \porous{\dimred{\tspacetest}},
	\end{aligned}
	\right.
\end{equation*}
where $\dimred{\velocity}_\subin \in H^1(\porous{\dimred{\Omega}})^2$ with $\dimred{\velocity}_\subin = 0$ on $\porous{\dimred{\Gamma}}_\subwall$ and $\dimred{\velocity}_\subin \cdot \normal = 0$ on $\porous{\dimred{\Gamma}}_\subdarcy$ as well as $\dimred{\temperature}_\subin \in H^1(\porous{\dimred{\Omega}})$, similarly to before.
Proceeding analogously to Section~\ref{sec:dimension_reduction}, we get the following weak form for the 2D Darcy model
\begin{equation}
	\label{eq:weak_darcy_dimension_state}
	\left\lbrace \quad 
	\begin{aligned}
		&\text{Find } \dimred{\solution} = (\dimred{\velocity}, \dimred{\pressure}, \dimred{\temperature}) \in \porous{\dimred{\statespace}} \text{ such that } \\
		&\qquad \integral{\porous{\dimred{\Omega}}_\subfluid} \frac{6}{5}\height\ \viscosity D \dimred{\velocity} : D \dimred{\testvelo} \dx{x} + \integral{\porous{\dimred{\Omega}}_\subdarcy} \height\ \viscosity D \dimred{\velocity} : D \dimred{\testvelo} \dx{x} - \integral{\porous{\dimred{\Omega}}} \height\ \dimred{\pressure}\ \divergence{\dimred{\testvelo}} + \height\ \dimred{\testpres}\ \divergence{\dimred{\velocity}} \dx{x} \\
		&\quad + \integral{\porous{\dimred{\Omega}}_\subfluid} \frac{12}{\height}\viscosity\ \dimred{\velocity} \cdot \dimred{\testvelo} \dx{x} + \integral{\porous{\dimred{\Omega}}_\subdarcy} \height\ \viscosity \permeability^{-1} \dimred{\velocity} \cdot \dimred{\testvelo} \dx{x} \\
		&\quad + \integral{\porous{\dimred{\Omega}}_\subfluid} \height\ \conductivity \grad \dimred{\temperature} \cdot \grad \dimred{\testtemp} \dx{x} + \integral{\porous{\dimred{\Omega}}_\subdarcy} \height\ \volumefrac\conductivity\ \grad \dimred{\temperature} \cdot \grad \dimred{\testtemp} \dx{x} + \integral{\porous{\dimred{\Omega}}} \height\ \density \hcapacity\ \dimred{\velocity}\cdot \grad \dimred{\temperature}\ \dimred{\testtemp} \dx{x} \\
		&\quad + \integral{\porous{\dimred{\Omega}}_\subdarcy} \height\ \hfs(\dimred{\temperature} - \temperature_\subwall) \dx{x} + \integral{\porous{\dimred{\Gamma}_\subwall}} \height\ \htc (\dimred{\temperature} - \temperature_\subwall) \dimred{\testtemp} \dx{s} + \integral{\porous{\dimred{\Omega}}_\subfluid} 2 \htc (\dimred{\temperature} - \temperature_\subwall) \dx{x} = 0 \\
		&\text{for all } \dimred{\testsolution} = (\dimred{\testvelo}, \dimred{\testpres}, \dimred{\testtemp}) \in \porous{\dimred{\adjointspace}},
	\end{aligned}
	\right.
\end{equation}
where we rescaled the equations to obtain the continuity of the pressure over the interface $\porous{\dimred{\Gamma}}_\subfd$.

\subsection{Application to the Optimization Problem}

For the shape optimization problem we again proceed as before and consider the cost functional
\begin{equation*}
	\porous{\dimred{\costfunction}}(\porous{\dimred{\Omega}}, \dimred{\solution}) = \weighttemp \Big( \porous{\dimred{\flux}}(\porous{\dimred{\Omega}}, \dimred{\solution}) - \fluxdes \Big)^2 + \weightvelo \integral{\porous{\dimred{\Omega}}_\subdarcy} \hspace{-0.5em} \height\ \abs{\dimred{\velocity} - \porous{\velocity}_\subdes}^2 \dx{x} + \weightreg \left(  \integral{\porous{\dimred{\Gamma}}} \hspace{-0.5em} \height \dx{s} + \integral{\porous{\dimred{\Omega}}} \hspace{-0.5em} 2 \dx{x} \right),
\end{equation*}
where 
\begin{equation*}
	\porous{\dimred{\flux}}(\porous{\dimred{\Omega}}, \dimred{\solution}) = \integral{\porous{\dimred{\Gamma}}_\subwall} \height\ \htc (\temperature_\subwall - \dimred{\temperature}) \dx{s} + \integral{\porous{\dimred{\Omega}}_\subfluid} 2 \htc (\temperature_\subwall - \dimred{\temperature}) \dx{x} +  \integral{\porous{\dimred{\Omega}}_\subdarcy} \height\ \hfs (\temperature_\subwall - \dimred{\temperature}) \dx{x}.
\end{equation*}
The corresponding (reduced) shape optimization problem reads
\begin{equation*}
	\min_{\porous{\dimred{\Omega}}}\ \porous{\dimred{\reducedcostfunction}}(\porous{\dimred{\Omega}}) = \porous{\dimred{\costfunction}}(\porous{\dimred{\Omega}}, \dimred{\solution}(\porous{\dimred{\Omega}})),
\end{equation*}
where $\dimred{\solution}(\porous{\dimred{\Omega}})$ denotes the solution of \eqref{eq:weak_darcy_dimension_state}. As in Section~\ref{sec:dimension_reduction}, we assume that the $z$-component of the vector field $\vectorfield$ vanishes (cf. \eqref{eq:ansatz_deformation}) in order to calculate the shape derivative, which then reads
\begin{equation*}
	\left\lbrace \quad
	\begin{aligned}
		&d\porous{\dimred{\reducedcostfunction}}(\porous{\dimred{\Omega}})[\dimred{\vectorfield}] \\
		=\ &2\weighttemp \Big( \porous{\dimred{\flux}}(\porous{\dimred{\Omega}}, \dimred{\solution}(\porous{\dimred{\Omega}})) - \fluxdes \Big) \integral{\porous{\dimred{\Gamma}_\subwall}} \height\ \htc (\temperature_\subwall - \dimred{\temperature}) \tdiv{\dimred{\vectorfield}} \dx{s} \\
		&+ 2\weighttemp \Big( \porous{\dimred{\flux}}(\porous{\dimred{\Omega}}, \dimred{\solution}(\porous{\dimred{\Omega}})) - \fluxdes \Big) \integral{\porous{\dimred{\Omega}}_\subfluid} 2\htc (\temperature_\subwall - \dimred{\temperature}) \divergence{\dimred{\vectorfield}} \dx{x} \\
		&+ 2\weighttemp \Big( \porous{\dimred{\flux}}(\porous{\dimred{\Omega}}, \dimred{\solution}(\porous{\dimred{\Omega}})) - \fluxdes \Big) \integral{\porous{\dimred{\Omega}}_\subdarcy} \height\ \hfs (\temperature_\subwall - \dimred{\temperature}) \divergence{\dimred{\vectorfield}} \dx{x} \\
		%
		%
		%
		&+ \weightvelo \integral{\porous{\dimred{\Omega}}_\subdarcy} \height\ \abs{\dimred{\velocity} - \porous{\velocity}_\subdes}^2 \divergence{\dimred{\vectorfield}} \dx{x} + \weightreg \left( \integral{\porous{\dimred{\Gamma}}} \height\ \tdiv{\dimred{\vectorfield}} \dx{s} + \integral{\porous{\dimred{\Omega}}} 2\  \divergence{\dimred{\vectorfield}} \dx{x} \right) \\
		&+ \integral{\porous{\dimred{\Omega}}_\subfluid} \frac{6}{5} \height\ \viscosity \left( D\dimred{\velocity} \left( \divergence{\dimred{\vectorfield}}I - 2\varepsilon(\dimred{\vectorfield}) \right) \right) : D\dimred{\advelo} + \frac{12}{\height}\ \viscosity\ \dimred{\velocity} \cdot \dimred{\advelo}\ \divergence{\dimred{\vectorfield}} \dx{x} \\
		&+ \integral{\porous{\dimred{\Omega}}_\subdarcy} \height\ \viscosity \left( D\dimred{\velocity} \left( \divergence{\dimred{\vectorfield}}I - 2\varepsilon(\dimred{\vectorfield}) \right) \right) : D\dimred{\advelo} \dx{x} - \integral{\porous{\dimred{\Omega}}} \height\ \dimred{\pressure}\ \tr \left( D\dimred{\advelo} \left( \divergence{\dimred{\vectorfield}}I - D\dimred{\vectorfield} \right) \right) \dx{x} \\
		& -\integral{\porous{\dimred{\Omega}}} \height\ \dimred{\adpres}\ \tr \left( D\dimred{\velocity} \left( \divergence{\dimred{\vectorfield}}I - D\dimred{\vectorfield} \right) \right) \dx{x} + \integral{\porous{\dimred{\Omega}}_\subdarcy} \height\ \viscosity \permeability^{-1} \dimred{\velocity} \cdot \dimred{\advelo}\ \divergence{\dimred{\vectorfield}} \dx{x}  \\
		&+ \integral{\porous{\dimred{\Omega}}_\subfluid} \height \conductivity \left( \left( \divergence{\dimred{\vectorfield}}I - 2\varepsilon(\dimred{\vectorfield}) \right) \grad \dimred{\temperature} \right) \cdot \grad \dimred{\adtemp} \dx{x} + \integral{\porous{\dimred{\Omega}}_\subdarcy} \height \volumefrac \conductivity \left( \left( \divergence{\dimred{\vectorfield}}I - 2\varepsilon(\dimred{\vectorfield}) \right) \grad \dimred{\temperature} \right) \cdot \grad \dimred{\adtemp} \dx{x} \\
		&+ \integral{\porous{\dimred{\Omega}}} \height\ \density \hcapacity \dimred{\velocity} \cdot \left( \left( \divergence{\dimred{\vectorfield}}I - D\dimred{\vectorfield}^T \right) \grad \dimred{\temperature}\right) \dimred{\adtemp} \dx{x} + \integral{\porous{\dimred{\Gamma}}_\subwall} \height\ \htc \left( \dimred{\temperature} - \temperature_\subwall \right) \dimred{\adtemp}\ \tdiv{\dimred{\vectorfield}} \dx{s} \\
		&+ \integral{\porous{\dimred{\Omega}}_\subfluid} 2 \htc \left( \dimred{\temperature} - \temperature_\subwall \right) \dimred{\adtemp}\ \divergence{\dimred{\vectorfield}} \dx{x} + \integral{\porous{\dimred{\Omega}}_\subdarcy} \height\ \hfs \left( \dimred{\temperature} - \temperature_\subwall \right) \dimred{\adtemp}\ \divergence{\dimred{\vectorfield}} \dx{x},
	\end{aligned}
	\right.
\end{equation*}
where $(\dimred{\velocity}, \dimred{\pressure}, \dimred{\temperature}) = \dimred{\solution}(\porous{\dimred{\Omega}})$ is the solution of \eqref{eq:weak_darcy_dimension_state} and $(\dimred{\advelo}, \dimred{\adpres}, \dimred{\adtemp}) = \dimred{\adsolution}(\porous{\dimred{\Omega}})$ solves the adjoint system
\begin{equation*}
	\left\lbrace \quad
	\begin{aligned}
		&\text{Find } \dimred{\adsolution} = (\dimred{\advelo}, \dimred{\adpres}, \dimred{\adtemp}) \in \porous{\dimred{\adjointspace}} \text{ such that } \\
		&\qquad \integral{\porous{\dimred{\Omega}}_\subfluid} \height \conductivity \grad \dimred{\adtemp} \cdot \grad \dimred{\testadtemp} \dx{x} + \integral{\porous{\dimred{\Omega}}_\subdarcy} \height \volumefrac \conductivity \grad \dimred{\adtemp} \cdot \grad \dimred{\testadtemp} \dx{x} + \integral{\porous{\dimred{\Omega}}} \height\ \density \hcapacity\ \dimred{\velocity} \cdot \grad \dimred{\testadtemp} \dimred{\adtemp} \dx{x} + \integral{\porous{\dimred{\Omega}}_\subdarcy} \height\ \hfs \dimred{\adtemp} \dimred{\testadtemp} \dx{x}  \\
		&\qquad\ \ +\integral{\porous{\dimred{\Gamma}}_\subwall} \height\ \htc \dimred{\adtemp} \dimred{\testadtemp} \dx{s} + \integral{\porous{\dimred{\Omega}}_\subfluid} 2\htc \dimred{\adtemp} \dimred{\testadtemp} \dx{x} + \integral{\porous{\dimred{\Omega}}_\subfluid} \frac{6}{5} \height\ \viscosity\ D\dimred{\advelo} : D\dimred{\testadvelo} + \frac{12}{\height}\ \viscosity \dimred{\advelo} \cdot \dimred{\testadvelo} \dx{x} \\
		&\qquad\ \ +\integral{\porous{\dimred{\Omega}}_\subdarcy} \height\ \viscosity\ D\dimred{\advelo} : D\dimred{\testadvelo} + \height\ \viscosity \permeability^{-1} \dimred{\advelo} \cdot \dimred{\testadvelo} \dx{x} - \integral{\porous{\dimred{\Omega}}} \height\ \dimred{\adpres}\ \divergence{\dimred{\testadvelo}} + \height\ \dimred{\testadpres}\ \divergence{\dimred{\advelo}} - \height\ \density \hcapacity \dimred{\testadvelo} \cdot \grad \dimred{\temperature} \dimred{\adtemp} \dx{x} \\
		&\qquad = 2\weighttemp \Big( \porous{\dimred{\flux}}(\porous{\dimred{\Omega}}, \dimred{\solution}(\porous{\dimred{\Omega}})) - \fluxdes \Big) \left( \integral{\porous{\dimred{\Gamma}}_\subwall} \height\ \htc \dimred{\testadtemp} \dx{s} + \integral{\porous{\dimred{\Omega}}_\subfluid} 2\htc \dimred{\testadtemp} \dx{x} + \integral{\porous{\dimred{\Omega}}_\subdarcy} \height\ \hfs \dimred{\testadtemp} \dx{x} \right) \\
		&\qquad\ \ -2 \weightvelo \integral{\porous{\dimred{\Omega}}_\subdarcy} \height\ (\dimred{\velocity} - \porous{\velocity}_\subdes) \cdot \dimred{\testadvelo} \dx{x} \\
		&\text{for all } \dimred{\testadsolution} = (\dimred{\testadvelo}, \dimred{\testadpres}, \dimred{\testadtemp}) \in \porous{\dimred{\adjointspace}}.
	\end{aligned}
	\right.
\end{equation*}

\section{Numerical Comparison of the Reduced Models}
\label{sec:numerical_state}

After introducing all reduced models in Sections~\ref{sec:dimension_reduction} to \ref{sec:darcy_2D}, we now investigate how well they approximate the original model from Section~\ref{sec:problem_formulation}. To distinguish between the models we use the following naming: we call the model from Section~\ref{sec:problem_formulation}, i.e., the three-dimensional model without further approximations, the \qe{full 3D} model. In analogy, we call the model from Section~\ref{sec:dimension_reduction}, i.e., the one arising from applying the dimension reduction technique to the full 3D model, the \qe{full 2D} model. The models of Section~\ref{sec:porous_medium} and \ref{sec:darcy_2D} are then called \qe{Darcy 3D} and \qe{Darcy 2D} model, respectively. 

\subsection{Discretization and Numerical Solution of the PDEs}

The computational meshes are generated with the help of FreeCAD~0.16 \cite{freecad} and Gmsh~4.1.0 \cite{gmsh}. To get a comparable discretization in the $x$--$y$ plane, we generate the three-dimensional meshes by extruding the ones for the two-dimensional models. Further, for the Darcy models we use a similar discretization on the in- and outlet domains as for the full models, whereas the homogenized part of the domain $\porous{\Omega}_\subdarcy$ is discretized significantly coarser than $\Omega_\subchannels$ to reduce computational cost. We discretize all PDEs with the finite element method using FEniCS~2018.1 (cf. \cite{fenics, fenics_book}). For the fluid velocity we use quadratic Lagrange elements, and for both the pressure and temperature we use linear Lagrange elements, i.e, we use the LBB-stable Taylor-Hood elements for the Stokes system (see, e.g., \cite{john, elman}). For the temperature we additionally use a streamline-upwind Petrov-Galerkin (SUPG) method to stabilize the convection-dominated system (see, e.g., \cite{brooks_hughes, elman, john}).


\begin{table}[b]
	\centering
	{\footnotesize
	\rowcolors{2}{\tablegray}{white}
	\begin{tabular}{l r r r r r r r r}
		\toprule
		Model & \hspace{1em} & Vertices & \hspace{1em} & DoF's & \hspace{1em} & Time [\si{\second}] & \hspace{1em} & Memory [GB] \\
		\midrule
		Full 3D & & \num{5.71e+5} & & \num{1.33e+7} & & \num{1610.7} & & \num{50.91} \\
		Darcy 3D & & \num{2.50e+5} & & \num{5.85e+6} & & \num{467.5} & & \num{22.28} \\
		Full 2D & & \num{1.13e+5} & & \num{1.10e+6} & & \num{33.5} & & \num{2.44} \\
		Darcy 2D & & \num{5.00e+4} & & \num{4.96e+5} & & \num{18.6} & & \num{1.39} \\
		\bottomrule
	\end{tabular}
	\caption{Comparison of number of vertices and DoF's, as well as time and memory for solving the state systems.}
	\label{table:comparison_efficiency}
	}
\end{table}

The linear systems arising from the finite element method are solved with the library PETSc, version 3.10.5 (cf. \cite{petsc-user-ref}). For all two-dimensional problems we use the direct solver MUMPS. For the three-dimensional Stokes equations we use GMRES as it showed a better convergence than the MINRES method. As preconditioner, we use a FIELDSPLIT method based on the Schur complement, consisting of an ILU preconditioner for the velocity block and the algebraic multigrid preconditioner BOOMERAMG for the pressure block. For the three-dimensional convection-diffusion equation we again use GMRES with an ILU factorization as preconditioner. All linear systems are solved to a relative tolerance of \num{1e-10}.

In Table~\ref{table:comparison_efficiency} we compare the computational efficiency of the models w.r.t.\ the size of the mesh and linear systems, as well as the resources needed for the solution of the corresponding PDEs. We see that all reduced models lead to smaller meshes and, thus, also to smaller linear systems. Especially the dimension reduction yields a significant decrease in system size that comes from the fact that we use quadratic elements for the velocity and that one velocity component vanishes for the 2D models. Further, both the solution time and memory requirements also go down quite a bit from the full 3D model to the 3D Darcy one, and they decrease substantially when considering the two-dimensional models. 

\subsection{Comparison of the Models}

To compare the models we proceed as follows. First, we solve the PDEs on their respective domains. For the full 2D model we then only have to \qe{reverse} the dimension reduction using the profiles given in \eqref{eq:parabolic_profile_ansatz} and \eqref{eq:constant_profile_ansatz} to obtain three-dimensional quantities. For the 3D Darcy model we do not alter the solution in $\porous{\Omega}_\subfluid$, since we did not change the model there, and compute from the mean velocity in $\porous{\Omega}_\subdarcy$ the corresponding physical velocity as solution of Poiseuille flow in the channels. For the pressure and temperature we interpolate the values on $\porous{\Omega}_\subdarcy$ to $\Omega_\subchannels$ as they describe values extended to the larger averaged domain $\porous{\Omega}_\subdarcy$. Finally, we combine both approaches for the 2D Darcy model.

\begin{figure}[b]
	\centering
	\includegraphics[width=0.9\textwidth]{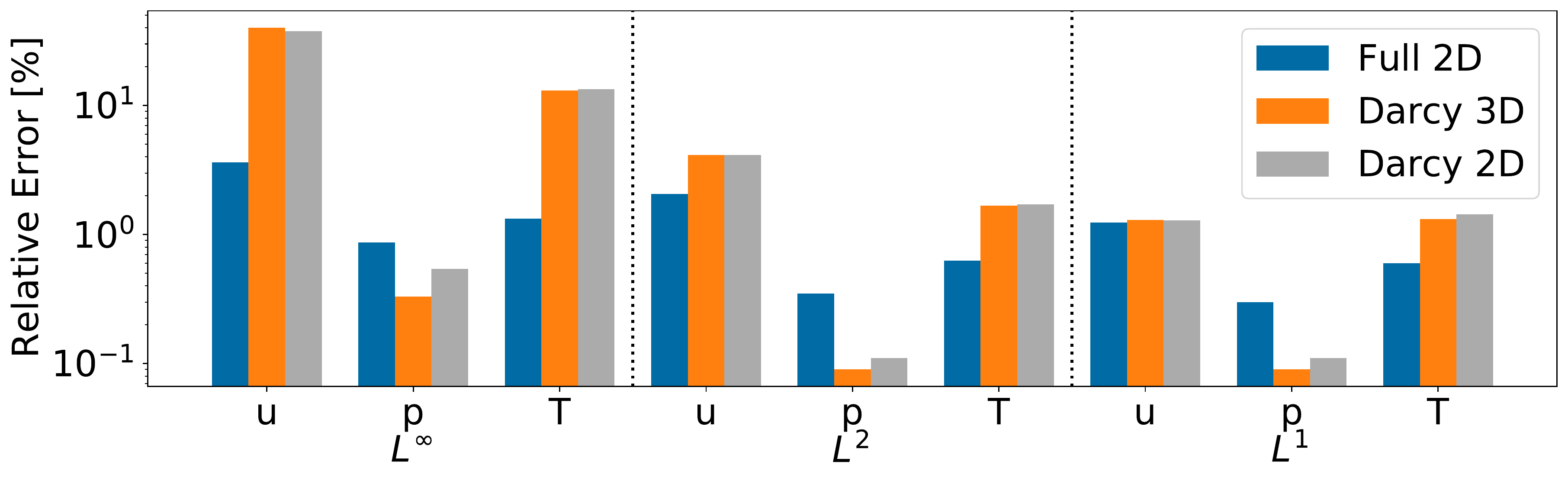}
	\caption{Comparison of the models' relative errors.}
	\label{figure:comparison_models}
\end{figure}

We computed both the absolute and relative errors of the reduced models to the full 3D model in three different norms: The $L^\infty(\Omega)$-norm (Table~\ref{table:Linf}), the $L^2(\Omega)$-norm (Table~\ref{table:L2}), and the $L^1(\Omega)$-norm (Table~\ref{table:L1}), and the relative errors are also shown in Figure~\ref{figure:comparison_models} in a logarithmic plot. Note, that the relative errors of the temperature are computed w.r.t.\ the reference temperature $\temperature_\subin$, i.e., we use
\begin{equation*}
	\frac{\norm{(\bar{\temperature} - \temperature_\subin) - (\temperature - \temperature_\subin)}{}}{\norm{\temperature - \temperature_\subin}{}} = \frac{\norm{\bar{\temperature} - \temperature}{}}{\norm{\temperature - \temperature_\subin}{}},
\end{equation*}
where $\temperature$ is the temperature computed by the full 3D model and $\bar{\temperature}$ is the temperature obtained from a reduced model.


\begin{table}[!t]
	\centering
	{\footnotesize
	\rowcolors{2}{\tablegray}{white}
	\begin{tabular}{l c r r c r r c r r}
		\toprule
		& \hspace{1.5em} & \multicolumn{2}{c}{velocity} & \hspace{1.5em} & \multicolumn{2}{c}{pressure} & \hspace{1.5em} & \multicolumn{2}{c}{temperature} \\
		\midrule
		Full 2D & & \num{2.40e-3} & (\num{3.63} \%) & & \num{1.36e+0} & (\num{0.87} \%) & & \num{9.75e-1} & (\num{1.33} \%) \\
		Darcy 3D & & \num{2.64e-2} & (\num{39.95} \%) & & \num{5.17e-1} & (\num{0.33} \%) & & \num{9.54e+0} & (\num{13.05} \%) \\
		Darcy 2D & & \num{2.50e-2} & (\num{37.81} \%) & & \num{8.46e-1} & (\num{0.54} \%) & & \num{9.76e+0} & (\num{13.35} \%) \\
		\bottomrule
	\end{tabular}
	\caption{Errors in the $L^\infty(\Omega)$ norm.}
	\label{table:Linf}
	}
\end{table}

\begin{table}[!t]
	\centering
	{\footnotesize
	\rowcolors{2}{\tablegray}{white}
	\begin{tabular}{l c r r c r r c r r}
		\toprule
		& \hspace{1.5em} & \multicolumn{2}{c}{velocity} & \hspace{1.5em} & \multicolumn{2}{c}{pressure} & \hspace{1.5em} & \multicolumn{2}{c}{temperature} \\
		\midrule
		Full 2D & & \num{3.48e-7} & (\num{2.07} \%) & & \num{3.21e-4} & (\num{0.35} \%) & & \num{2.45e-4} & (\num{0.63} \%) \\
		Darcy 3D & & \num{6.96e-7} & (\num{4.13} \%) & & \num{8.38e-5} & (\num{0.09} \%) & & \num{6.55e-4} & (\num{1.67} \%) \\
		Darcy 2D & & \num{6.95e-7} & (\num{4.12} \%) & & \num{1.05e-4} & (\num{0.11} \%) & & \num{6.69e-4} & (\num{1.71} \%) \\
		\bottomrule
	\end{tabular}
	\caption{Errors in the $L^2(\Omega)$ norm.}
	\label{table:L2}
	}
\end{table}

\begin{table}[!t]
	\centering
	{\footnotesize
	\rowcolors{2}{\tablegray}{white}
	\begin{tabular}{l c r r c r r c r r}
		\toprule
		& \hspace{1.5em} & \multicolumn{2}{c}{velocity} & \hspace{1.5em} & \multicolumn{2}{c}{pressure} & \hspace{1.5em} & \multicolumn{2}{c}{temperature} \\
		\midrule
		Full 2D & & \num{2.00e-10} & (\num{1.24} \%) & & \num{2.59e-7} & (\num{0.30} \%) & & \num{2.17e-7} & (\num{0.60} \%) \\
		Darcy 3D & & \num{2.10e-10} & (\num{1.30} \%) & & \num{7.34e-8} & (\num{0.09} \%) & & \num{4.79e-7} & (\num{1.32} \%) \\
		Darcy 2D & & \num{2.09e-10} & (\num{1.29} \%) & & \num{9.32e-8} & (\num{0.11} \%) & & \num{5.20e-7} & (\num{1.43} \%) \\
		\bottomrule
	\end{tabular}
	\caption{Errors in the $L^1(\Omega)$ norm.}
	\label{table:L1}
	}
\end{table}

The largest difference between all models can be seen for the $L^\infty(\Omega)$-norm. There, the errors in velocity and temperature for the full 2D model are ten times as small as the ones for both Darcy models. This comes from the fact that the full 2D model is posed on the domain $\dimred{\Omega}$ that still includes the microchannels, whereas the Darcy models utilize an averaged domain. In particular, the flow of the coolant into and out of the channels is only resolved by the full 2D model, which explains the large $L^\infty$ error for the other ones. However, this is remedied when we consider the $L^2$ and $L^1$ norms. There, we observe that the Darcy models are still worse compared to the full 2D model, but the difference of the errors is now considerably smaller. 
Further, we observe that the pressure is approximated very well by all models since the corresponding relative errors are below \num{1}~\% in all considered norms. Altogether, except for the $L^\infty$ norm, all relative errors are below \num{5}~\% for all quantities, which indicates that our reduced models work rather well.


In addition to this, we see that both Darcy models have very similar errors in all considered norms. This suggests that most of the error of the 2D Darcy model comes from the modeling as a porous medium, and not from the subsequent dimension reduction. From this, we conclude that we can use the 2D Darcy model instead of the 3D one as their errors are nearly identical, but the 2D model is significantly more efficient (cf. Table~\ref{table:comparison_efficiency}). Finally, the full 2D model shows the best overall performance of the reduced models due to its combination of efficiency and accuracy.


\section{Numerical Results for the Optimization Problem}
\label{sec:numerical_optimization}

Now, we describe the numerical solution of shape optimization problems using \eqref{eq:def_opt_problem} as an example and then explain the details and modifications we use for the reduced models. Afterwards, we discuss the results obtained from the different models and compare them.

\subsection{Numerical Solution of the Shape Optimization Problem}

In the previous sections we presented the shape derivatives for our optimization problems. However, for numerical purposes this is not yet sufficient as the shape derivative only gives us the sensitivity of the cost functional for a given deformation generated by the flow of a vector field $\vectorfield$. In particular, the shape derivative $d\reducedcostfunction(\Omega)[\vectorfield]$ is a linear functional acting on $\vectorfield$. To solve the shape optimization problem numerically, we compute a descent direction $\searchdirection$, i.e., a vector field whose associated flow gives a descent in the cost functional, from the shape derivative. To do so, we choose a symmetric, continuous and coercive bilinear form $a\colon H(\Omega) \times H(\Omega) \to \R$, where $H(\Omega)$ denotes some suitable Hilbert space on the current geometry $\Omega$ which is specified later on. Then, the shape gradient $\shapegradient$ w.r.t. $a$ is defined as the solution of the variational equation
\begin{equation}
	\label{eq:variational_shape_gradient}
	\text{Find } \shapegradient \in H(\Omega) \text{ such that } \quad a(\shapegradient, \vectorfield) = d\reducedcostfunction(\Omega)[\vectorfield] \quad \text{ for all } \vectorfield\in H(\Omega).
\end{equation}
Thanks to the Lax-Milgram Lemma this has a unique solution $\shapegradient$. 

As in the case of nonlinear optimization, the negative shape gradient is a descent direction as
\begin{equation*}
	d\reducedcostfunction(\Omega)[-\shapegradient] = -a(\shapegradient, \shapegradient) \leq 0,
\end{equation*}
due to the coercivity of $a$. Therefore, we use the negative shape gradient in a line search method to solve problem \eqref{eq:def_opt_problem}. After computing the search direction $\searchdirection = - \shapegradient$, we deform the domain numerically with the so-called perturbation of identity given by
\begin{equation}
	\label{eq:perturbation_of_identity}
	\Omega_t = (I + t\searchdirection)\Omega := \Set{x + t\searchdirection(x) | x\in \Omega}.
\end{equation}
Starting from an initial guess, the step size $t$ is accepted if it satisfies the Armijo condition
\begin{equation}
	\label{eq:armijo}
	\reducedcostfunction(\Omega_t) \leq \reducedcostfunction(\Omega) + \sigma t\ d\reducedcostfunction(\Omega)[\searchdirection] = \reducedcostfunction(\Omega) + \sigma t\ a(\shapegradient, \searchdirection),
\end{equation}
where the last equality holds due to \eqref{eq:variational_shape_gradient}, and the parameter $\sigma \in (0,1)$ is chosen to be $\sigma = \num{1e-4}$ (cf. \cite{kelley, nocedal_wright}). If this is not satisfied, the step size is halved and the procedure is repeated. Additionally, we do a mesh quality control based on conditions given in \cite{herzog} to avoid step sizes that lead to excessively large deformations. We only accept step sizes that satisfy
\begin{equation}
	\label{eq:quality}
	\frac{1}{2} \leq \det\left( I + t\ D\searchdirection \right) \leq 2 \quad \text{ as well as } \quad t \norm{D\searchdirection}{F} \leq 0.3,
\end{equation}
for each element of the mesh, where $\norm{\cdot}{F}$ denotes the Frobenius norm of a matrix. In case the step size is accepted, we double it after the deformation of the domain to get a good initial guess for the next iteration.

For the stopping criterion we make the following considerations. As $a$ is symmetric, continuous, and coercive, we see that it defines an inner product on $H(\Omega)$. Hence, for the norm of the shape gradient we use the one induced by this scalar product, i.e., $\norm{\shapegradient}{H(\Omega)} := \sqrt{a(\shapegradient, \shapegradient)}$ (cf. \cite{welker, herzog}). We denote the shape gradient on the initial geometry $\Omega_0$ by $\shapegradient_0$ and terminate the algorithm if the relative stopping criterion
\begin{equation}
	\label{eq:stopping_criterion}
	\frac{\norm{\shapegradient}{H(\Omega)}}{\norm{\shapegradient_0}{H(\Omega_0)}} \leq \varepsilon,
\end{equation}
is satisfied. For all of our numerical experiments we choose the relative tolerance as $\varepsilon = \num{1e-3}$. Additionally, we also stop the optimization after a fixed amount of iterations or in case the conditions \eqref{eq:armijo} and \eqref{eq:quality} cannot be satisfied.

The general numerical method is summarized in Algorithm~\ref{algo:descent}, where problem \eqref{eq:def_opt_problem} is used as an example. The only thing left to do before we can apply it is the specification of the bilinear form $a$ for all models, which is done in the subsequent section.

\begin{algorithm2e}
	\KwIn{Initial domain $\Omega_0$, initial step size $t$, tolerance $\varepsilon$, maximum number of iterations $m$}
	\For{k=0,1,2,\dots, $m$}{
		Compute the solution of the state system, e.g., \eqref{eq:weak_state} \\
		Compute the solution of the adjoint system, e.g., \eqref{eq:weak_adjoint} \\
		Compute the shape gradient $\shapegradient_k$ via \eqref{eq:variational_shape_gradient}, e.g., with $d\reducedcostfunction(\Omega_k)[\vectorfield]$ as in \eqref{eq:shape_derivative}\\
		\If{Condition \eqref{eq:stopping_criterion} is satisfied for $\shapegradient = \shapegradient_k$}{
			Stop with approximate solution $\Omega_k$
		}
		Define the search direction as $\searchdirection_k = -\shapegradient_k$ \\
		\While{Either \eqref{eq:armijo} or \eqref{eq:quality} with $\Omega = \Omega_k$ and $\searchdirection = \searchdirection_k$ are not satisfied}{
			Decrease the step size: $t = \nicefrac{t}{2} $
		}
		Update the geometry via \eqref{eq:perturbation_of_identity}: $\Omega_{k+1} = (I + t \searchdirection_k)\Omega_k$ \\
		Increase the step size for the next iteration: $t = 2 t$
	}
	\caption{Numerical solution of shape optimization problems.}
	\label{algo:descent}
\end{algorithm2e}

\subsection{Choice of the Bilinear Form for the Shape Gradient}

The bilinear forms used for computing the shape gradient in \eqref{eq:variational_shape_gradient} are based on the equations of linear elasticity which are used, e.g., in \cite{welker, hohmann, berggren}. We modify this approach and use equations of anisotropic, inhomogeneous, linear elasticity for our problems. 

\subsubsection*{Full 3D Model}
Let us start with the model from Section~\ref{sec:problem_formulation}. We define the Hilbert space as
\begin{equation*} 
H(\Omega) = \Set{\vectorfield \in H^1(\Omega)^3 | \vectorfield = 0 \text{ on } \Gamma_\subin \cup \Gamma_\subout,\ \vectorfield \cdot \normal = 0 \text{ on } \Set{z=0} \cup \Set{z=\height} \cup \Gamma_\subchannels}.
\end{equation*}
Note, that this choice is consistent for our optimization problem, i.e., it respects the geometrical constraints stated in Section~\ref{sec:optimization_problem}. For $\vectorfield \in H(\Omega)$ we get by the perturbation of identity \eqref{eq:perturbation_of_identity} for $\Omega_t = (I + t\vectorfield) \Omega$ that the boundaries $\Gamma_\subin$ and $\Gamma_\subout$ are fixed for all $t\in [0,\tau]$. Furthermore, thanks to the slip condition $\vectorfield \cdot \normal = 0$ on $\Set{z=0} \cup \Set{z=\height} \cup \Gamma_\subchannels$, the height of the geometry remains fixed and the microchannels can only change their length. For the computation of the shape gradient, we use the bilinear form $a\colon H(\Omega)\times H(\Omega) \to \R$ given by
\begin{equation}
	\label{eq:anisotropic_full}
	\left\lbrace 
	\begin{aligned}
		a(U, V) =& \integral{\Omega} \nu(x)\ \Big( \sigma(U) : \varepsilon(V) + \delta\ U \cdot V \Big) \dx{x}, \qquad \sigma(U) = \lambda_\text{elas}\ \tr(\varepsilon(U))I + 2\mu_\text{elas}\ E(U),\\
		E(U) =& \left[ \begin{array}{r r r}
			\varepsilon(U)_{1,1} & \varepsilon(U)_{1,2} & C\ \varepsilon(U)_{1,3} \\
			\varepsilon(U)_{2,1} & \varepsilon(U)_{2,2} & C\ \varepsilon(U)_{2,3} \\
			C\ \varepsilon(U)_{3,1} & C\ \varepsilon(U)_{3,2} & C\ \varepsilon(U)_{3,3}
		\end{array} \right].
	\end{aligned}
	\right.
\end{equation} 
Here, $\lambda_\text{elas}$ and $\mu_\text{elas}$ are the so-called Lam\'e parameters and $\delta \geq 0$ is a damping parameter. Further, $C \gg 1$ is a numerical constant that leads to an anisotropic strain tensor $E(U)$, which we choose as $C := \num{1e5}$. Hence, the $z$-component of the shape gradient is small in comparison to its other components. This is done to avoid unnecessarily small step sizes due to a deformation in $z$-direction that does not have an actual effect on the geometry as its height is already fixed. Finally, the term $\nu(x)$ models an inhomogeneous stiffness of the geometry which can only be defined after discretization of the geometry with a triangulation $T_h$. Then, it is given by
\begin{equation*}
	\nu(x) = \frac{\max_{T_k \in T_h} \abs{T_k}}{\abs{T_n}} \quad \text{ for } x \in T_n,
\end{equation*}
i.e., it is one over the relative ($d$-dimensional) volume of the considered element. This idea from \cite{berggren} ensures that large elements have a lower stiffness than small ones since they absorb large deformations better. As in \cite{berggren} we compute $\nu$ once on the initial mesh and do not update it during the optimization process so that elements cannot become arbitrarily stiff. We found that using this approach significantly increased the mesh quality during the optimization and that we did not have to remesh at all.

\subsubsection*{Full 2D Model}
For the dimension reduction model we proceed similarly to Section~\ref{sec:dimension_reduction}. We already assumed that the $z$-component of the deformation vanishes to derive the two-dimensional shape derivative (cf. \eqref{eq:ansatz_deformation}). This is now also used to derive a dimension-reduced formulation of the bilinear form $a$. To this end, we introduce the Hilbert space
\begin{equation*}
	\dimred{H}(\dimred{\Omega}) = \Set{\dimred{\vectorfield} \in H^1(\dimred{\Omega})^2 | \dimred{\vectorfield} = 0 \text{ on } \dimred{\Gamma}_\subin \cup \dimred{\Gamma}_\subout,\ \dimred{\vectorfield} \cdot \normal = 0 \text{ on } \dimred{\Gamma}_\subchannels}.
\end{equation*}
Using condition \eqref{eq:ansatz_deformation} in \eqref{eq:anisotropic_full} then yields the following bilinear form on $\dimred{H}(\dimred{\Omega})$
\begin{equation*}
	\dimred{a}(\dimred{U}, \dimred{V}) = \integral{\dimred{\Omega}} \height\ \nu(x)\  \Big( \dimred{\sigma}(\dimred{U}) : \varepsilon(\dimred{V}) + \delta \dimred{U} \cdot \dimred{V} \Big) \dx{x}, \quad \text{ where } \quad \dimred{\sigma}(\dimred{U}) = \lambda_\text{elas}\ \tr(\varepsilon(\dimred{U}))I + 2\mu_\text{elas}\ \varepsilon(\dimred{U}).
\end{equation*}

\subsubsection*{3D Darcy Model}

For the Darcy model we use the Hilbert space
\begin{equation*}
	\porous{H}(\porous{\Omega}) = \Set{\vectorfield\in H^1(\porous{\Omega})^3 | \vectorfield = 0 \text{ on } \porous{\Gamma}_\subin \cup \porous{\Gamma}_\subout,\ \vectorfield \cdot \normal = 0 \text{ on} \Set{z=0} \cup \Set{z=\height} \cup \porous{\Gamma}_\subdarcy},
\end{equation*}
which acts as an equivalent of $H(\Omega)$ for the domain $\porous{\Omega}$. As before, the perturbation of identity with a vector field from $\porous{H}$ leaves both the height $\height$ of $\porous{\Omega}$ and the boundaries $\porous{\Gamma}_\subin$ and $\porous{\Gamma}_\subout$ fixed. We choose the bilinear form
\begin{align*}
	\porous{a}(U, V) &= \integral{\porous{\Omega}_\subfluid} \nu(x)\ \Big( \sigma(U) : \varepsilon(V) + \delta\ U\cdot V \Big) \dx{x} + \integral{\porous{\Omega}_\subdarcy} \nu(x)\ \Big( \sigma_\subdarcy(U) : \varepsilon(V) + \delta\ U\cdot V \Big) \dx{x}, \\
	\sigma_\subdarcy(U) &= \lambda_\text{elas}\ \tr(\varepsilon(U))I + 2\mu_\text{elas}\ E_\subdarcy(U), \quad E_\subdarcy(U) = \left[ \begin{array}{r r r}
		C\ \varepsilon(U)_{1,1} & \varepsilon(U)_{1,2} & C\ \varepsilon(U)_{1,3} \\
		\varepsilon(U)_{2,1} & \varepsilon(U)_{2,2} & C\ \varepsilon(U)_{2,3} \\
		C\ \varepsilon(U)_{3,1} & C\ \varepsilon(U)_{3,2} & C\ \varepsilon(U)_{3,3}
	\end{array} \right],
\end{align*}
and $\sigma(U)$ is given in \eqref{eq:anisotropic_full}. As before, the Darcy model coincides with the full model in $\porous{\Omega}_\subfluid$, and in $\porous{\Omega}_\subdarcy$ we have the following modification. 
The strain tensor $E_\subdarcy$ in $\porous{\Omega}_\subdarcy$ now also exhibits a large stiffness in the $x$-direction. This models the influence of the microchannels on the shape gradient, due to the following reason. As we have lots of narrow channels in $\Omega_\subchannels$ and use the slip condition $\vectorfield \cdot \normal = 0$ on $\Gamma_\subchannels$, the geometry of the channels is only allowed to stretch or compress along the $y$-axis. Therefore, we can neglect the $x$-component of the shape gradient in $\Omega_\subchannels$. The anisotropic strain tensor $E_\subdarcy$ achieves a similar effect by increasing the stiffness in $x$-direction in $\porous{\Omega}_\subdarcy$.

\subsubsection*{2D Darcy Model}
The 2D Darcy model again combines the dimension reduction technique with the porous medium model. As before, we assume that the $z$-component of the deformation vanishes, and now use the Hilbert space
\begin{equation*}
	\porous{\dimred{H}}(\porous{\dimred{\Omega}}) = \Set{\dimred{\vectorfield} \in H^1(\porous{\dimred{\Omega}})^2 | \dimred{\vectorfield} = 0 \text{ on } \porous{\dimred{\Gamma}}_\subin \cup \porous{\dimred{\Gamma}}_\subout,\ \dimred{\vectorfield}\cdot \normal = 0 \text{ on } \porous{\dimred{\Gamma}}_\subdarcy}.
\end{equation*}
The bilinear form for the 2D Darcy model then reads
\begin{align*}
	\porous{\dimred{a}}(\dimred{U}, \dimred{V}) &= \integral{\porous{\dimred{\Omega}}_\subfluid} \height\ \nu(x)\ \Big( \dimred{\sigma}(\dimred{U}) : \varepsilon(\dimred{V}) + \delta \dimred{U} \cdot \dimred{V} \Big) \dx{x} + \integral{\porous{\dimred{\Omega}}_\subdarcy} \height\ \nu(x)\ \Big( \dimred{\sigma}_\subdarcy(\dimred{U}) : \varepsilon(\dimred{V}) + \delta \dimred{U} \cdot \dimred{V} \Big) \dx{x}, \\
	\dimred{\sigma}_\subdarcy(\dimred{U}) &= \lambda_\text{elas}\ \tr(\varepsilon(\dimred{U}))I + 2\mu_\text{elas}\ \dimred{E}_\subdarcy(\dimred{U}), \quad \dimred{E}_\subdarcy(\dimred{U}) = \left[ \begin{array}{r r}
	C\ \varepsilon(\dimred{U})_{1,1} & \varepsilon(\dimred{U})_{1,2} \\
	\varepsilon(\dimred{U})_{2,1} & \varepsilon(\dimred{U})_{2,2}
	\end{array} \right].
\end{align*}

For all models we choose $\mu_\text{elas} = \num{1}$ as well as $\lambda_\text{elas} = \num{1e-1}$ for the Lam\'e parameters. The damping parameter was chosen to be $\delta = \num{1e-1}$. For the numerical solution of \eqref{eq:variational_shape_gradient} we again use FEniCS. We discretize the equations using linear Lagrange elements and solve the resulting linear systems using MUMPS for the two-dimensional problems, and a conjugate gradient method preconditioned with BOOMERAMG for the three-dimensional ones.

\subsection{Numerical Results of the Optimal Shape Design Problem}

\begin{figure}[t]
	\centering
	\begin{subfigure}{0.49\textwidth}
		\centering
		\includegraphics[width=0.8\textwidth]{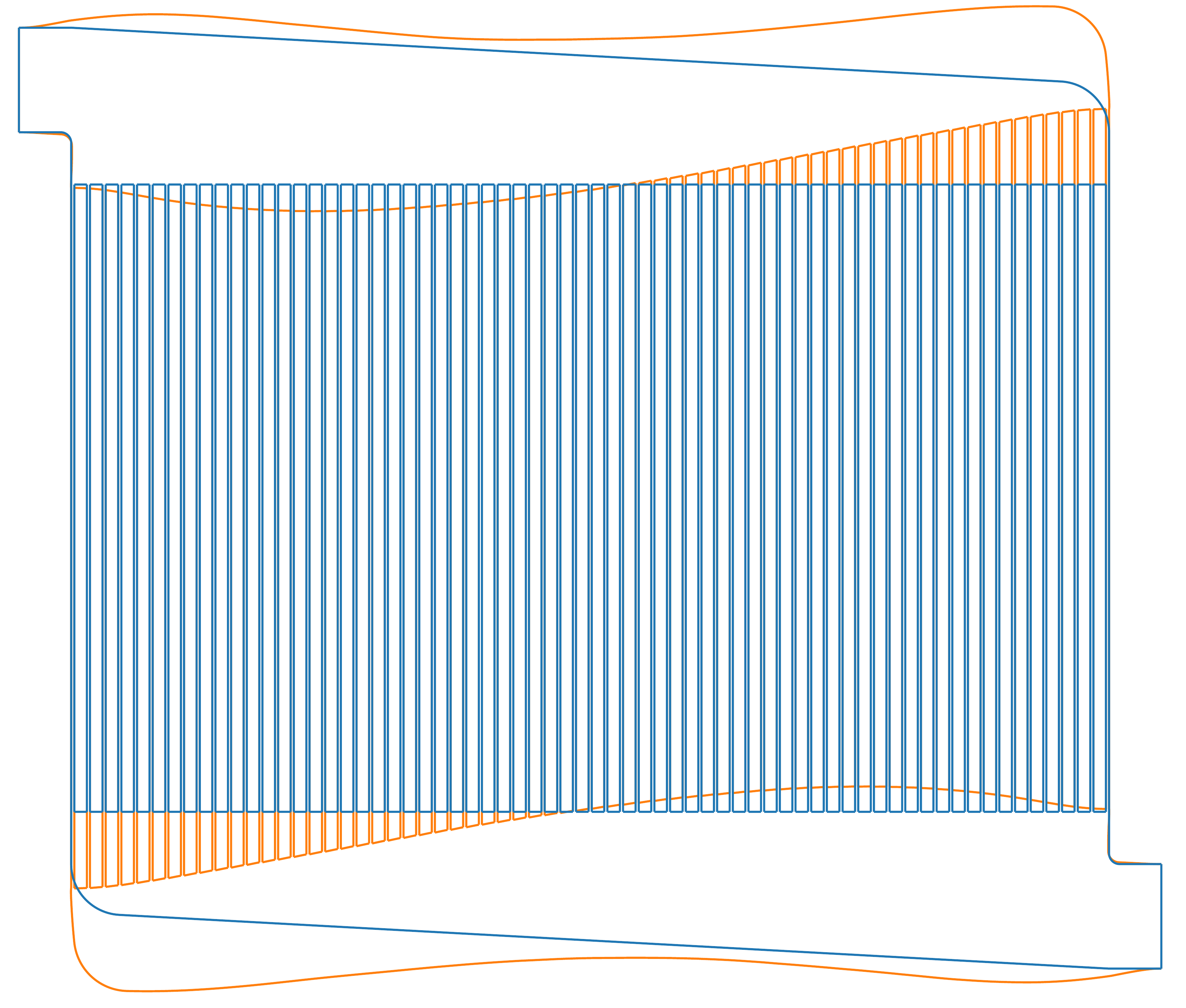}
		\caption{Comparison of initial shape (blue) and optimized one (orange).}
		\label{figure:difference}
	\end{subfigure}
	\hfill
	\begin{subfigure}{0.49\textwidth}
		\includegraphics[width=\textwidth]{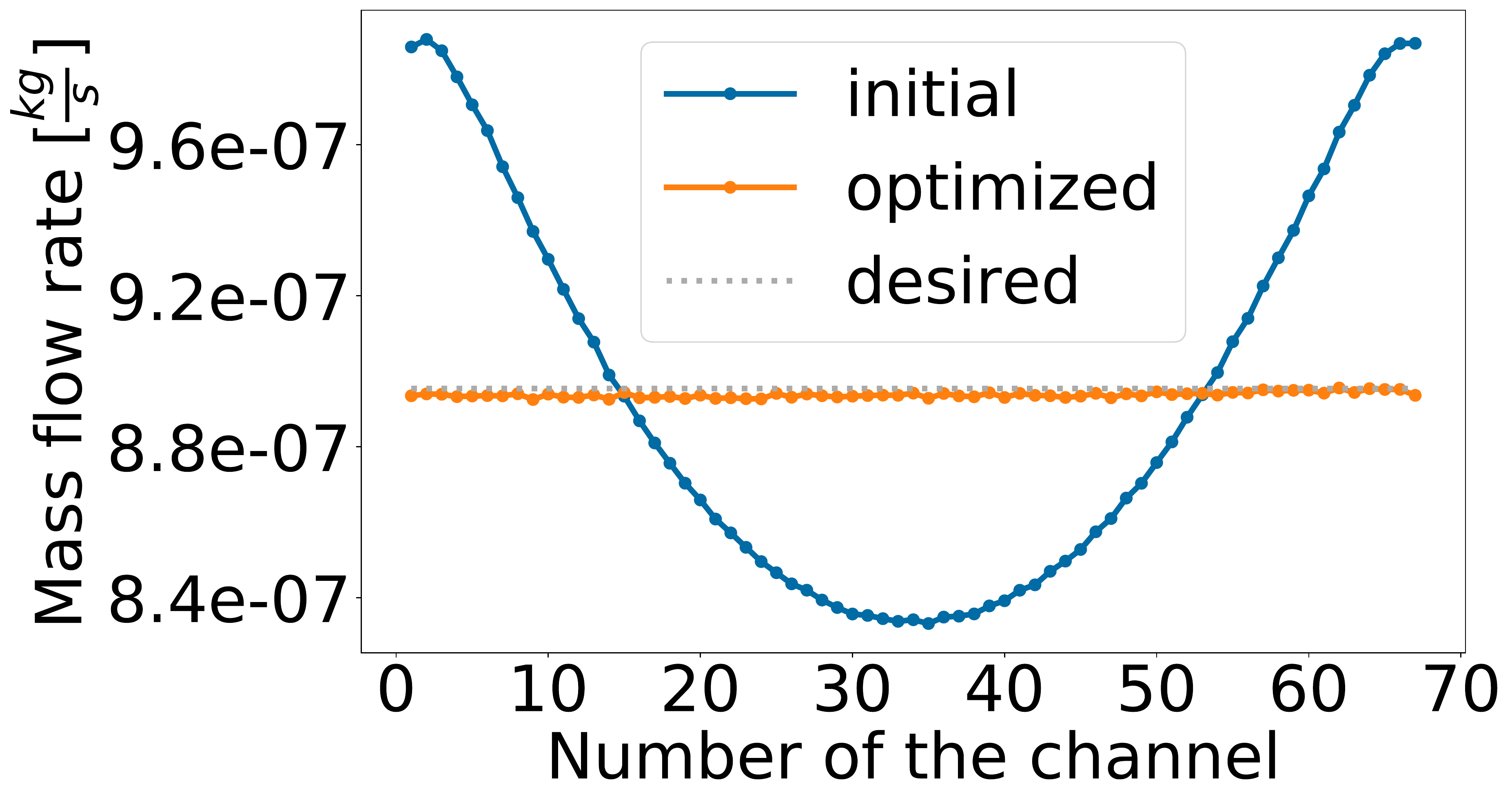}
		\caption{Distribution of the flow into the channels.}
		\label{figure:flow}
	\end{subfigure}
	\caption{Results of the optimization for the full 3D model.}
	\label{figure:full_3D}
\end{figure}

After giving the details of the numerical solution of the shape optimization problems as well as the choice of the bilinear forms for computing the shape gradient for all models, let us now investigate the results we obtained. We choose the weights for the cost functionals as
\begin{equation*}
	\weighttemp = \frac{1}{\costfunction_1(\Omega_0, \solution(\Omega_0))}, \qquad \weightvelo = \frac{\num{1}}{\costfunction_2(\Omega_0, \solution(\Omega_0))}, \qquad \text{ and } \qquad \weightreg = \frac{\num{1e-2}}{\costfunction_3(\Omega_0, \solution(\Omega_0))},
\end{equation*}
where $\Omega_0$ denotes the initial geometry. The scaling of the cost functional is chosen so that we weight the functions $\costfunction_1$ and $\costfunction_2$ equally, while having only a slight regularization from $\costfunction_3$. For the reduced models we choose an analogous scaling. Finally, we discuss the choice of $\velocity_\subdes$ and $\fluxdes$. For the former, we choose the velocity corresponding to the case of uniformly distributed flow among the microchannels, which we compute numerically as Poiseuille flow. For the latter, we have the following considerations. The cooling system absorbs \num{6.74}~\si{\watt} of heat with its initial shape. For the optimized cooler we want to increase this and, hence, choose $\fluxdes = \num{7.1}~\si{\watt}$, which corresponds to an increment of about \num{5}~\%.


\begin{figure}[b]
	\centering
	\begin{subfigure}[b]{0.49\textwidth}
		\centering
		\includegraphics[width=\textwidth]{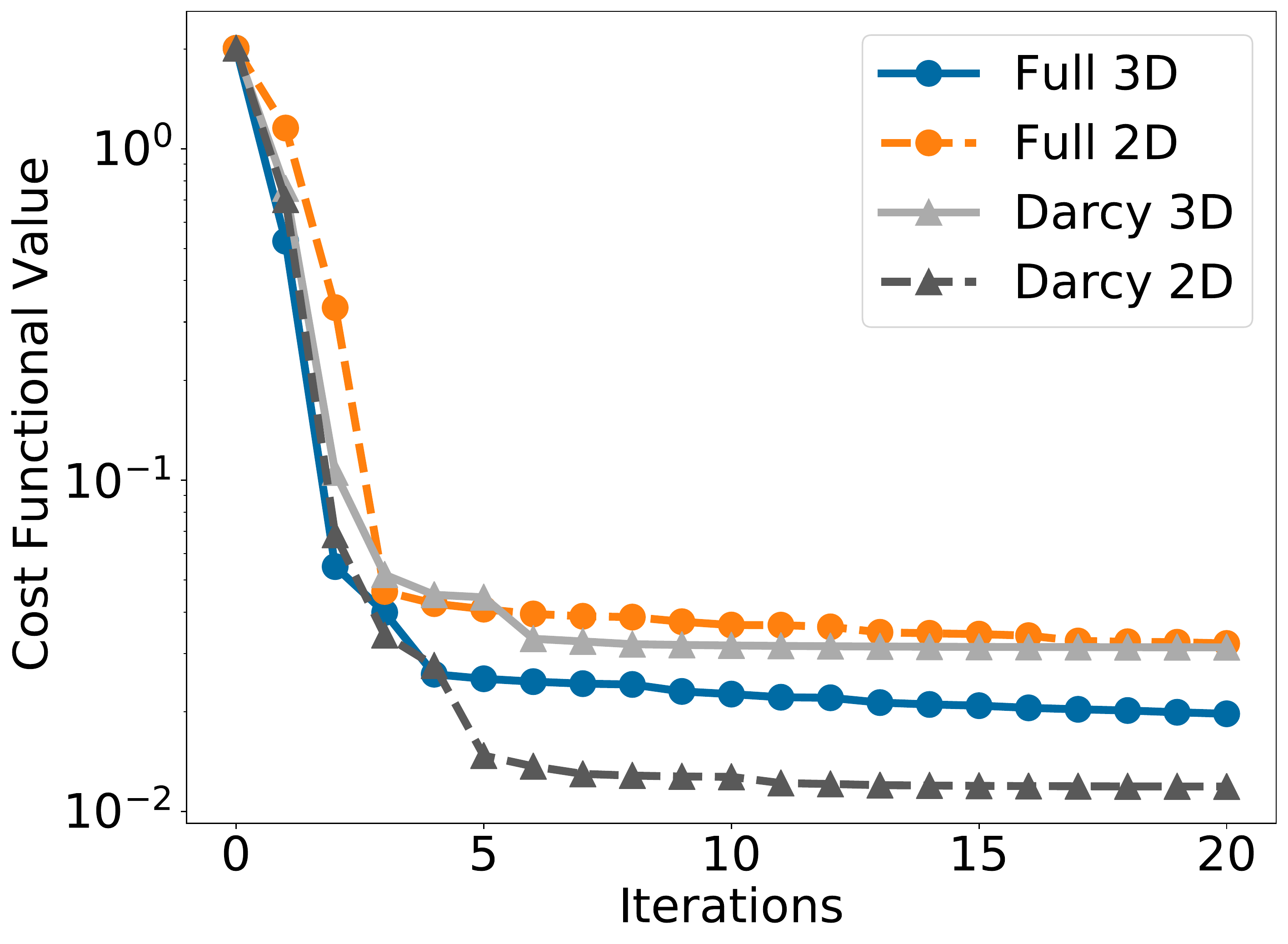}
		\caption{Value of the cost functional.}
	\end{subfigure}
	\hfil
	\begin{subfigure}[b]{0.49\textwidth}
		\centering
		\includegraphics[width=\textwidth]{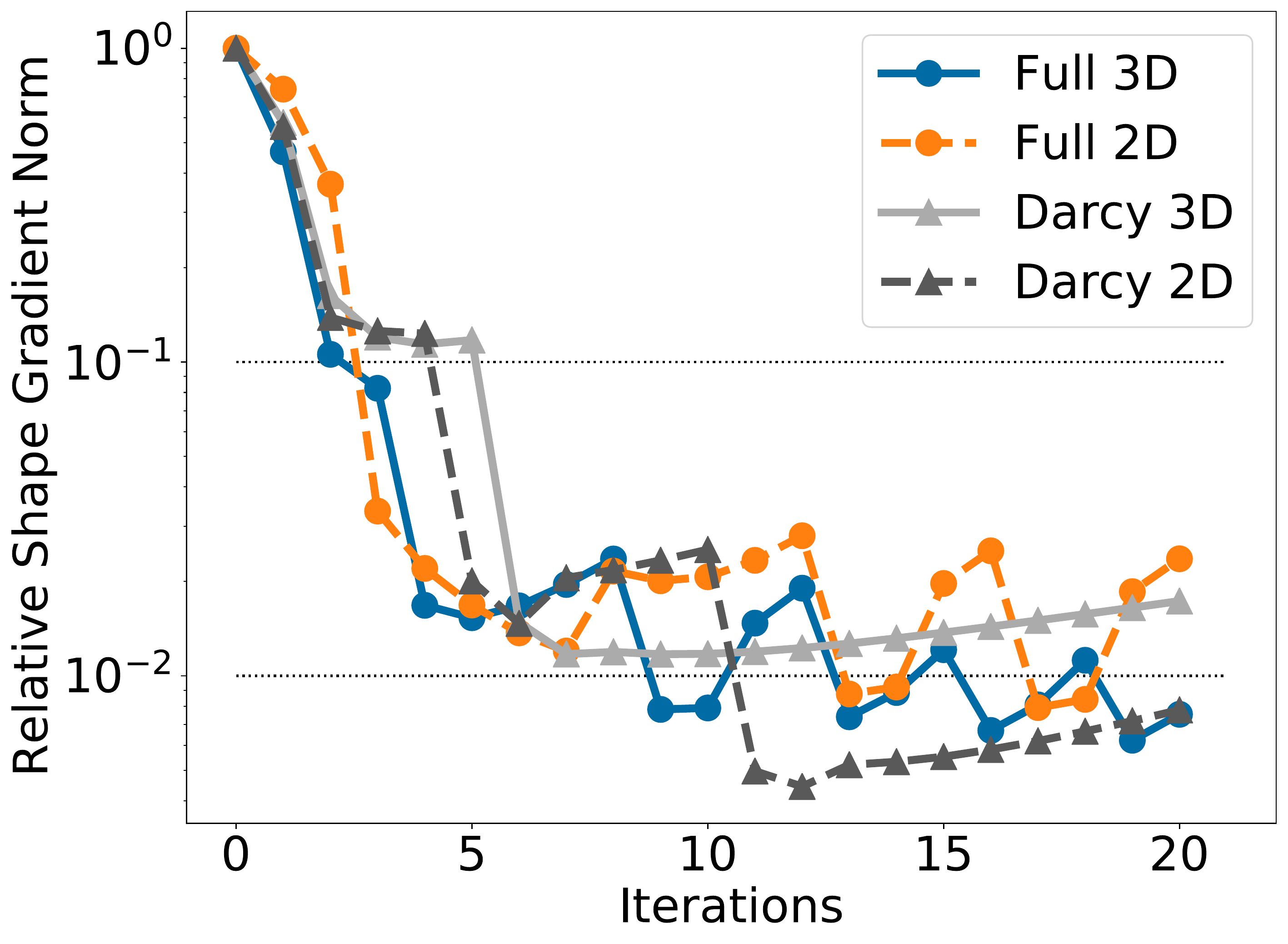}
		\caption{Relative norm of the gradient.}
	\end{subfigure}
	\caption{History of the optimization process.}
	\label{figure:history}
\end{figure}

The results of the shape optimization for the full 3D model can be seen in Figure~\ref{figure:full_3D}. In Figure~\ref{figure:difference} the initial and optimized geometries are shown. First, we observe that both the initial and optimized geometry are (nearly) point-symmetric to the center of the geometry. For the optimized shape, we see that the in- and outlet domains are pushed to the outside. Additionally, they are dented in the middle, creating a kind of U-shape (on the top). The length of the channels changed accordingly, to balance the pressure differences generated in the in- and outlet domains. The results of this can be seen in Figure~\ref{figure:flow}, where the mass flow rate of the coolant in the channels is depicted. We see that for the initial geometry the flow through the channels resembles a U-shape. In particular, the outer channels get the most amount of fluid and the middle ones get the least amount. This discrepancy is removed on the optimized domain, where we observe a nearly uniform flow distribution among all channels, achieving one goal of the optimization. Additionally, the heat absorbed by the cooler increases from \num{6.74}~\si{\watt} on the initial geometry to \num{7.094}~\si{\watt} on the optimized one, nearly reaching the desired amount of \num{7.1}~\si{\watt}. Moreover, the size of the geometry did not increase substantially. Hence, the optimized cooler could be used in place of the initial configuration. As the reduced models perform nearly indistinguishably with regards to these two objectives, we do not show the analogous results for them for the sake of brevity.


The history of the optimization process is shown in Figure~\ref{figure:history} for all four models, where both the value of the cost functional and the relative norm of the gradient are shown. We observe that the function value does not decrease much further after five iterations for all models which indicates that all geometries converge very quickly to the optimal one. Note, that the cost functional value is comparatively similar for all four models, again suggesting that they perform similarly. We terminated the optimization algorithm after 20 iterations since the norm of the gradient nearly decreased by two orders of magnitude for all models, which is sufficient for industrial applications. These results indicate that Algorithm~\ref{algo:descent} converged and found a (local) minimizer. 


\begin{figure}[t]
	\centering
	\begin{subfigure}{0.49\textwidth}
		\centering
		\includegraphics[width=0.9\textwidth]{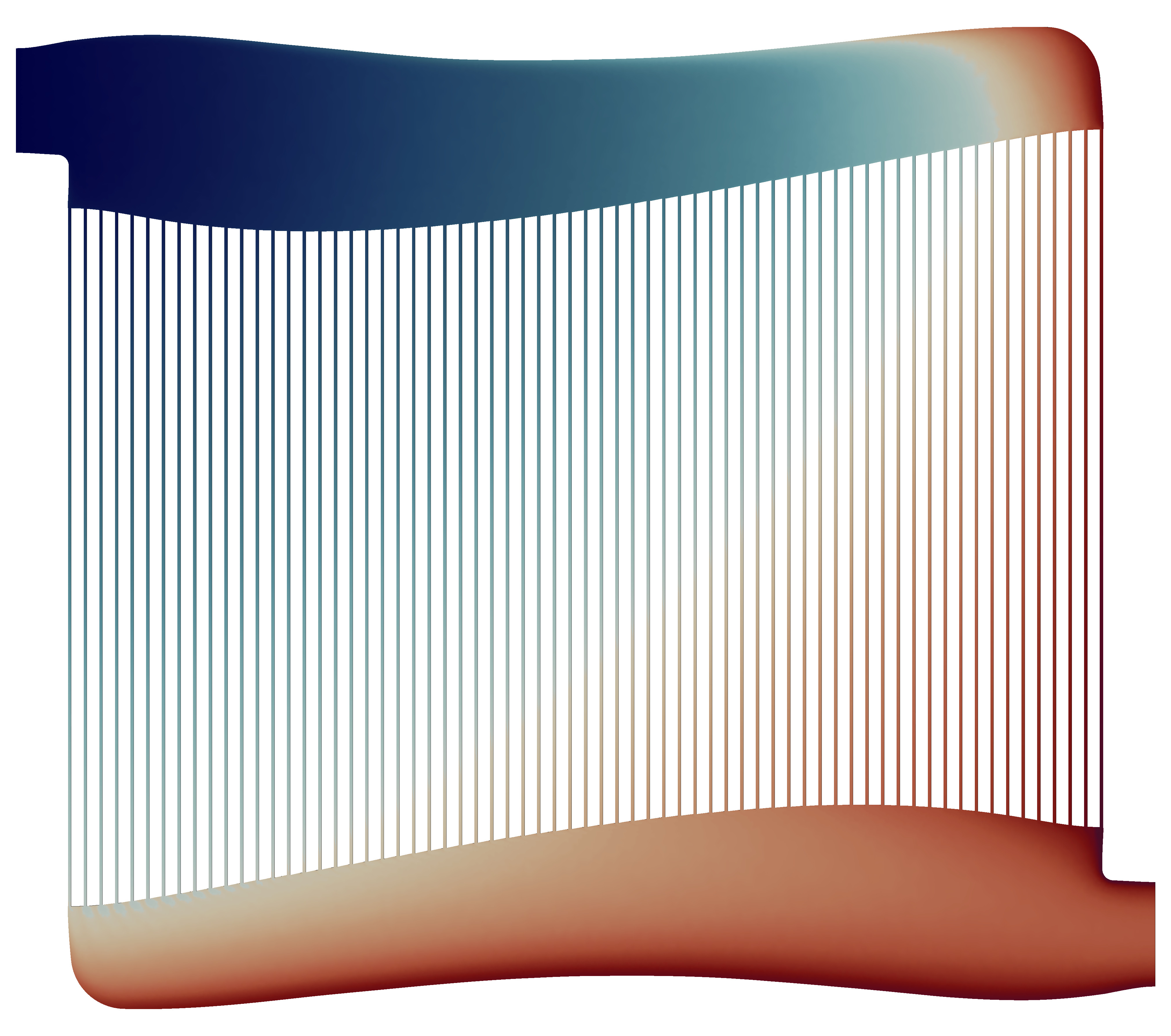}
		\caption{Full 3D Model.}
	\end{subfigure}
	\hfill
	\begin{subfigure}{0.49\textwidth}
		\centering
		\includegraphics[width=0.9\textwidth]{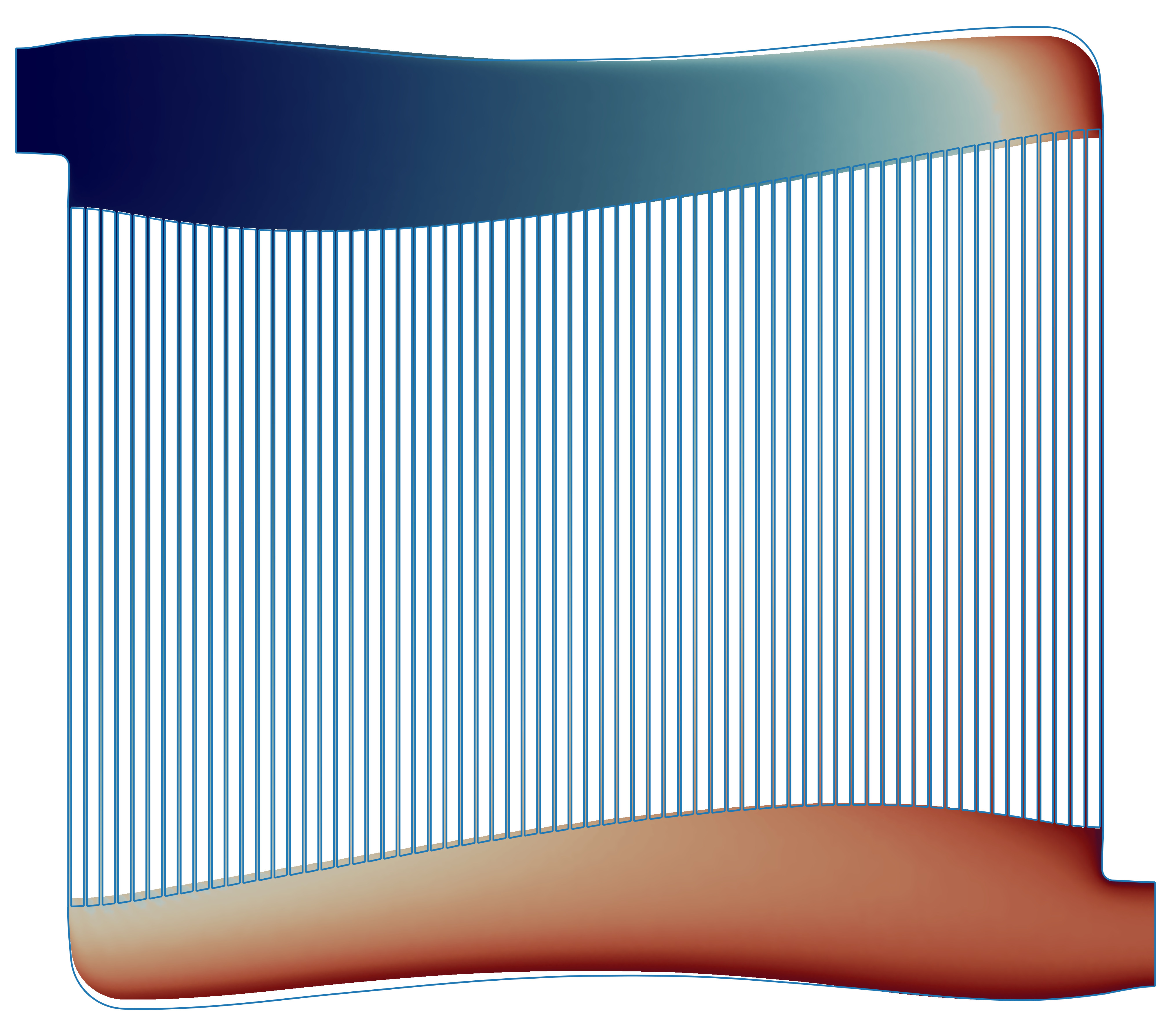}
		\caption{Full 2D Model.}
	\end{subfigure}
	\\
	\begin{subfigure}{0.49\textwidth}
		\centering
		\includegraphics[width=0.9\textwidth]{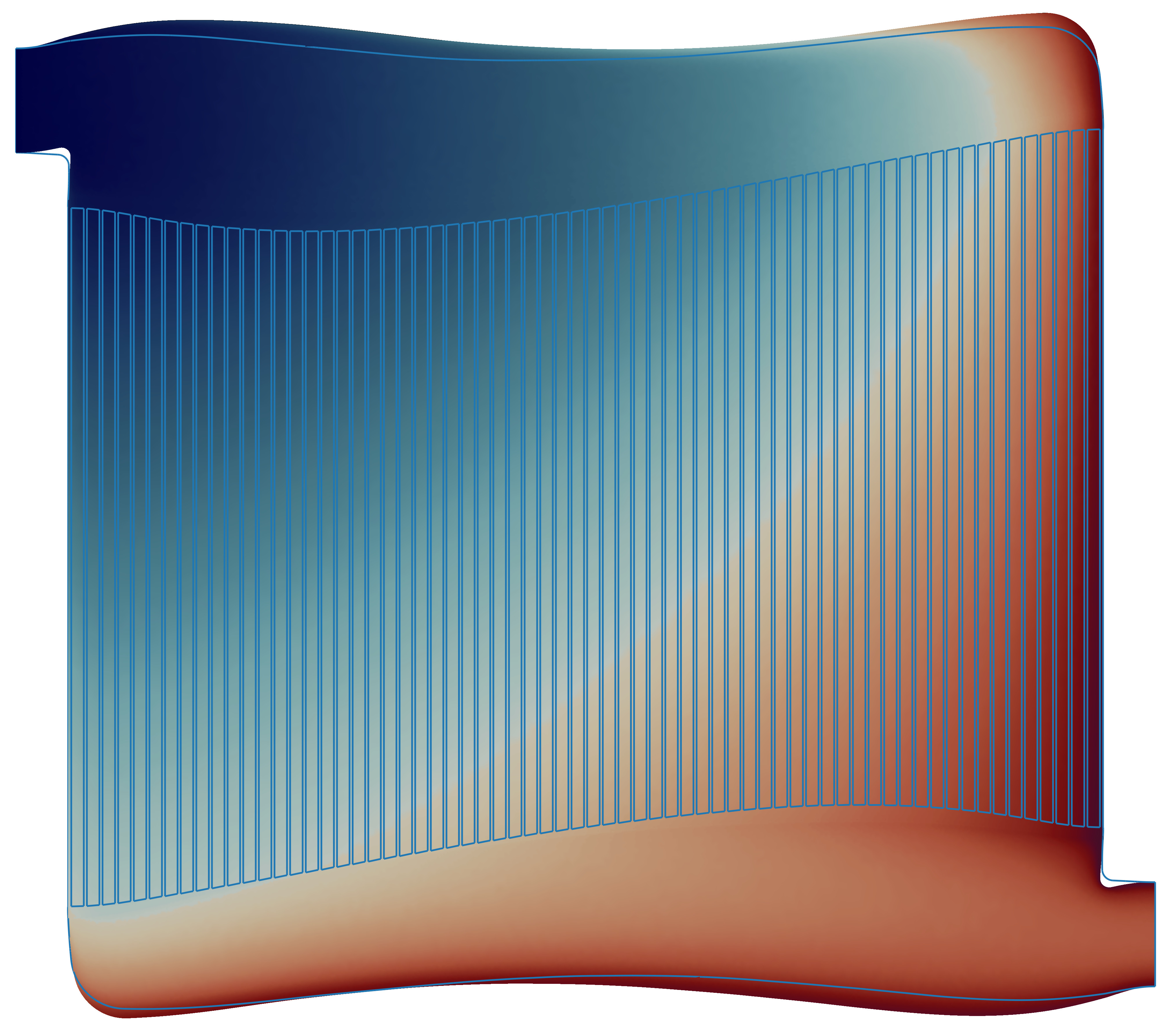}
		\caption{Darcy 3D Model.}
	\end{subfigure}
	\hfil 
	\begin{subfigure}{0.49\textwidth}
		\centering
		\includegraphics[width=0.9\textwidth]{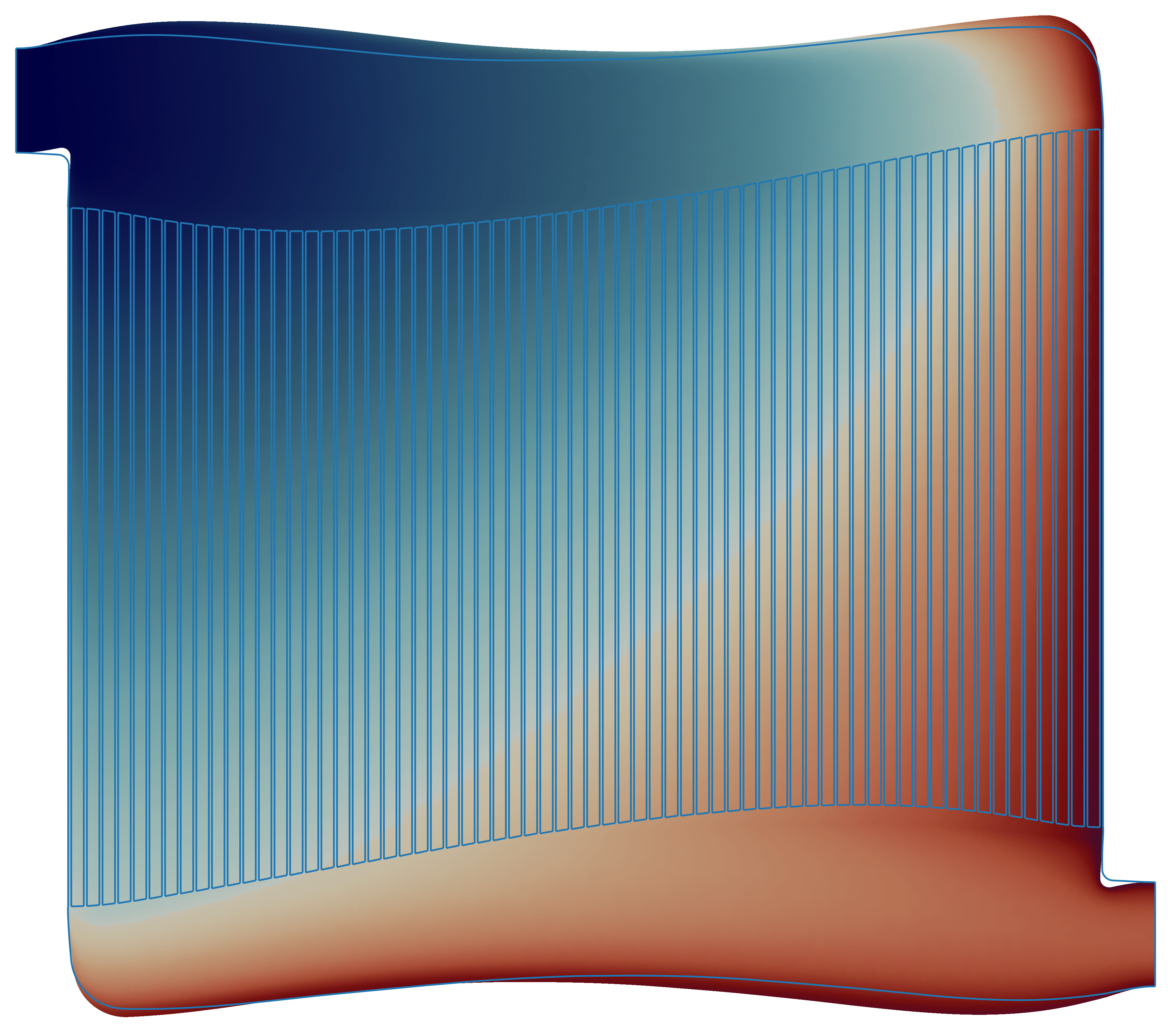}
		\caption{Darcy 2D Model.}
	\end{subfigure}
	\includegraphics[width=\textwidth]{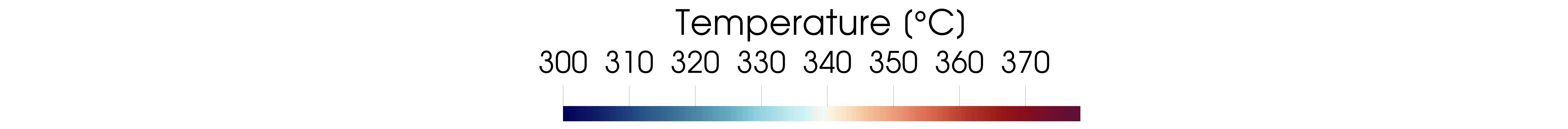}
	\caption{Optimized geometries for all models, depicting the temperature distribution.}
	\label{figure:optimized_geometries}
\end{figure}

The optimized geometries for all four models can be seen in Figure~\ref{figure:optimized_geometries}, where additionally the corresponding temperature distribution is shown. To compare them, we also incorporated the boundaries of the optimized geometry for the full 3D model into the plots for the other models. We observe that all four geometries look very similar. In particular, there are only minor differences between the full 3D and the full 2D model. The full 2D model has some slightly longer or shorter channels as well as slight differences in the outer boundaries of the in- and outlet domains. Compared to the optimized geometry of the full 3D model, both Darcy models show larger deviations. The in- and outlet domains for the Darcy models are consistently further to the outside than their equivalents. Additionally, the curves near the in- and outlet are pushed deeper into the geometry. However, these differences are to be expected due to the changes in modeling. In contrast, the region $\porous{\Omega}_\subdarcy$ behaves similarly to its counterpart, $\Omega_\subchannels$. The curves describing the in- and outlets of the channels are very similar between the full 3D model and the Darcy ones (not visualized). Moreover, comparing the optimized geometries of both Darcy models reveals only very subtle differences. These results suggest that the differences in the optimized geometries for the Darcy models mainly arise from the porous modeling and not from the dimension reduction.

\begin{table}[tb]
	\centering
	\rowcolors{2}{\tablegray}{white}
	\begin{tabular}{r r r r r r r r r}
		\toprule
		& \hspace{1em} & Full 3D & \hspace{1em} & Darcy 3D & \hspace{1em} & Full 2D & \hspace{1em} & Darcy 2D \\
		\midrule
		time (speedup) & & 149416 \si{\second} & & 37040 \si{\second} (4) & & 1440 \si{\second} (104) & & 869 \si{\second} (172) \\
		\bottomrule
	\end{tabular}
	\caption{Comparison of the wall time for the solution of the optimization problem.}
	\label{table:times}
\end{table}

In conclusion, all reduced models also work very well in the shape optimization context. Additionally, since the main numerical work of Algorithm~\ref{algo:descent} consists of solving the state and adjoint systems, the reduced models demand significantly less computational resources. This is depicted in Table~\ref{table:times}, where the time for the solution of the shape optimization problems as well as the speedup for the reduced models relative to the full 3D model is shown. Again, we observe that the two-dimensional models are over 100 times faster than the full 3D models, making them particularly attractive. Finally, we note that, as for the state system, the full 2D model shows the best performance of all reduced models due to its combination of accuracy and efficiency.

\section{Conclusion and Outlook}
\label{sec:conclusion_outlook}

In this work, we introduced four different models for a microchannel cooling system using both porous medium modeling and a dimension reduction technique. A numerical comparison showed that all reduced models approximate the original one quite well while requiring substantially less computational resources. Further, we presented a shape optimization problem based on heat absorbed by the cooler and the uniform distribution of coolant among the microchannels which we adapted to all reduced models. For all our models, we presented both the shape derivative and the adjoint systems we derived with a material derivative free Lagrangian approach, which we rigorously analyzed theoretically in our earlier work \cite{blauth}.
We solved the shape optimization problems numerically using a gradient descent method, where we computed the shape gradient with equations of anisotropic, inhomogeneous, linear elasticity. The numerical results for the shape optimization show that all models behave similarly and that they yield similar optimized geometries.

For future work, the shape optimization of a heat source, e.g., a chemical reactor or an electronic device, coupled to the cooling system can be investigated. This would yield more realistic models and optimization problems. Additionally, the models introduced in this work can be used as building blocks for models describing such a coupled system. Finally, considering a topology optimization before the shape optimization could further increase the quality of the geometries.


\section*{Acknowledgments}

S. Blauth gratefully acknowledges financial support from the Fraunhofer Institute for Industrial Mathematics (ITWM).


\bibliographystyle{siam}
\bibliography{lit.bib}

\begin{thebibliography}{10}

\bibitem{perimeter_regularization}
{\sc G.~{Allaire} and A.~{Henrot}}, {\em {On some recent advances in shape
  optimization.}}, {C. R. Acad. Sci., Paris, S\'er. IIb, M\'ec.}, 329 (2001),
  pp.~383--396.
\newblock \url{https://doi.org/10.1016/S1620-7742(01)01349-6}.

\bibitem{fenics}
{\sc M.~S. Aln{\ae}s, J.~Blechta, J.~Hake, A.~Johansson, B.~Kehlet, A.~Logg,
  C.~Richardson, J.~Ring, M.~E. Rognes, and G.~N. Wells}, {\em The fenics
  project version 1.5}, Archive of Numerical Software, 3 (2015).
\newblock \url{https://doi.org/10.11588/ans.2015.100.20553}.

\bibitem{andreasen}
{\sc C.~S. Andreasen, A.~R. Gersborg, and O.~Sigmund}, {\em Topology
  optimization of microfluidic mixers}, Internat. J. Numer. Methods Fluids, 61
  (2009), pp.~498--513.
\newblock \url{https://doi.org/10.1002/fld.1964}.

\bibitem{petsc-user-ref}
{\sc S.~Balay, S.~Abhyankar, M.~F. Adams, J.~Brown, P.~Brune, K.~Buschelman,
  L.~Dalcin, A.~Dener, V.~Eijkhout, W.~D. Gropp, D.~Karpeyev, D.~Kaushik, M.~G.
  Knepley, D.~A. May, L.~C. McInnes, R.~T. Mills, T.~Munson, K.~Rupp, P.~Sanan,
  B.~F. Smith, S.~Zampini, H.~Zhang, and H.~Zhang}, {\em {PETS}c users manual},
  Tech. Rep. ANL-95/11 - Revision 3.11, Argonne National Laboratory, 2019.
\newblock \url{https://www.mcs.anl.gov/petsc}.

\bibitem{berggren}
{\sc E.~{B\"angtsson}, D.~{Noreland}, and M.~{Berggren}}, {\em {Shape
  optimization of an acoustic horn.}}, {Comput. Methods Appl. Mech. Eng.}, 192
  (2003), pp.~1533--1571.
\newblock \url{https://doi.org/10.1016/S0045-7825(02)00656-4}.

\bibitem{blauth}
{\sc S.~Blauth, C.~Leith\"{a}user, and R.~Pinnau}, {\em Shape sensitivity
  analysis for a microchannel cooling system}, J. Math. Anal. Appl., 492
  (2020), p.~124476.
\newblock \url{https://doi.org/10.1016/j.jmaa.2020.124476}.

\bibitem{dimred}
{\sc T.~Borrvall and J.~Petersson}, {\em Topology optimization of fluids in
  {S}tokes flow}, Internat. J. Numer. Methods Fluids, 41 (2003), pp.~77--107.
\newblock \url{https://doi.org/10.1002/fld.426}.

\bibitem{brinkman}
{\sc H.~C. Brinkman}, {\em A calculation of the viscous force exerted by a
  flowing fluid on a dense swarm of particles}, Flow, Turbulence and
  Combustion, 1 (1949), p.~27.
\newblock \url{https://doi.org/10.1007/BF02120313}.

\bibitem{brooks_hughes}
{\sc A.~N. Brooks and T.~J.~R. Hughes}, {\em Streamline
  upwind/{P}etrov-{G}alerkin formulations for convection dominated flows with
  particular emphasis on the incompressible {N}avier-{S}tokes equations},
  Comput. Methods Appl. Mech. Engrg., 32 (1982), pp.~199--259.
\newblock FENOMECH ''81, Part I (Stuttgart, 1981),
  \url{https://doi.org/10.1016/0045-7825(82)90071-8}.

\bibitem{methanation}
{\sc K.~P. Brooks, J.~Hu, H.~Zhu, and R.~J. Kee}, {\em Methanation of carbon
  dioxide by hydrogen reduction using the sabatier process in microchannel
  reactors}, Chemical Engineering Science, 62 (2007), pp.~1161 -- 1170.
\newblock \url{https://doi.org/10.1016/j.ces.2006.11.020}.

\bibitem{bruus}
{\sc H.~Bruus}, {\em Theoretical microfluidics}, Oxford University Press,
  United States, 2007.

\bibitem{chen2007forced}
{\sc C.-H. Chen}, {\em Forced convection heat transfer in microchannel heat
  sinks}, International Journal of Heat and Mass Transfer, 50 (2007), pp.~2182
  -- 2189.
\newblock \url{https://doi.org/10.1016/j.ijheatmasstransfer.2006.11.001}.

\bibitem{chentopo}
{\sc X.~Chen and T.~Li}, {\em A novel passive micromixer designed by applying
  an optimization algorithm to the zigzag microchannel}, Chemical Engineering
  Journal, 313 (2017), pp.~1406 -- 1414.
\newblock \url{https://doi.org/10.1016/j.cej.2016.11.052}.

\bibitem{delfour_zolesio}
{\sc M.~C. Delfour and J.-P. Zol\'{e}sio}, {\em Shapes and geometries}, vol.~22
  of Advances in Design and Control, Society for Industrial and Applied
  Mathematics (SIAM), Philadelphia, PA, second~ed., 2011.
\newblock Metrics, analysis, differential calculus, and optimization,
  \url{https://doi.org/10.1137/1.9780898719826}.

\bibitem{elman}
{\sc H.~C. Elman, D.~J. Silvester, and A.~J. Wathen}, {\em Finite elements and
  fast iterative solvers: with applications in incompressible fluid dynamics},
  Numerical Mathematics and Scientific Computation, Oxford University Press,
  Oxford, second~ed., 2014.
\newblock \url{https://doi.org/10.1093/acprof:oso/9780199678792.001.0001}.

\bibitem{ern_guermond}
{\sc A.~Ern and J.-L. Guermond}, {\em Theory and practice of finite elements},
  vol.~159 of Applied Mathematical Sciences, Springer-Verlag, New York, 2004.
\newblock \url{https://doi.org/10.1007/978-1-4757-4355-5}.

\bibitem{herzog}
{\sc T.~Etling, R.~Herzog, E.~Loayza, and G.~Wachsmuth}, {\em First and second
  order shape optimization based on restricted mesh deformations}, SIAM J. Sci.
  Comput., 42 (2020), pp.~A1200--A1225.
\newblock \url{https://doi.org/10.1137/19M1241465}.

\bibitem{ferraris}
{\sc A.~Ferraris, A.~G. Airale, D.~Berti~Polato, A.~Messana, S.~Xu, P.~Massai,
  and M.~Carello}, {\em City car drag reduction by means of shape optimization
  and add-on devices}, in Advances in Mechanism and Machine Science, T.~Uhl,
  ed., Cham, 2019, Springer International Publishing, pp.~3721--3730.
\newblock \url{https://doi.org/10.1007/978-3-030-20131-9_367}.

\bibitem{foli2006optimization}
{\sc K.~Foli, T.~Okabe, M.~Olhofer, Y.~Jin, and B.~Sendhoff}, {\em Optimization
  of micro heat exchanger: Cfd, analytical approach and multi-objective
  evolutionary algorithms}, International Journal of Heat and Mass Transfer, 49
  (2006), pp.~1090 -- 1099.
\newblock \url{https://doi.org/10.1016/j.ijheatmasstransfer.2005.08.032}.

\bibitem{gangl_topo}
{\sc P.~Gangl and U.~Langer}, {\em Topology optimization of electric machines
  based on topological sensitivity analysis}, Comput. Vis. Sci., 15 (2012),
  pp.~345--354.
\newblock \url{https://doi.org/10.1007/s00791-014-0219-6}.

\bibitem{gangl_shape}
{\sc P.~Gangl, U.~Langer, A.~Laurain, H.~Meftahi, and K.~Sturm}, {\em Shape
  optimization of an electric motor subject to nonlinear magnetostatics}, SIAM
  J. Sci. Comput., 37 (2015), pp.~B1002--B1025.
\newblock \url{https://doi.org/10.1137/15100477X}.

\bibitem{sigmund}
{\sc A.~Gersborg-Hansen, O.~Sigmund, and R.~B. Haber}, {\em Topology
  optimization of channel flow problems}, Struct. Multidiscip. Optim., 30
  (2005), pp.~181--192.
\newblock \url{https://doi.org/10.1007/s00158-004-0508-7}.

\bibitem{gmsh}
{\sc C.~Geuzaine and J.-F. Remacle}, {\em Gmsh: {A} 3-{D} finite element mesh
  generator with built-in pre- and post-processing facilities}, Internat. J.
  Numer. Methods Engrg., 79 (2009), pp.~1309--1331.
\newblock \url{https://doi.org/10.1002/nme.2579}.

\bibitem{review_channels}
{\sc J.~S. Goodling}, {\em {Microchannel heat exchangers: a review}}, in High
  Heat Flux Engineering II, A.~M. Khounsary, ed., vol.~1997, International
  Society for Optics and Photonics, SPIE, 1993, pp.~66 -- 82.
\newblock \url{https://doi.org/10.1117/12.163830}.

\bibitem{pipe2}
{\sc E.~Helgason and S.~Krajnovi{\'c}}, {\em Aerodynamic shape optimization of
  a pipe using the adjoint method}, in ASME International Mechanical
  Engineering Congress and Exposition, vol.~45233, American Society of
  Mechanical Engineers, 2012, pp.~259--267.
\newblock \url{https://doi.org/10.1115/IMECE2012-89396}.

\bibitem{hinze_pinnau_ulbrich2}
{\sc M.~Hinze, R.~Pinnau, M.~Ulbrich, and S.~Ulbrich}, {\em Optimization with
  {PDE} constraints}, vol.~23 of Mathematical Modelling: Theory and
  Applications, Springer, New York, 2009.
\newblock \url{https://doi.org/10.1007/978-1-4020-8839-1}.

\bibitem{hohmann}
{\sc R.~Hohmann and C.~Leith\"{a}user}, {\em Shape {O}ptimization of a
  {P}olymer {D}istributor {U}sing an {E}ulerian {R}esidence {T}ime {M}odel},
  SIAM J. Sci. Comput., 41 (2019), pp.~B625--B648.
\newblock \url{https://doi.org/10.1137/18M1225847}.

\bibitem{hohmann_erosion}
{\sc R.~Hohmann and C.~Leithäuser}, {\em Gradient-based shape optimization for
  the reduction of particle erosion in bended pipes}, 2019.
\newblock \url{https://arxiv.org/abs/1908.04712}.

\bibitem{john}
{\sc V.~John}, {\em Finite element methods for incompressible flow problems},
  vol.~51 of Springer Series in Computational Mathematics, Springer, Cham,
  2016.
\newblock \url{https://doi.org/10.1007/978-3-319-45750-5}.

\bibitem{kaviany}
{\sc M.~Kaviany}, {\em Principles of Heat Transfer in Porous Media}, Springer
  US, New York, NY, 1991.
\newblock \url{https://doi.org/10.1007/978-1-4684-0412-8}.

\bibitem{kelley}
{\sc C.~T. Kelley}, {\em Iterative methods for optimization}, vol.~18 of
  Frontiers in Applied Mathematics, Society for Industrial and Applied
  Mathematics (SIAM), Philadelphia, PA, 1999.
\newblock \url{https://doi.org/10.1137/1.9781611970920}.

\bibitem{khan2011review}
{\sc M.~G. Khan and A.~Fartaj}, {\em A review on microchannel heat exchangers
  and potential applications}, International Journal of Energy Research, 35
  (2011), pp.~553--582.
\newblock \url{https://doi.org/10.1002/er.1720}.

\bibitem{kim2000local}
{\sc S.~Kim, D.~Kim, and D.~Lee}, {\em On the local thermal equilibrium in
  microchannel heat sinks}, International Journal of Heat and Mass Transfer, 43
  (2000), pp.~1735 -- 1748.
\newblock \url{https://doi.org/10.1016/S0017-9310(99)00259-8}.

\bibitem{kim1999forced}
{\sc S.~J. Kim and D.~Kim}, {\em {Forced Convection in Microstructures for
  Electronic Equipment Cooling}}, Journal of Heat Transfer, 121 (1999),
  pp.~639--645.
\newblock \url{https://doi.org/10.1115/1.2826027}.

\bibitem{etching}
{\sc G.~Kolb, T.~Baier, J.~Schürer, D.~Tiemann, A.~Ziogas, H.~Ehwald, and
  P.~Alphonse}, {\em A micro-structured 5kw complete fuel processor for
  iso-octane as hydrogen supply system for mobile auxiliary power units: Part
  i. development of autothermal reforming catalyst and reactor}, Chemical
  Engineering Journal, 137 (2008), pp.~653 -- 663.
\newblock \url{https://doi.org/10.1016/j.cej.2007.06.033}.

\bibitem{kubo_topology}
{\sc S.~Kubo, K.~Yaji, T.~Yamada, K.~Izui, and S.~Nishiwaki}, {\em A level
  set-based topology optimization method for optimal manifold designs with flow
  uniformity in plate-type microchannel reactors}, Struct. Multidiscip. Optim.,
  55 (2017), pp.~1311--1327.
\newblock \url{https://doi.org/10.1007/s00158-016-1577-0}.

\bibitem{leith3}
{\sc C.~Leith\"{a}user and R.~Pinnau}, {\em The production of filaments and
  non-woven materials: the design of the polymer distributor}, in Math for the
  digital factory, vol.~27 of Math. Ind., Springer, Cham, 2017, pp.~321--340.
\newblock \url{https://doi.org/10.1007/978-3-319-63957-4_15}.

\bibitem{leith2}
{\sc C.~Leith\"{a}user, R.~Pinnau, and R.~Fe\ss~ler}, {\em Designing polymer
  spin packs by tailored shape optimization techniques}, Optim. Eng., 19
  (2018), pp.~733--764.
\newblock \url{https://doi.org/10.1007/s11081-018-9396-3}.

\bibitem{leith}
{\sc C.~Leith\"{a}user, R.~Pinnau, and R.~Fe{\ss}ler}, {\em Shape design for
  polymer spin packs: modeling, optimization and validation}, J. Math. Ind., 8
  (2018), pp.~Paper No. 13, 17.
\newblock \url{https://doi.org/10.1186/s13362-018-0055-2}.

\bibitem{fenics_book}
{\sc A.~Logg, K.-A. Mardal, G.~N. Wells, et~al.}, {\em Automated Solution of
  Differential Equations by the Finite Element Method}, Springer, 2012.
\newblock \url{https://doi.org/10.1007/978-3-642-23099-8}.

\bibitem{lyu}
{\sc Z.~Lyu and J.~R. R.~A. Martins}, {\em Aerodynamic design optimization
  studies of a blended-wing-body aircraft}, Journal of Aircraft, 51 (2014),
  pp.~1604--1617.
\newblock \url{https://doi.org/10.2514/1.C032491}.

\bibitem{naqiuddin2018overview}
{\sc N.~H. Naqiuddin, L.~H. Saw, M.~C. Yew, F.~Yusof, T.~C. Ng, and M.~K. Yew},
  {\em Overview of micro-channel design for high heat flux application},
  Renewable and Sustainable Energy Reviews, 82 (2018), pp.~901 -- 914.
\newblock \url{https://doi.org/10.1016/j.rser.2017.09.110}.

\bibitem{nield_bejan}
{\sc D.~A. Nield and A.~Bejan}, {\em Convection in porous media}, Springer,
  Cham, 2017.
\newblock Fifth edition, \url{https://doi.org/10.1007/978-3-319-49562-0}.

\bibitem{nocedal_wright}
{\sc J.~Nocedal and S.~J. Wright}, {\em Numerical optimization}, Springer
  Series in Operations Research, Springer-Verlag, New York, 1999.
\newblock \url{https://doi.org/10.1007/b98874}.

\bibitem{novotny_sokolowski}
{\sc A.~A. Novotny and J.~Soko{\l}owski}, {\em Topological derivatives in shape
  optimization}, Interaction of Mechanics and Mathematics, Springer,
  Heidelberg, 2013.
\newblock \url{https://doi.org/10.1007/978-3-642-35245-4}.

\bibitem{othmer}
{\sc C.~Othmer}, {\em Adjoint methods for car aerodynamics}, J. Math. Ind., 4
  (2014), pp.~Art. 6, 23.
\newblock \url{https://doi.org/10.1186/2190-5983-4-6}.

\bibitem{pan2008optimal}
{\sc M.~Pan, Y.~Tang, L.~Pan, and L.~Lu}, {\em Optimal design of complex
  manifold geometries for uniform flow distribution between microchannels},
  Chemical Engineering Journal, 137 (2008), pp.~339 -- 346.
\newblock \url{https://doi.org/10.1016/j.cej.2007.05.012}.

\bibitem{itakura}
{\sc E.~M. Papoutsis-Kiachagias, V.~G. Asouti, K.~C. Giannakoglou, K.~Gkagkas,
  S.~Shimokawa, and E.~Itakura}, {\em Multi-point aerodynamic shape
  optimization of cars based on continuous adjoint}, Struct. Multidiscip.
  Optim., 59 (2019), pp.~675--694.
\newblock \url{https://doi.org/10.1007/s00158-018-2091-3}.

\bibitem{ihtc1}
{\sc D.~L. Penha, S.~Stolz, J.~Kuerten, M.~Nordlund, A.~Kuczaj, and B.~Geurts},
  {\em Fully-developed conjugate heat transfer in porous media with uniform
  heating}, International Journal of Heat and Fluid Flow, 38 (2012), pp.~94 --
  106.
\newblock \url{https://doi.org/10.1016/j.ijheatfluidflow.2012.08.007}.

\bibitem{printing}
{\sc A.~Radwan, S.~Ookawara, S.~Mori, and M.~Ahmed}, {\em Uniform cooling for
  concentrator photovoltaic cells and electronic chips by forced convective
  boiling in 3d-printed monolithic double-layer microchannel heat sink}, Energy
  Conversion and Management, 166 (2018), pp.~356 -- 371.
\newblock \url{https://doi.org/10.1016/j.enconman.2018.04.037}.

\bibitem{rao2016dimensional}
{\sc R.~Rao, K.~More, J.~Taler, and P.~Ocłoń}, {\em Dimensional optimization
  of a micro-channel heat sink using jaya algorithm}, Applied Thermal
  Engineering, 103 (2016), pp.~572 -- 582.
\newblock \url{https://doi.org/10.1016/j.applthermaleng.2016.04.135}.

\bibitem{freecad}
{\sc J.~Riegel, W.~Mayer, and Y.~van Havre}, {\em Freecad (version 0.16)}.
\newblock \url{https://www.freecadweb.org/}.

\bibitem{roos}
{\sc H.-G. Roos, M.~Stynes, and L.~Tobiska}, {\em Numerical methods for
  singularly perturbed differential equations}, vol.~24 of Springer Series in
  Computational Mathematics, Springer-Verlag, Berlin, 1996.
\newblock Convection-diffusion and flow problems,
  \url{https:doi.org/10.1007/978-3-662-03206-0}.

\bibitem{salimpour2013constructal}
{\sc M.~R. Salimpour, M.~Sharifhasan, and E.~Shirani}, {\em Constructal
  optimization of microchannel heat sinks with noncircular cross sections},
  Heat Transfer Engineering, 34 (2013), pp.~863--874.
\newblock \url{https://doi.org/10.1080/01457632.2012.746552}.

\bibitem{gauger}
{\sc S.~Schmidt, C.~Ilic, V.~Schulz, and N.~R. Gauger}, {\em Airfoil design for
  compressible inviscid flow based on shape calculus}, Optim. Eng., 12 (2011),
  pp.~349--369.
\newblock \url{https://doi.org/10.1007/s11081-011-9145-3}.

\bibitem{schmidt2013three}
\leavevmode\vrule height 2pt depth -1.6pt width 23pt, {\em Three-dimensional
  large-scale aerodynamic shape optimization based on shape calculus}, AIAA
  Journal, 51 (2013), pp.~2615--2627.
\newblock \url{https://doi.org/10.2514/1.J052245}.

\bibitem{schmidt_horn}
{\sc S.~Schmidt, E.~Wadbro, and M.~Berggren}, {\em Large-scale
  three-dimensional acoustic horn optimization}, SIAM J. Sci. Comput., 38
  (2016), pp.~B917--B940.
\newblock \url{https://doi.org/10.1137/15M1021131}.

\bibitem{pipe1}
{\sc A.~Schulz}, {\em The optimal shape of a pipe}, Z. Angew. Math. Phys., 64
  (2013), pp.~1177--1185.
\newblock \url{https://doi.org/10.1007/s00033-012-0277-x}.

\bibitem{welker}
{\sc V.~H. Schulz, M.~Siebenborn, and K.~Welker}, {\em Efficient {PDE}
  constrained shape optimization based on {S}teklov-{P}oincar\'{e}-type
  metrics}, SIAM J. Optim., 26 (2016), pp.~2800--2819.
\newblock \url{https://doi.org/10.1137/15M1029369}.

\bibitem{soghrati2013computational}
{\sc S.~Soghrati, A.~R. Najafi, J.~H. Lin, K.~M. Hughes, S.~R. White, N.~R.
  Sottos, and P.~H. Geubelle}, {\em Computational analysis of actively-cooled
  3d woven microvascular composites using a stabilized interface-enriched
  generalized finite element method}, International Journal of Heat and Mass
  Transfer, 65 (2013), pp.~153 -- 164.
\newblock \url{https://doi.org/10.1016/j.ijheatmasstransfer.2013.05.054}.

\bibitem{sokolowski_zolesio}
{\sc J.~Soko{\l}owski and J.-P. Zol\'{e}sio}, {\em Introduction to shape
  optimization}, vol.~16 of Springer Series in Computational Mathematics,
  Springer-Verlag, Berlin, 1992.
\newblock Shape sensitivity analysis,
  \url{https://doi.org/10.1007/978-3-642-58106-9}.

\bibitem{sturm}
{\sc K.~Sturm}, {\em Minimax {L}agrangian approach to the differentiability of
  nonlinear {PDE} constrained shape functions without saddle point assumption},
  SIAM J. Control Optim., 53 (2015), pp.~2017--2039.
\newblock \url{https://doi.org/10.1137/130930807}.

\bibitem{sturm_new_trends}
\leavevmode\vrule height 2pt depth -1.6pt width 23pt, {\em Shape
  differentiability under non-linear {PDE} constraints}, in New trends in shape
  optimization, vol.~166 of Internat. Ser. Numer. Math.,
  Birkh\"{a}user/Springer, Cham, 2015, pp.~271--300.
\newblock \url{https://doi.org/10.1007/978-3-319-17563-8_12}.

\bibitem{ihtc2}
{\sc F.~E. Teruel and L.~Díaz}, {\em Calculation of the interfacial heat
  transfer coefficient in porous media employing numerical simulations},
  International Journal of Heat and Mass Transfer, 60 (2013), pp.~406 -- 412.
\newblock \url{https://doi.org/10.1016/j.ijheatmasstransfer.2012.12.022}.

\bibitem{troeltzsch}
{\sc F.~Tr\"{o}ltzsch}, {\em Optimal control of partial differential
  equations}, vol.~112 of Graduate Studies in Mathematics, American
  Mathematical Society, Providence, RI, 2010.
\newblock Theory, methods and
  applications,\url{https://doi.org/10.1090/gsm/112}.

\bibitem{yan2019topology}
{\sc S.~Yan, F.~Wang, J.~Hong, and O.~Sigmund}, {\em Topology optimization of
  microchannel heat sinks using a two-layer model}, International Journal of
  Heat and Mass Transfer, 143 (2019), p.~118462.
\newblock \url{https://doi.org/10.1016/j.ijheatmasstransfer.2019.118462}.

\end{thebibliography}

\end{document}